\documentclass[twoside,12pt,reqno]{amsart}
\usepackage{amsfonts, amsbsy, amsmath, amsthm, amssymb, latexsym, verbatim, enumerate}

\usepackage{comment}
\usepackage[english]{babel}

\newtheorem{propo}{Proposition}[section]
\newtheorem{defi}[propo]{Definition}

\newtheorem{lemma}[propo]{Lemma}
\newtheorem{corol}[propo]{Corollary}

\newtheorem{theo}[propo]{Theorem}

\newcommand{\ld}{,\ldots ,}

\newcommand{\ra}{ \rightarrow }

\newcommand{\lan}{ \langle }
\newcommand{\ran}{ \rangle }

\newcommand{\diag}{\mathop{\rm diag}\nolimits}

\newcommand{\Id}{\mathop{\rm Id}\nolimits}

\newcommand{\Om}{\Omega}

\newcommand{\CC}{\mathbb{C}}
\newcommand{\FF}{\mathbb{F}}
\newcommand{\NN}{\mathbb{N}}
\newcommand{\QQ}{\mathbb{Q}}
\newcommand{\RR}{\mathbb{R}}
\newcommand{\ZZ}{\mathbb{Z}}

\newcommand{\al}{\alpha}
\newcommand{\be}{\beta}

\newcommand{\ep}{\varepsilon}

\newcommand{\lam}{\lambda }

\newcommand{\up}{^{-1}}

\newcommand{\om}{\omega }
\newcommand{\si}{\sigma }
\newcommand{\med}{\medskip }

\newcommand{\GG}{{\mathbf G}}
\newcommand{\mar}{\marginpar}

\newcommand{\bl}{\begin{lemma}\label}

\newcommand{\bp}{\begin{proof}}
\newcommand{\enp}{\end{proof}}

\newcommand{\AGL}{\operatorname{AGL}}

\def\d12{{_{12}}}

\def\ac{{algebraically closed }}
\def\acf{{algebraically closed field }}

\def\au{{automorphism }}

\def\ccc{{constituent }}

\def\sg{subgroup }

\def\ei{{eigenvalue }}
\def\eis{{eigenvalues }}

\def\f{{following }}

\def\ho{{homomorphism }}

\def\ii{{if and only if }}

\def\ir{{irreducible }}

\def\irt{{irreducible. }}
\def\irr{{irreducible representation }}

\def\itf{{It follows that }}

\def\mult{{multiplicity }}

\def\hw{{highest weight }}

\def\po{{polynomial }}
\def\pos{{polynomials }}

\def\rep{{representation }}
\def\reps{{representations }}
\def\rept{{representation. }}

\def\sg{{subgroup }}

\def\syl{{Sylow $p$-subgroup }}

\def\SU{{\rm SU}}
\def\GL{{\rm GL}}
\def\Sp{{\rm Sp}}

\def\F{{\rm F}}
\def\Sym{{\rm Sym}}

\def\Out{{\rm Out}}

\newcommand{\el}{\end{lemma}}

\begin{document}
\title  [Eigenvalue 1]{Unisingular  subgroups of symplectic group\\ $Sp_{2n}(2)$ for $2n< 250$}

\author[A. Zalesski]{Alexandre Zalesski}

\address{Department of  Mathematics, Department of Mathematics, University Brasilia, Brasilia-DF 70910-900, Brazil}
\email{alexandre.zalesski@gmail.com}
\maketitle

Abstract. A linear group is called {\it unisingular} if every element of it has eigenvalue 1.
A certain aspect of the theory of abelian varieties requires the knowledge
of   unisingular irreducible subgroups of the symplectic groups over the field of two elements. A more special, but an important question is on the existence of such subgroups in the symplectic groups of particular degree. We answer this question
for almost all degrees $2n<250$, specifically, the question remains open only 7 values of $n$. Additionally, the paper contains  results of general nature on the structure of unisingular  irreducible  linear groups.\footnote{Mathematics subject classification: 11G10, 11F80, 20C33, 20H30}\,\footnote{Keywords: eigenvalue 1, symplectic group, irreducible subgroups, finite linear groups, finite group representations, abelian varieties, torsion points}


\section{Introduction}

\def\ch{characteristic }

Let $G$ be a finite group and $\phi$  a \rep of $G$ over a field $F$. 
We say that $\phi$ is {\it unisingular} if $\phi(g)$ has \ei 1 for every $g\in G$. Unisingular \ir \reps
were first considered in \cite{Z90} where the Steinberg \reps of most simple groups of Lie type
were shown to be unisingular (for $F$ to be of defining characteristic), see also \cite{Z88}.
The term "unisingular" was introduced in \cite{GT}, where the authors focused on the
classification of "unisingular" groups, that is, those with all \reps unisingular.

The question of the existence of \ei 1 for a particular element of a linear group
is essential for applications due to its geometric nature: this means that the element in question fixes a non-zero vector at the underlying vector space. Attempts of systematic study of the \ei 1 occurance were made in  \cite{Z87} and \cite{Z90} for groups of Lie type, see also \cite[Problem $1'$, p. 207]{Z88}.   In \cite{GT} the authors  look for a bound for the number of fixed point free elements in certain linear groups. This kind of questions also arises in a more uniform way as follows:

\med
Problem 1. Determine finite \ir linear groups whose every element has \ei 1.

\med
This problem is explicitly stated in \cite{gl3} and, for \reps of simple groups of Lie type in
their defining characteristic, in  \cite{z18} and \cite{z21}. Note that Problem 1 cannot have any explicit answer in full generality but it describes an area which various more special problems  belong to. A less universal version of Problem 1 is 

\med
Problem 2. Given a prime $p$ determine finite \ir linear groups whose every $p$-element has \ei 1.

\med
For quasi-simple finite linear groups over the complex numbers and elements of  prime  order   a full solution to Problem 2 is obtained in \cite{z08}. A more precise version of Problem  2
requires, given an element $g$ of prime power order in a group $G$, determine the \ir \reps $\phi$ of $G$ over a field $F$ such that 1 is not an \ei of $\phi(g)$.  
See \cite{TZ8} and \cite{TZ22} for the case of cross-characeristic \reps of groups of Lie type.

 Unisingular \reps over fields of characteristic 2 are of particular interest for a certain aspect of the study of abelian varieties in algebraic geometry, as explained in detail in \cite{cu10} and \cite{cu12},  see also   \cite[Corollary 1.9]{CZ}; this is  based on the fundamental work by Katz \cite{katz}. 
The \f problem  forms a general frame of research activity in this direction:

\med
Problem 3. Determine \ir unisingular subgroups of $Sp_{2n}(2)$, that is, the groups  whose every  element has \ei $1$.

\med
The cases of $n=1 ,2$ were settled in \cite{katz} and \cite{cu9}, for $n=3$ see \cite[Lemma 3]{cu9}.
At the initial stage of the study of Problem 3 it was probably expected that 
unisingular subgroups of $Sp_{2n}(2)$ are very rare, and the first example was suggested by J.-P. Serre for $n=4$, see a discussion in \cite[page 1833]{cu9}. In \cite{CZ} the authors  produce a full list of maximal unisingular \ir subgroups of $Sp_{8}(2)$. In fact, Serre's example is a special case of the \f theorem:

\begin{theo}\label{se2} Let $G$ be a finite simple group of Lie type in defining characteristic $2$. Suppose that $G$ is not isomorphic to $PSL_2(q)$ for $q$ even. Then $G$ is isomorphic to a unisingular \ir subgroup of $Sp_{2n}(2)$ with $2n=|G|_2$, where $|G|_2$ is the $2$-part of $|G|$. \end{theo}
 
Special cases of this result appeared in \cite{Z90} and complemented in \cite[Theorem 1.8]{CZ}. Moreover, in \cite[Theorem 1.10]{CZ}  are obtained sufficient conditions of unisingularity for   $2$-modular  \ir \reps  of simple symplectic and orthogonal groups over $\mathbb{\FF}_q$, $q$ even. For groups $GL_k(2)$
and $Sp_{2k}(2)$ all unisingular $2$-modular  \ir \reps have been determined in \cite{z18} and \cite{z21}, respectively. The abundance of examples makes it evident that Problem 3 cannot have an explicit solution for arbitrary $n$.  One can narrow the content   of Problem 3 by trying to determine only maximal unisingular  subgroups of $Sp_{2n}(2)$. However, even this version is too ambitious as obtaining a solution requires too much cases-by-case analysis. Note that the situation cannot allow any uniform treatment in terms of $n$. Moreover, a feature of Problem 3 (as well as of Problem 1) is that the number of subgroups in question for a given $n$ crutially depends 
on the  factorisation of $n$ as a product of primes. Our experience leads us to single out 
the following extremal special case   of Problem 3:
 
 \med
 Problem 4. Determine the set  $\NN_0$ of  integers $n\geq 1$ such that 
$Sp_{2n}(2)$ contains no unisingular \ir subgroup.

\med
Note that a similar problem is easy for subgroups of $GL_n(r)$ with $r$ odd, see Theorem \ref{fo8}.

\med
Conjecture 1. $\NN_0  $ is   infinite.
\med

   In this paper we
 examine Conjecture 1 for $n< 250$ and arrive at the following conclusion:

  \begin{theo}\label{th1} Let $n<125$. 

 $(1)$ $n\in \NN_0$ for $n\in\{1,2,3,5,9,27,29,43,53,106,113\}$. 

 $(2)$  $n\notin \NN_0$ unless $n$ appears in $(1)$ and possibly $n\in\{47,58,67,83,86,103,107\}$. 
 \end{theo}

 Thus, $\NN_0$ contains at least 11 and at most 18 values in the range $1\leq n<125$. 
The cases with $n\in\{47,58,67, 83,86,103,107\}$ remain open. Some information on 
unisingular subgroups of $Sp_{2n}(2)$ is accumulated in Table 5
at the end of this paper.  

One easily observes that if $n\notin \NN_0$ then  $nk\notin \NN_0 $ for every integer $k>1$. This fact shows that the cases of $n$ prime are of particular interest. So we state
the \f special case of Conjecture 1:

\med
Conjecture 2. $\NN_0$ contains   infinitely many primes.  

\med
For $n$ prime we have the following alternative:

\bl{2po3}  Let $p>2$ be a prime and let  $G$ be a unisingular \ir subgroup of $Sp_{2p}(2)$ and let N be a minimal normal subgroup of G. 

$(1)$ N is either a elementary abelian normal selfcentralizing  $3$-group or a simple non-abelian group.

$(2)$ If N is abelian then G contains an \ir elementary abelian-by-cyclic subgroup X of order $3^kp$ for $p>3$ and some integer $k>2$. If N is simple then N is isomorphic to an absolutely \ir subgroup of $Sp_{2p}(2)$ or $U_p(4)$.\el

Lemmas \ref{2po3} reduces the study of Conjecture 2 (and of Problem 3  for $n$ a prime) to groups $G$ that either non-abelian simple or have an elementary abelian normal subgroup $A$, say, with $C_G(A)=A$. They are sometimes called groups of affine type (usually one assumes that $G$ splits over $A$).  

In general, there is no efficient tool for deciding, given a natural number $n$, whether any simple group $N$ has a unisingular 2-modular \irr of degree $n$ and $2n$. A more promissing approach is, given a series of simple group $G$, to determine unisingular
absolutely \ir \reps in characteristic 2, and then decide which of them is unisingular and realizes over $\FF_2$.  In Section  4  we do this for the series $PSL_2(q)$. In particular, we observe that there infinitely many integers $n$ such that $Sp_{2n}(2)$ contains an
absolutely \ir subgroup isomorphic to some   $PSL_2(q)$. Specifying $n$ to be a prime here leads to a question, probably non-trivial, whether $\NN\setminus \NN_0$ contains  infinitely many primes. 
  
 Hiss and Malle \cite{HM,HM1} provide a full list of simple groups that have an absolutely \irr of degree at most 250. This result with some additional efforts reduces Conjecture 2 for $n<125$ to groups of affine types.  

For distinct primes $p,r$ denote by $G_{r,p}$ a non-abelian group of minimal order $r^dp$. Such a group has a unique non-trivial normal   $r$-subgroup $A$ of  index $p$. Obviously, $A$ is an elementary abelian $r$-group of rank $d$.  Then the conjugation action of $G$ on   $A$ is irreducible in the sense that there is no non-trivial proper $G$-invariant subgroup of $A$. One observes that there is only one such a group for given primes $p,r$.  
For such groups we state

\med
Conjecture 3. Let $r>2$. Then there are infinitely many primes $p$ such that every non-trivial \irr of   $G=G_{r,p}$   is not unisingular.   

\med
We have the following result toward Conjecture 3:

\bl{mi2}  Let $p,r>2$ be distinct primes  such that $(r^i-1,p)=1$   for every $0<i<p-1$.
 Let $G=G_{r,p}$.  Then G has no non-trivial unisingular \ir \rep over  any field of characteristic $\ell\neq r$.\el

(Lemma \ref{mi2} can be deduced from \cite[Theorem 6.5]{AGJ}, we provide an alternative proof.) Note that the case with $r=2$ is not of interest as every \irr of $G_{2,p}$ is unisingular. 

Observe that in case (2) of Lemma \ref{2po3} $G$ contains an \ir subgroup isomorphic to $G_{3,p}$ (Lemma \ref{2po}).

Conjecture 3 follows from Lemma \ref{mi2} and Artin's famous conjecture: 


\med\noindent
Artin's Conjecture (1927). For every prime $r$ there are infinitely many primes $p$ such that  $p|(r^m-1)$, $m<p$, implies $m= p-1$. 

\med
Artin's Conjecture  is still open. It is proved in \cite{HB} that there are at most two primes $r$ for which Artin's conjecture fails. However, the method of \cite{HB}  does not allow  to prove the conjecture  for any fixed single prime $r$. 

Lemma \ref{mi2} is not always true if $p|(r^i-1)$ with $1<i<p-1$. We show that it fails 
if $r=3$ and $p\in\{11,23,41,73\}$, see Lemma \ref{46v}. This is based on computer computations performed by Eamonn O'Brian. In general, the question of unisingularity
for \ir \reps of groups $G_{r,p}$ is quite challenging.

Observe that this question is  equivalent to a problem of permutation group theory. 
Let $G$ be a finite group and let $\Om$ be a permutational $G$-set. We say that  
$\Om$ is $p$-{\it restricted} if for every elementary abelian $p$-subgroup $A$ of $G$ all $A$-orbits are of size at most $p$.
  
\bl{eq9} Let  $p,r$ be distinct primes, F a field of characteristic $\ell\neq r$ and let
   $G=G_{r,p}$. Then the \f are equivalent:

$(1)$ there exists a faithful unisingular \irr of G;

$(2)$   every r-element of G fixes a point on some faithful $p$-restricted permutational G-set $\Delta$.\el

Observe that the condition of $\Delta$ to be $p$-restricted is essential. (Indeed, one can 
define a structure of a $G$-set on the regular $A$-set.)

Lemma \ref{mi2} yields a necessary condition for (1) and hence for (2) to be true, but 
it is not known when that condition is sufficent. 

Statement (2) of Lemma \ref{eq9} is related with a general problem of permutation group theory:

\med
Problem 5. Given  a natural number $n>1$ and a prime $p|n$ determine   transitive subgroups of $S_n$ whose every $p$-element fixes a point. 

\med
This problem is a closed analog of Problem 2.  They are not equivalent but Lemma \ref{eq9} shows that there are some common roots of these problems. According with our knowledge, Problem 5 were not discussed in the literature in full generality. 
 
\med
Problem 6. For every fixed prime $p$ determine the set  $\mathbb{P}_p$ of integers $n>1$ such that every transitive permutation group of degree $n$ contains a fixed point free $p$-element. 

\med
Some contribution to this problem is obtained in \cite{Be}. 
In particular, it is shown that $p^a(p+1)\in \mathbb{P}_p$ for $a>1$.  
Further results can be found in \cite{sp} and \cite{Ca}. In fact, essential efforts in these works are directed to solving Isbell's problem:

\med
Problem 7.  Let $G$ be a group of order $p^kb$, where $p$ a fixed prime and $b>1$ is a fixed integer. Is it true that the degree of a transitive permutation \rep of $G$  whose every $p$-elements fixes a point is bounded in terms of $ b$?

\med
This problem was stated by Isbell for $p=2$ and for arbitrary $p$ in  \cite{cfk}. The case of $p=2$ is motivated by the \f result:

\begin{theo}\label{isb} {\rm (Isbell \cite{Is})} Let $n>0$ be even. Then there exists an $n$-player homogeneous  game \ii there exists a transitive permutation group of degree $n$ that  contains no fixed point free $2$-element (in other words, every $2$-element fixes a point).\end{theo}

In fact, there is a bijective correspondence between $n$-player homogeneous  games
and transitive permutation groups of degree $n$ that  contain no fixed point free $2$-element. This fact leads to Problem 5 for $p=2$ as it is equivalent to the classification of $n$-player homogeneous  games.  See \cite[p. 243]{Ca} for more comments.

\med
Problem 5 is very difficult, and possibly is not treatable in full generality. 
It is natural to single out the \f special case:

\med
Problem 8. Given a prime $p$ determine integers  $m>1$ such that 
every transitive subgroup $G\subset S_{pm}$ has fixed point free $p$-element.

\med
One   observes that this is the case whenever $m$  is   a $p$-power.

Understanding the \ir \reps of groups $G_{r,p}$ introduced above seems to be essential for making a progress with Problem 8. In a sense, these  are simplest groups for which   Problem  8 and  Problem 2  remain open. To emphasize this we state the \f 

\med
Problem 9. Which groups $G_{r,p}$ with $r>2$ have a non-trivial unisingular \ir representation?

\med
The answer depends on ${\rm ord}_p( r)$, the order of $r$ modulo $p$, Lemma \ref{mi2} deals with the case  ${\rm ord}_p( r)=p-1$. 
However we have no general result for ${\rm ord}_p( r)=(p-1)/2$.  Note that Problem 9 is not interesting for $r=2$ as every faithful \rep of $G_{2,p}$ is unisingular.

\med
 Problem 9 is related to a combinatorial problem of linear algebra discussed in a few
publications, especially see \cite[Section 4]{Ca}.

\med
Problem 10. Let $G\subset GL_n(q)$ be a subgroup and $V $ the underlying space of $ GL_n(q)$. Let $W\subset V$ be subspace. When $V=\cup _{g\in G}\, gW$?

\med We are interested with a special case of Problem 10 where $G$ is irreducible
and $\dim W=n-1$. Moreover, we single out here a minimal (in a sense) version   
of Problem 10 which is of significance  in this paper: 

\med
Problem 11. Let $p>r>2$ be  primes and  $G\subset GL_n(r)$ be an \ir subgroup 
of order $p$. Let $V$ be the underlying space for $GL_n(r)$ and $W$ a subspace of dimension $n-1$. When $V=\cup _{g\in G}\, gW$?

\med
In fact, Problem 11 is equivalent to Problem 9:  

\bl{11w}  Let $p>r>2$ be  primes, let $G=G_{r,p}$ and $H=G/O_r(G)$. 
Then the \f are equivalent:

$(1)$  there exists a  faithful unisingular \irr of G over a field of characteristic $\ell\neq r$;
 
$(2)$ $V=\cup _{g\in G}\, gW$ for $V\cong \FF_r^p$   and   some subspace W of V of codimention $1$.
\el

We are not aware whether there can exist unisingular and non-unisingular faithful \ir \reps of  $G=G_{r,p}$. 

\med
This paper mainly focuses on \reps over fields of characteristic 2. Most of questions discussed above deserve to be clarified also for \reps over fields of characteristic $\ell\neq 2$. However, this is not our principal goal. Nonetheless, in Section 8 we prove the \f result:

\begin{theo}\label{fo8} Let $F$ be a field of characteristic $r\neq 2$ and  $H=GL_n(F)$.  Then $H$ has an absolutely \ir unisingular subgroup, unless $n=2,4$ or $n=8,r=3$. 
\end{theo} 

This hints that problems of the existence of unisingular \ir subgroups   in classical groups 
over fields of odd characteristic are simpler that for characteristic 2.

The paper is organized as follows.
In section 2 we collect some elementary facts on unisingular linear groups and finite group \reps and quote some known results. 
 
Section 3 we recall and prove some results of general nature and comment the vector space covering problem.

In section 4 we discuss unisingular \reps of groups $PSL_2(q)$. 

In Section 5 we prove Theorem \ref{se2} and discuss certain aspects of \rep theory of 
algebraic groups for further use in the paper. 

In section 6 we provide examples of unisingular \ir subgroups of $Sp_{2n}(2)$ for $n<125$.

In Section 7 we prove  our results on the absence of unisingular \ir subgroups of $Sp_{2n}(2)$ for some values of $n$ with $1<n<125$.

In Section 8 we discuss unisingular \reps over fields of odd characteristic and prove Theorem \ref{fo8}.
 
\vskip1cm
{\it Notation and some definitions} 
The cardinality of a set $S$ is denoted by $|S|$. We also use $|g|$ to denote the order of a group element $g$. 
$\ZZ$ denotes the ring of integers,  $\ZZ^+$ is the set of non-negative integers and $\NN$  is the set of natural numbers.  
 $\CC$ is the complex number field, $\QQ$ is the field of rational,  $\FF_q$ denotes the finite field of $q$ elements, $\FF_q^\times$ is the multilicative group of $\FF_q$. We write for $\FF_q^n$ the set of column vectors with coordinates in $\FF_q$.
By $\overline{\FF}_q$ is denoted the algebraic closure of $\FF_q$.

For  $m,n\in \NN$ we write $m|n$ to say that $m$ divides $n$ and $(m,n)$ for the greatest integer $k$ dividing both $m,n$.     We write   
$n_p$ for the $p$-part of $n$, the greatest $p$-power dividing $n.$ 
Let $q$ be a prime power and $n>0$ an integer.  We denote by ${\rm ord}_ q(p)$ the minimal integer $k\geq 1$ such that $q^k\equiv 1\pmod p$ (or $p|(q^k-1)$). Equivalently,  $n={\rm ord}_ q( p)$ \ii  $\FF_{q^n}$
is the smallest field extension of $\FF_q$ such that $\FF_q^\times$ has an element of order $p$, and \ii $n$ is the smallest number $k>1$ such that  $GL_k(q)$   has an element of order $p$. Note that  ${\rm ord}_ q( p) \leq  p-1$, in fact, $p|(q^{p-1}-1)$. A prime $p$ is called a {\it primitive prime divisor of} $q^n-1$ if $n={\rm ord}_ q( p)$. 

We denote by $C_p$ the cyclic group of order $p$ and $C_p^d=C_p\times \cdots \times C_p$ the direct product of $d$ copies of $C_p$. 

We denote by $\mathcal{A}_n$ the alternating group on $n$ letters and by $S_n$ the symmetric group. For sporadic simple groups we follow notation in \cite{Atlas}. We use a standard notation for classical groups $GL_n(F)$, $SL_n(F)$, $Sp_n(F)$,  $O_n(F)$; if $F=\FF_q$, we replace $F$ by $q$: $GL_{n}(q)$, $SL_{n}(q)$, etc. In addition,   $\Om^\pm_{2n}(q)$ and $\Om_{2n+1}(q)$ is the subgroup of the respective orthogonal group over $\FF_q$ formed by elements of spinor norm 1.   $U_n(q)$ is the unitary group over $\FF_{q^2}$. We denote 
by $PSL_n(q)$, $PSU_n(q)$, $n>2$, the simple groups obtaining  from linear groups $SL_n(q)$, $SU_n(q)$ by factorization over the center.  To simplify notation, the group $PSL_2(q)$ is denoted here by $L_2(q)$.  
 
To specify the underlying vector space of the natural realization  of classical groups we 
 we often write $GL(V)=GL_n(F)$, $Sp(V)=Sp_{2n}(F)$ etc.; we also use $GL_n(F)$
to denote the group of non-degenerate $(n\times n)$-matrices over $F$. 

Let $g\in GL(V)$, where $V$ is vector space over a field. Then $V^g=\{v\in V: gv=v\}$.  We say that $g$ is {\it fixed point free on} $V$ if $V^g=0$, that is, $gv=v$ for $v\in V$ implies $v=0$, equivalently, 1  is not an \ei of $g$.  A subgroup $G\subset GL_n(F)$ is called {\it unisingular} if every element $g\in G$ has \ei 1. A \rep $\phi:G\ra GL_n(F)$ of a group $G$ is called {\it unisingular} if the group $\phi(G)$ is unisingular.  

If  $G\subset GL_n(V)$ and $W$ is a subspace of $V$ then $GW=\{gw: g\in G,w\in W\}$.

Let   $G$ be a group. We denote by $G^m$ the direct product of $m$ copies of $G$, by  $G'$ the derived subgroup of $G$ and by $Z(G)$ the center of $G$. If $p$ is a prime, a $p'$-element of $G$ is one of order coprime to $p$. By $1_G$ we denote the trivial one-dimensional \rep of $G$ or the trivial character (the grounf field is usually clear from the context.) 

Let $H$ be a subgroup of a group $G$. We write $|G:H|$ for the index of $H$ in $G$.
If $\lam$ is a \rep of $H$, we use $\lam^G$ for the induced \rept If $\mu$ is a \rep of $G$ we write $\mu|_H$ for the restriction of $\mu$ to $H$. Let $F$ be a field and  $M$ be an $FG$-module. The sum of all \ir submodules of $M$ isomorhpic to each other is called a {\it homogeneous component} of $M$. A {\it quasi-homogeneous component}
is the  sum of all \ir submodules of $M$ with the same kernel. (So a quasi-homogeneous component is a direct sum of homogeneous components.) 

If $G,H$ are groups then $G\wr H$ denotes the wreath product of $G$ and $H$.
If $G\subset GL_n(F)$ and $H\subset S_m$ then we keep the notation $G\wr H$ for an imprimitive subgroup of $GL_{mn}(F)$ such that $G^m$ stabilizes $m$ imprimitivity blocks permuted transitively by  $H\subset S_m$. In most situations below $H$ is cyclic
and $m=|H|$.  
A $G$-set is a set $S$, say, on which a group $G$ acts by permutations. For $g\in G$ we set $S^g=\{s\in S: g(s)=s\}$. 

A representation $\phi:G\ra GP_n(F)$ of a group $G$ is called {\it tensor-decomposable} if $\phi(G)$
is contained in a subgroup $GL_k(F)\otimes  GL_m(F)\subset \GL_n(F)$ for some integers
$k,m>1$ with $km=n$. Otherwise, we call $\phi$ tensor-indecomposable. 

We write $G=N\rtimes H$ to denote the semidirect product of groups $N$ and $H$ with
$N$ normal in $G$. Sometimes we simply write $NH$ if the structure of $NH$ is clear from the context.

\newpage
\section{Preliminaries}

\def\hw{highest weight }

We first record the \f well known fact:

\bl{c21}  \cite[8.2, 8.3]{Asch}  Let G be a finite group and N  a minimal normal subgroup of G.  Then  $N$ is either elementary abelian or a direct product of non-abelian simple groups isomorphic to each other.\el

\bl{gw1} Let $F$ is an \ac field,  $g\in GL(V)\cong GL_n(F)$, and let $V=V_1\oplus \cdots\oplus V_m$ be a direct sum decomposition such that g transitively permutes $V_1,..., V_m$. 
 Then the \eis of g are all $m$-roots of the \eis of 
   $g^m$, in fact of those on $V_1$.  \el

\bp Observe that $g^m$ stabilizes every  $V_i$ for $i=1\ld m$. In addition, the $V_i$'s are $F\lan g^m \ran$-modules isomorphic to each other. Therefore, the \eis of $g^m$ on $V_i$ are the same as on $V_1$. Let $\lam$ be an \ei of  $g^m$ on $V_1$ and $v\in V_1$ is an $\lam$-eigenvector. Then vectors $v,gv\ld g^{m-1}v$ are linear independent and $W:=\lan v,gv\ld g^{m-1}v\ran$ is a $g$-invariant subspace of dimension $m$. \itf the characteristic \po of $g$ on  $W$ is $x^m-\lam$,
whose all roots are the \eis of $\lam$ on $W$.  \enp

\bl{gv1} Let $F$ be an arbitrary field, $g\in GL(V)\cong GL_n(F)$ and let $V=V_1\oplus \cdots\oplus V_m$ be a direct sum decomposition such that g transitively permutes $V_1,..., V_m$. 
 Then  g   has \ei $1$   \ii $g^m$ has \ei $1$ on $V_1$.\el

\bp This follows from Lemma \ref{gw1} as $g$ has \ei \ii $g$ has \ei 1 as an element of $GL_n(\overline{F})$, where $\overline{F}$ is an \ac field containing $F$.\enp

We recall Clifford's theorem (see for instance \cite[Theorem 12.13]{Asch} or \cite[Theorem 49.7]{CR}.

\bl{ct3} Let $G\subset GL(V)=GL_n(F)$ be an \ir subgroup and N a normal subgroup of G. 

$(1)$ V is a direct sum of \ir $FN$-modules of the same dimension;

$(2)$ G permutes transitively the homogeneous components of $N$ on V;

$(3)$ Let $V_1$ be a homogeneous components of $N$ on V and $G_1=\{g\in G: gV_1=V_1\}$. Then $V_1$ is an \ir $FG$-module.\el

The following lemma is trivial but we state it for reader's convenience; this is used throughout the paper with no explicit reference. 

\bl{af7} Let F be a field and $\overline{F}$ the algebraic closure of F. Let $G\subset GL_n(F)$ be a finite group. Then G is unisingular \ii G is so as a subgroup of  $ GL_n(\overline{F})$.\el

Let $G\subset GL_n(q)$ be an \ir subgroup, and $P$ an extension field of $\FF_q$. Then $P$ is called a {\it splitting field for} $G$ if all the composition factors   of $G$ as a subgroup of $GL_n(P)$ are absolutely \irt There exists a unique minimal splitting field and this is $\FF_{q^k}$, where $k$ is the number of  the composition factors in question. 
Moreover, if $\tau$ is any composition factor then $\FF_{q^k}$ is the minimal field such that $\tau(G)$ is equivalent to a representation into $GL_{n/k}(\F_{q^k})$. In fact, $\FF_{q^k}$ is obtained from $\FF_q$ by adding all traces of $\tau(g)$ for $g\in G$. (This does not depend on the choice of $\tau$.) See \cite[Theorem 19.4]{Fe}. Some details are in the \f lemma:
    
\bl{na3}  Let $G\subset GL_n(q)$ be an \ir unisingular subgroup. Suppose that G is not absolutely \irt Then $G$ is isomorphic to a unisingular  absolutely \ir subgroup of $GL_{n/k}(q^{k})$ for some integer $k|n$, $k>0$. More precisely, $\FF_{q^k}$ is the minimal splitting field for G.\el

\bp Let $V$ be the underlying space for $GL_n(q)$ and $M={\rm Mat}_n(\FF_q)$. By Schur's lemma, the $\FF_q$-enveloping algebra $[G]$ of $G$ is simple,  the center of it is a field isomorphic to $\FF_{q^k}$ for some integer $k>1$ with $k|n$ and $[G]={\rm Mat}_{n/k}( \FF_{q^k})$, see \cite[Ch. V, Theorem 19]{KL}.
Let $V_1=V\otimes \FF_{q^k}$. Then $V_1$ is a direct sum of ${\rm Mat}_{n/k}(\FF_{q^{k}})$-modules permuted transitively by ${\rm Gal}( \FF_{q^k}/\FF_q)$. This yields $k$ absolutely \ir \reps $G\ra {\rm Mat}_{n/k}(\FF_{q^{k}})$, which are Galois conjugate to each other. As $g\in G$ has \ei 1  on $V$ \ii $g$ has \ei 1 on $V_1$, the result follows. \enp

 \bl{rs1} Let $G\subset Sp_{2n}(2)$ be a unisingular \ir  subgroup which is not absolutely irreducible. Then $G\subset H\cong U_{n/k}(2^k)$ for some divisor $k$ of n, and G is unisingular and absolutely \ir as a subgroup of H.\el

\bp Let $X=C_{GL_{2n}(2)}(G)$. We have seen in the proof of Lemma \ref{na3} that
$X$ is the multiplicative group of a field $\FF_{2^k}$, where $k| 2n$. In addition,  $C_{GL_{2n}(2)}(X)\cong GL_{2n/k}( 2^k)$ and the $\FF_2$-span of $G$
in $ {\rm Mat}_{2n}(\FF_2)$ is isomorphic to ${\rm Mat}_{2n/k}(\FF_{2^k})$.  
Therefore, $G$ is absolutely \ir as a subgroup of ${\rm Mat}_{2n/k}(\FF_{2^k})$.
Let $\si$ be an \au of $GL_{2n}(2)$ such that $Sp_{2n}(2)=\{h\in GL_{2n}(2): \si(h)=h\up\}$. Then $\si(X)=X$. Note that $V$ is a homogeneous $\FF_2X$-module
(as $G$ is irreducible).   
 By \cite[Lemma 6.6(2)]{EZ}, we have $G\subset H\cong U_{n/k}(2^ {k})$. In fact, the underlying spaces of $ Sp_{2n}(2)$ and  $U_{n/k}(2^{k})$ coincide so $G$ is unisingular as a subgroup of $U_{n/k}(2^ {k})$.  \enp

For a finite group $X$ let $O^\ell(X)$ be the minimal normal subgroup $N$ of $X$ such that $X/N$ is an $\ell$-group. Note that a subgroup of $GL_n(2)$ is usingular \ii every odd order element of has \ei 1. 

\bl{tf1}  Let M be a finite group, $\ell$ a prime and let F be a field of characteristic $\ell$. Let $\phi_1,\phi_2$ be \ir $F$-representations  of M. Then $\phi_1$ is equivalent to $\phi_2$ \ii their restrictions to $O^\ell(M)$ are equivalent. In addition, $\phi_1$ is unisingular \ii so is $\phi(O^\ell(M))$. \el

\bp By \cite[Corollary 17.10]{CR1}, $\phi_1$ is equivalent to $\phi_2$ \ii their Brauer characters  coincide. By the definition of $O^\ell(M)$, the $\ell'$-elements of $M$ lies in $O^\ell(M)$, whence the result. This also implies the additional claim. \enp

 \bl{tt1}  Let $\be$ be the Brauer character of an absolutely \ir $2$-modular \rep $\rho$ of a finite group $G$, and $d=\be(1)>1$. Then $\rho(G)$ is equivalent to a \rep into $Sp_d(2)$ \ii $\be(g)\in \RR$  and
 $\be(g)\pmod 2\in \FF_2$  for every odd order element $g\in G$. In particular,
if the values of $\be$ are integers then $\rho(G)$ is equivalent to a \rep into $Sp_d(2)$. \el

\bp By \cite[Ch. I,  Theorem 19.3]{Fe}, $\rho(G)$ is equivalent to a \rep into $GL_d(2)$ \ii $\be(g)\pmod 2\in \FF_2$  for every $2'$-element $g\in G$.
In addition, $\be(g)\in \RR$ for every $2'$-element $g\in G$ \ii  $\rho$ is self-dual \cite[Ch. IV,  Lemma 2.1]{Fe}. In turn, this is equivalent to the inclusion $G\subset Sp_d(2)$ \cite[Ch. IV,  Theorem 11.1 and Corollary 11.2]{Fe}.\enp

The  following lemma states well known facts, see also \cite[Lemma 3.13]{z22}. 

\bl{ek1} $(1)$ Let F be a field and $g\in GL_n(F)=GL(V)$. Let $\si:F\ra F$ be a field \au extended to $GL_n(F)$ in the natural way. Then $\si(g) $ has \ei $1$ \ii $g$ has.   

$(2)$ Let $\mathbf{G}$ be an algebraic group, $\si$ a standard Frobenius endomorphism of $\mathbf{G}$ and $\phi:\mathbf{G}\ra GL(V)$ be a rational \rept Let $g\in \mathbf{G}$.
Then $\phi(g)$ has \ei $1$ \ii $\phi^\si(g)$ has. 
\el

\bp (1) Extend $\si$ to $V$. Then every subspace of $V$ is $\si$-invariant. Let $W$ be a maximal $g$-stable subspace of $V$ such that $g$ has \ei 1 on $V/W$. 
Then $\si(g)W=W$, and hence it suffices to prove the lemma for $V/W$ in place of $V$.
Then the 1-eigenspace of $g$ is $\si$-invariant and hence $\si(g)$-invariant. Whence the claim. 

(2) It is well known that $\phi(\si(g))=\si(\phi(g))$, see \cite[\S 1.17]{C}. So the claim follows from (1).  \enp

\bl{ct1}  Let V be a vector space over an arbitrary  field,  let $G\subset GL_n(V)$ be an \ir group and N a normal subgroup of G. Suppose that $V|_N$ is reducible and $|G/N|$ is prime. Then $V|_N$ is the sum of pairwise non-isomorphic \ir $N$-modules.\el

\bp By Clifford's theorem, $V|_N$ is a direct sum the homogeneous components, $V_1\ld V_k$, say, which are transitively permuted by $G.$ Let $G_1=\{g\in G:gV_1=V_1\}$.
Then $G_1$ acts in $V_1$ irreducibly (again by Clifford's theorem). As $|G/N|$ is prime, we have $G_1=N,$ so $V_1|_N$ is \irt \enp

\bl{zz1} Let F be a field, $G\subset GL_n(F)$  a unisingular \ir subgroup, and let $A$ a non-trivial abelian normal subgroup. 

$(1)$  A is not cyclic, and hence  reducible.  

$(2)$ Suppose that G is primitive. Then G has no non-trivial abelian normal subgroup.  \el

The following lemma is often used with no reference. 

\bp  
(1) Let   $V$ be the underling space of $GL_n(F)$, $1\neq a\in A$ and let $V^a$ be the 1-eigenspace of $a$ on $V$. Then $GV^a=V^a$, a contradiction unless $V^a=V$, but this means that $a=\Id$. If $A$ is \ir  then $A$ is cyclic by Schur's lemma. 

(2) Suppose the contrary, and let $A$ be such a group. By Clifford's theorem, $V|_A$ is homogeneous, so it is a direct sum of isomorphic \ir $A$-modules. Let $W$ be one of them. Then $V|_A$ is unisingular \ii so is $W$. Then the claim follows from Lemma \ref{zz1}. 
\enp

\bl{od1} Let $G$ be a finite group  and $V=V_1\oplus V_2$, where $V_1,V_2$ are dual FG-modules. (Here F is an arbitrary field.)

$(1)$  $V$ is unisingular \ii so is $V_1$.

$(2)$ Suppose that $V_1$ is homogeneous. Then $V$ is unisingular \ii so is every \ir constituent of $V$.

$(3)$ Suppose that $V_1$ is homogeneous and G is abelian. Then V is not unisingular.
 \el

\bp (1), (2) are obvious. (3) Let $V'$ be an \ir $FG$-submodule of $V$, and $K$ is the kernel of $V_1$. By Schur's lemma, $G/K$ is cyclic, so $V'$ is not unisingular by Lemma \ref{zz1}. So (2) implies (3).\enp

\bl{in5a}  Let $F$ be an arbitrary field, let $G\subset GL_{n}(F)$ be a unisingular  \ir  subgroup and $m=kn$. Then $H=G\wr C_k$ is a unisingular \ir subgroup of $GL_{m}(F)$. If n is even and $G\subset Sp_{n}(F)$ then $H\subset Sp_{m}(F)$.\el
 
\bp Note that $G\wr C_k$ has a normal subgroup $N\cong G_1\times \cdots\times G_k$, where $G_i\cong G$ for $i=1\ld k$.
Let $V=V_1\oplus \cdots \oplus V_k$, where $V_i$ is the natural $FGL_{n}(F)$-module, viewed as an $F G_i$-module due to the isomorphism $G_i\ra G$. Then $V$ can be turned to $FH$-module, which is obviously
absolutely irreducible. 

If $G\subset Sp_{n}(F)$ then $H\subset Sp_{m}(F)$, where we  regard $V_1\ld V_k$ as non-degenerate subspaces of $V$. We show that $V$ is a unisingular $F H$-module. Let $h\in H$, and let $m\geq 0$ be the minimal integer such that $h^m\in N$. As $C_k$ is cyclic and transitively permutes $V_1\ld V_k$, it follows that the $h$-orbits on $V_1\ld V_k$ are of size $m$.  
Then $h^m=\diag(g_1\ld g_k)$, where $g_i\in G_i$. Obviously, $V^{h^m}=V_1^{g_1}+\cdots +V_k^{g_k}$. As $V^{h^m}$ is $h$-stable, it follows that $h$ permutes  $V_1^{g_1}\ld V_k^{g_k}$  and the $h$-orbits  on $V_1^{g_1}\ld V_k^{g_k}$ are of size $m$.  As $G$ is unisingular, we have $V_1^{g_1}\neq 0$. Let $0\neq v\in V_1^{g_1}$. Then the vectors $v,hv\ld h^{m-1}v$ are linearly independent and $h^mv=v$. Then $v+hv+\cdots+ h^{m-1}v  \in V^h$, as required. \enp

\begin{corol}\label{in5} If $Sp_{2n}(2)$ contains a unisingular \ir (respectively, absolutely irreducible) subgroup then $Sp_{2nk}(2)$ contains a unisingular \ir (respectively, absolutely irreducible)  subgroup for every integer $k>1$.  In addition, $Sp_{8k}(2)$ contains a unisingular \ir (respectively, absolutely irreducible)  subgroup.\end{corol}
 
\bp The additional claim follows   as $Sp_{8}(2)$ has a unisingular absolutely \ir  subgroup \cite{CZ}.\enp

\bl{a22} {\rm \cite[Theorem 1.1]{W}} Let $G=\mathcal{A}_n,n>8$, be the alternating group. Then the minimal degree of a non-trivial $2$-modular   \rep   is $n-1$ if n is odd
and $n-2$ if $ n$ is   even. These \reps are realizable  over $\FF_2$. \el

 Note that the minimal degree of a non-trivial $2$-modular   \rep of $A_8$, $A_7$, $A_6$ and $A_5$  is $4,4,4,2$, respectively. 

\bl{aa3} Let $n\geq 4$ be even. Then $A_{2n+1}$ and $A_{2n+2}$  have  \ir $2$-modular \reps  of degree $2n$, and elements of order $2n+1$ do not have \ei $1$ in an \irr of  degree $2n$. In addition, $A_{2n+2}$ has no \ir $2$-modular \reps  of degree $2n+1$.\el

\bp This follows from Lemma \ref{a22} for $n>8$, for $n\leq 8$ see \cite{JLPW}. For the additional statement use    \cite[Corollary 2.4]{zk}.  \enp
 
\med
I am ibdebted to Pablo Spiga for the \f lemma:

\bl{pp1}
Let $G$ be a finite group acting transitively on a set $\Omega $ and P a \syl of G. Then every P-orbit size   is a multiple to  $|\Omega|_p$. In particular,
if  $|\Omega|$ is a p-power then $P$ is transitive.\el

\bp Let $\al\in \Omega$ and let $G_1,P_1$ be the stabilizers of $\al$ in $G,P$, respectively.
Let $P\al$ be the $P$-orbit containing $\al$, so $|P\al|=|P:P_1|$. Then $|G:P_1|=|G:P|\cdot |P:P_1|=|G:G_1|\cdot |G_1:P_1|=|\Omega|\cdot |G_1:P_1|$. As $|G:P|$ is coprime to $p$, we conclude that $|P:P_1|$ is a multiple of $|\Omega|_p$.

By Sylow's theorem, $P_1$ is a \syl of $M$ for some \syl $P$ of $G$. Then $|P\al|=|P:P_1|=|G:G_1|_p=|\Omega|_p$. As Sylow $p$-subgroups are conjugate,
$|P\al|=|\Omega|_p$ for every \syl $P$ of $G$.

One observes that at least one of the $P$-orbit has size $|\Omega|_p$. To see this, choose $P$ so that $P_1$ is a \syl of $G_1$. Then $(|G_1:P_1|, p)=1$, whence $|G:G_1|_p=|P:P_1|$.\enp

\bl{a33} Let $G$ be a transitive subgroup on a set $J$ and $|J|=n$. Let O be an orbit of G on the unordered pairs  $(a,b)$, $a,b\in J$, $a\neq b$. Then either $|O|\geq n$ or $|O|=n/2$.\el

\bp Let $O_1=G(a,b)$ be the $G$-orbit on the ordered pairs $(a,b)$, $a,b\in J$, $a\neq b$ and  $G_{a,b}$ the stabilizer of $a$ and $b$ in $G$. Let $G_{(a,b)}$ be the stabilizer of the unordered pair  $(a,b)$ in $G$. Then $G_{a,b}\subseteq G_{(a,b)}$ is a subgroup of index at most 2. So $|O_1|=|O|\cdot |G_{(a,b)}:G_{a,b}|$. Observe that $|O_1|=|G:G_a|\cdot |G_a:G_{a,b}|=n\cdot |G_a:G_{a,b}|$. So  

$$|O|=\frac{|G:G_a|\cdot |G_a:G_{a,b}|}{|G_{(a,b)}:G_{a,b}|}=\frac{n\cdot |G_a:G_{a,b}|}{|G_{(a,b)}:G_{a,b}|}.$$

So $|O|\geq n$  unless   $|G_a:G_{a,b}|<|G_{(a,b)}:G_{a,b}|=2$. This latter  $|O_1|= n$ and $|O|= n/2$. \enp  

\bl{si2}
Let $G\subset GL_n(q)$ be a subgroup with  \ir normal $r$-subgroup R for some prime r. Then $G$ is not unisingular.\el

\bp By Schur's lemma, $Z(R)$ is a non-trivial cyclic normal subgroup of $G$  by Schur's lemma. So the result follows from Lemma \ref{zz1}. \enp
  
\bl{am1} 
Let $G\subset GL(V)=GL_n(F)$ be a  unisingular \ir subgroup,  and $N$  a normal subgroup of $G$. Let $V|_N=V_1+\cdots +V_l$, where $V_1\ld V_l$ are the homogeneous components of $N$ on $V$.   Suppose that $l$ is a p-power and let S be a \syl of $G$. Then  $Z(N)\cap Z(S)=1$. \el

\bp Suppose the contrary, that $Z=Z(N)\cap Z(S)\neq 1$. By Lemma \ref{pp1},  $S$  transitively permutes  $V_1\ld V_l$; these are isomorphic to each other as $FZ$-modules (since $Z\subset Z(S)$). In addition, $V_1|_Z$ is homogenious, hence so is $V|_Z$. 
This is a contradiction, as $Z$ is abelian (Lemma \ref{od1}).  \enp
 
\bl{it4} Let $G\subset GL_n(F)=GL(V)$ be an  \ir subgroup, and $A$ an abelian normal subgroup of G.  Let $V=W_1+\cdots +W_l$, where $W_i$ are quasi-homogeneous components of $V|_A$. 
Suppose that A is unisingular and $l$ is a p-power.  Then $(|A|,p)=1$.\el

\bp  Suppose the contrary, and let $B$ be the Sylow $p$-subgroup of $A$.
Then $B$ is normal in $G$ and in $S$, where $S$ is a \syl of $G$. This implies  $[b,S]=1$ for some $1\neq b\in B$. As $l$ is an $p$-power, $S$ transitively permutes  $W_1\ld W_l$ (Lemma \ref{pp1}). We can assume that $b$ is non-trivial on $W_1$. Then $b$ is non-trivial on every \ir constituent of $W_1|_A$, and hence $b$ acts  fixed point freely on $W_1$. Let $s\in S$ be such that $sW_1=W_i$. Then $bs=sb$ implies that $b$ acts  fixed point freely on $W_i$, and hence on $V$ as $i\in \{1\ld l\}$ is arbitrary. This is a contradiction. \enp 

\bl{d8d} Let $r>2$ be a prime, F be an arbitrary field of characteristic $\ell \neq r$,  and let $A\subset GL_n(F)=GL(V)$ be a unisingular abelian r-group of rank  $d$. Set $m=F(\zeta):F$, where $\zeta$ is a primitive r-root of unity.  
If $A$ fixes no vector $0\neq v\in V$ then   $1<d<(n-m$)/m.  \el

\bp  Let $A_0$ be the subgroups of $A$ of elements of order at most $p$. Then the rank of $A_0$ equals $d$. It suffices to prove the lemma for $A=A_0$, so we assume that $A$ is elementary abelian. As $C_V(A)=0$, $V$ is a direct sum of non-trivial \ir $FA$-submodules,
each of dimension $m$, so $n/m$ is an integer and $A$ is isomorphic to a subgroup of $GL_{n/m}(F(\zeta))$. So we can assume that $F=F(\zeta)$, and then we can assume that $A$ is a subgroup of the group $D$  of  diagonal matrices  $x\in GL_n(F)$
with $x^p=1$.   Then $D$ is of rank $n$. Clearly, $A$ contains no non-identity scalar matrix.

Suppose first that $d=n$. Then   $d=n$ and   $ A=D$. Then $A$ contains a non-identity scalar matrix, a contradiction. 

Suppose that $d=n-1$.   Let $\nu$ be the natural \ho of $GL(V)$ onto $PGL(V)$. Then the mapping $A\ra \nu(A)$ is injective, whereas the rank of $\nu(D)$ equals $n-1$.
So  $\nu(A)=\nu(D)$. Therefore, for every $x\in D$ there exists a scalar matrix $z$ such that $zx\in A$. 

Suppose first that   $x=\diag(\eta, 1\ld 1)\in A$. Let $X=\diag(1,GL_{n-1}(F))$ and $D_1= D\cap X$.  Set $A_1=A\cap D_1=A\cap X $. We can view $A_1$ as a subgroups of $X_1:=GL_{n-1}(F)$; then $A_1$ as a subgroups of $X_1$ is unisingular, as otherwise $xg_1$ does not have \ei 1 whenever $g_1\in A_1\subset X_1$ does not have. In addition, the rank of $A_1$ equals $n-2$.  
The case $n=2$ is trivial so we can use induction on $n$. By induction assumption,
 $A_1$ is not unisingular in $X_1$, and we are done in this case. 

Suppose that $x\notin A_1$. Then $z=\eta\up\cdot\Id_n$ is the only scalar matrix such that $xz$ has \ei 1, so $y:=xz\in A_1$ by the above. So
$y=\diag(1,\eta\up\cdot\Id_{n-1})\in A_1$. As  above, we can assume that   $x'=\diag( 1\ld 1,\eta)\notin A_1$ and $y'=\diag(\eta\up\cdot\Id_{n-1}, 1)\in A_1$. Then $yy'\in A_1$ does not have \ei 1, a contradiction.  \enp

Note that the assumption $p>2$ cannot be dropped. 

\med The \f lemma is well known.

\bl{mm9} Let $p,r$ be distint primes and $k>0$ and integer. Then the \f are equivalent:

$(1)$ $k={\rm ord }_ r( p);$

$(2)$  $\FF_{r^k}$ is the minimal field of characteristic $3$ whose multiplicative group contains an elt of order $p$.

$(3)$ k is the minimal positive integer such that $GL_k(r)$ contains an element of order p.
\el

 Let  $3<p<250$ be a prime such that $(3^i-1,p)=1$   for every $0<i<p-1$. Then \med  $p\in\{5,7,17,19,29,31,43,53,79,89,101,113,127,137,139,149,163,173,197,199,\\211,223, 233\},$
see \cite[Table A062117]{OL1}. 
In Table 1, for further use, we tabulate some data extracted from  \cite[Table A062117]{OL1}.

\medskip 
\begin{table}[ht]\label{tab1}
\begin{center}
\caption{Order of $3$ modulo some primes $p$}  

\bigskip
\begin{tabular}{|c|c|c|c|c|c|c|c|c|c|c|c|c|c|c|c|} 
\hline
$p$&11&23&29&41&43&47&53&67&73&83&89 &103&107&113\\
\hline
ord$\,3\mod p$&5  &11&28&8&42&23&52&22 &12&41&88&34&53&112  \\
\hline
\end{tabular}
\end{center}
\end{table}
\med

\bp[Proof of Lemma {\rm \ref{mi2}}]  Let $A$  be an elementary abelian normal $r$-subgroup of $G$, and let $d$ be the rank of $A$. We can identify $A$ with the additive group of a vector space over $\FF_r$ of dimension $d$. Then $G/A$ can be identified with a subgroup $C$, say, of $GL(V)=GL_d(\FF_r)$.  As $A$ contains no non-trivial $G$-invariant subgroup other than $A$ itself, we conclude that $C$ is an \ir subgroup of $GL_d(\FF_r)$. By Lemma \ref{mm9}(3), $d={\rm ord}_r(p)$. By assumption, this equals $p-1$. 

Let $\phi$ be an \irr of $G$ over a field $F$ of characteristic $\ell\neq r$.   Then either $\phi$ is faithful or $A$ is the kernel of $\phi$. In the latter case the lemma is true due to Lemma \ref{zz1}. Suppose that $\phi$ is faithful. We can assume $F$ algebraically closed. In addition, as every \irr of a solvable group lifts to characteristic 0, see \cite[Ch. X, \S 2, Theorem 2.1]{Fe}, we can assume  $F$  of characteristic 0.  Then $\dim\phi=p$
by Ito's theorem \cite[Corollary 53.18]{CR}. As $\phi$ is irreducible,  $\phi(A)$ has no trivial \ir constituent. By Lemma \ref{d8d}, the rank of $A$ is at most $p-2$, which is a contradiction.   \enp
 
\newpage
\section{Some general results}

Denote by $\AGL_n(q)$ the semidirect product of $GL_n(q)$ and the additive group
of the natural $\FF_qGL_n(q)$-module. 

\begin{theo}\label{af1}  {\rm \cite[Theorem 1.6]{CZ}} Let $q$ be odd and $m=q^n-1$.  Suppose that $n>1$ or $n=1$ and $q$ is not a prime. Then there exists a unisingular  absolutely \ir  subgroup of $\Sp_{m}(2)$ isomorphic to $\AGL_n(q)$. \end{theo}
  


Let $N$ be a group and $M$ a completely reducible $KN$-module over some field $K$. 
Recall that homogeneous component of $M$ is a the sum of all \ir submodules of $M$
isomorphic to each other. So $M$ is a sum of its homogenious component. 
A {\it quasi-homogeneous component of} $M$ is the sum of all \ir submodules of $M$
with the same kernel. Then there exists a unique decomposition $M=M_1+\cdots +M_k$
such that $M_1\ld M_k$ are quasi-homogeneous components of   $M$. We call such a decomposition the {\it quasi-homogeneous decomposition}  of $M$.  These notions are not customary   but very useful when $N$ is abelian. 

\subsection{Remarks on Clifford's theory}


\bl{ct2}  Let V be a symplectic space over a field F and $G\subset Sp(V)$ be an \ir subgroup. Let N be a reducible normal subgroup of G,   and let $V_1\ld V_k$ be the homogeneous components of N on V. Then  $V_1\ld V_k$ are transitively permuted by G and either

$(1)$ all $V_1\ld V_k$ are non-degenerate and orthogonal to each other or  

$(2)$ all $V_1\ld V_k$ are totally isotropic, $k=2l$ is even, and after a suitable reordering the subspaces  $V_{2i-1}+V_{2i}$ $(i=1\ld l)$ are non-degenerate, orthogonal to each other  and transitively permuted by $G$. In addition, $V_{2i-1},V_{2i}$ are dual.

$(3)$  Quasi-homogeneous components of $V|_N$ are non-degenerate subspaces of V of the same dimension dividing $\dim V$, and they are transitively permuted by G. \el

\bp Statement (1),(2) is a refinement of Clifford's theorem for subgroups of classical groups \cite[Proposition 5]{Z71} specified for symplectic groups.


(3) In (2) $V_{2i-1}$ and $V_{2i}$ are dual, so they are in the same quasi-homogeneous component. Therefore, the quasi-homogeneous components are non-degenerate subspaces of $V.$  As $V_i=g_iV_1$ for some $g_i\in G$, it follows that $g_iK_1g_i\up$ is the kernel of $V_i$ for $i=1\ld k$. This easily implies the claim.  
\enp

Lemma \ref{2po3} is contained in the \f lemma:

\bl{2po}  Let $p>2$ be a prime, $G$  a unisingular \ir subgroup of $Sp_{2p}(2)$ and N a minimal normal subgroup of G. Then one of the \f holds:

$(1)$ $N$ is  simple and either \ir or has two \ir constituents which are dual to each other;

$(2)$ N is an elementary abelian $3$-group,  
 the \ir constituents of N are non-isomorphic and of dimension $2$ each.
 
$(3)$ G contains an \ir subgroup isomorphic to $G_{3,p}$. 
\el

\bp 
Let $V$ be the underlying space of $Sp_{2p}(2)$.  Suppose first that $N$ is \irt  
  If $N$ is  abelian then, by Schur's lemma,  $N$ is cyclic, which contradicts  Lemma \ref{zz1}. So $N$ is non-abelian, and hence  $N$ is a direct product of non-abelian simple groups. Let $S$ be one of them.
  If $S=N$, we are done, so we assume that $S\neq N$. Then $N=S\times N_1$, where
$N_1\neq 1$ is a sirect product of copies of $S$. By Schur's lemma,  $S$ is reducible. By Clifford's theorem, $V$ is a direct sum of \ir $FS$-modules $V_1\ld V_l$ of equal dimension. So $l|2p$, and $l\neq p$ as $\dim V_1=2$ implies $S$ to be  solvable. So $l=2$ and $V_1,V_2$ are of dimension $p$. These are isomorphic $FS$-modules  (since $[S,N_1]=1$. 
As $\dim V$ is a prime, $S$ is absolutely \ir on $V_1$ and $V_2$ (see Lemma \ref{na3}). 
Therefore, there is an embedding $N_1\ra GL_2(\overline{\FF}_2)$ (see \cite[Theorem 11.20]{CR1}). Simple subgroups of $GL_2(\overline{\FF}_2)$ are known to be isomorphic to $SL_2(2^t)$ for some integer $t>1$. So $N/S\cong SL_2(2^t)$, and hence $S\cong SL_2(2^t)$. However, the degrees of absolutely irreducible \reps of $SL_2(2^t)$  over  $\overline{\FF}_2)$ are well known to be 2-powers (see Lemma \ref{sl2} below or elsewhere).   
This is a contradition. 




Next suppose that $N$ is  reducible. Let $V_1\ld V_k$ be the homogeneous components of  $V|_N$.  Then $k|2p$, so $k\in\{1,2,p\}$. If $k=p$ then $\dim V_i=2$ for $i=1\ld p$. Then we arrive at case (2) of the statement.

  Let $k=1$. Then $V|_N$ is homogeneous, that is, a sum of $\FF_2N$-modules isomorphic to each other. Let  $m$ be the dimension of any of them. Then $m>1$ and $m|2p$, so $m=p$ or $m=2$. In the latter case $|N|=3$, hence $N$ is cycllc, violating Lemma \ref{zz1}.    So $m=p$.
  Then     $V|_N=W_1\oplus W_2$, where $W_1, W_2$ are isomorphic \ir $\FF_2N$-modules. Let $\be$ be the Brauer character of $V|_N$. By Lemma \ref{tt1},  $\be$ is integrally  valued.  The Brauer character of $W_1$ is $\be/2$. As every value of a Brauer character is   an algebraic integer, we have $\be(g)/2$ is a rational algebraic integer. This implies $\be(g)/2$ is integer. By \cite[Ch. IV, \S 11, Corollary 11.2]{Fe}, $\dim W_1$ is even, a contradiction.

 Let $k=2$. So  $V_1,V_2$ are non-equivalent \ir $\FF_2N$-modules. By Lemma \ref{ct2}, $V_1,V_2$ are either totally isotropic and dual to each other, or both  non-degenerate. The latter is ruled out as non-degenerate subspaces have even dimensions.  So  (1) holds.

(3)
Let $V_1\ld V_p$ be the homogeneous  components of $V|_N$. Then $\dim V_i=2$ for $i=1\ld p$,
so each $V_1\ld V_p$ is \ir $\FF_2N$-module. As $G$ acts transitively on $V_1\ld V_p$, it follows that 
there is a $p$-element $h\in G$   which transitively permutes  $V_1\ld V_p$. Then $h^p$ stabilizes 
$V_1\ld V_p$, and hence is of exponent at most $6=|GL_2(2)|$. So $h^p=1$. By  \cite[Corollary 45.5]{CR}, the group $\lan N, h \ran$ is \irt Let $N_1\subseteq N$ be a minimal non-trivial $h$-invariant subgroup of $N$. We show that $X:=\lan N_1, h\ran$ is \irt For this it suffices to show that the $\FF_2N_1$-modules
$V_i|_{N_1}$ are not isomorphic. Clearly, $X$ is normal in $\lan N, h \ran$. If $X$ is reducible then 
$V$ is a direct sum of non-trivial \ir  $\FF_2X$-modules of the same dimension $d$, say, and $d|2p$. We have $d>2$ as 
$(|h|,6)=1$. Hence $d=p$ or $2p$, in the latter case we are done. Let $d=p$. If $U$ is an \ir $\FF_2X$-module then $U|_{N_1}$ is a sum of   \ir $N_1$-modules of dimension 2, which is a contradiction.  \enp

 \bl{im4a}  Let $G\subset GL_{n}(\overline{\F}_p)=GL(V)$ be a primitive unisingular subgroup, and $N$   the product of minimal normal subgroups of G. Then   $N=S_1\times \cdots \times S_k$, where $S_1\ld S_k$ are non-abelian simple groups for some integer $k\geq1$ and
  $V|_N$ is a homogeneous $FN$-module.   Moreover, if W is an \ir constituent of $V|_N$ then W   is an external tensor product $W_1\otimes \cdots \otimes W_k$, where $W_i$ for  $i\in\{1\ld k\}$ is a faithful \ir unisingular $FS_i$-module. \el

\bp By Lemma \ref{c21},  $N$ is a direct product of simple groups, $S_1\times\cdots  \times S_k$, say.  As $G$ is primitive, by Clifford'd theorem, the module $V|_N$ is homogeneous, that is, a direct sum of copies of $W$.   As $V|_N$ is unisingular, so is $W$.

 By Clifford's theorem applied to  $W$, this is an external
tensor product of $S_i$-modules $W_i,$ $i=1\ld k$. As $W$ is a unisingular $N$-module, it follows as above that each $W_i$ is a unisingular $S_i$-module. \enp

\subsection{Unisingular groups and subspace covering problems}

\def\hc{homogeneous component}
\def\hcs{homogeneous components}

Let $G$ be a group, $F$ s field and $V$ an $FG$-module. Recall that a quasi-\hc  \, of $V$ means the sum of all \ir submodules with the same kernel. This notion is more useful for $G$ abelian and  
$V$ completely reducible. In this case $V$ has a unique decomposition as a direct sum of its quasi-homogeneous components.

\bl{aa8} Let $A\subset  GL_n(F)=GL(V)$  be a finite  abelian group  of odd order, 
where $F$ is a field of characteristic $\ell$ coprime to $|A|$. Let $W_1,\ldots ,W_l$ be the quasi-\hcs \, of V, and let $K_i$ is the kernel of $W_i $ for $i=1\ld l$.  

$(1)$  $A$ is unisingular \ii $K_1\cup\cdots  \cup K_l=A;$ 

$(2)$  If A is a unisingular then $l>2$; if $A$ is a unisingular  p-group of rank   r    and $C_V(A)=1$ then     $p\leq \ell-1$ and $r\leq l-2$.   \el

\bp Let $a\in A$ and let $U_i$ be an \ir submodule of $W_i$. As $A$ is abelian and $U_i$ is irreducible, $a$ fixes a non-zero vector of $W_i$ \ii $a\notin K_i$. \itf    $V^a=0$ \ii $a\notin K_1\cup\cdots\cup K_l$.  This implies (1). If $l=2$ then $|A|=|K_1\cup  K_2|$; as $|K_i|\leq |A|/2$ for $i=1,2$ and $|K_1\cap  K_2|\geq 1$, it follows that  $|A|>|K_1\cup  K_2|$, a contradiction.

(2).  The inequality $\ell\geq p+1$ is trivial. Indeed, we can assume that $A$ is elementary abelian. 
Then $|K_i|=p^{r-1}$, $|K_i\cap K_j|=p^{r-2}$ for $1\leq i<j\leq l$. Therefore, $|K_1\cup\cdots  \cup K_l|\leq p^{r-1}-(l-1)p^{r-2}$, which is less than $|A|=p^r$ for $l\leq p$, contrary to (1).  

We can assume $F$ to be algebraically closed. Let $U_i\subset W_i$ be a one-dimensional $FA$-submodule. Then $K_i$ is the kernel of $U_i$. Set $V'=U_1\oplus \cdots \oplus U_l$. Then $A$ acts faithfully on $V'$ and  $\dim V'=l$. By Lemma \ref{d8d}, $d<l-1$, as claimed.
\enp

 For an arbitrary group  $H$ if $H=K_1\cup ... \cup K_l$ for some subgroups $K_1\ld K_l$ of $H$ then $|H:K_i|\leq m$ for some $i\in\{1\ld l\}$, see   \cite[Theorem  3.3A]{Di}.

\bl{bb8} Let $G\subset GL_n(F)=GL(V)$ be an \ir subgroup, A a minimal non-trivial abelian normal subgroup of G. 
Let W be a quasi\hc of $V|_A$, $U\subset W$ an \ir $FA$-module and K the kernel of $U$. 
Then A  is unisingular \ii $A= \cup_{g\in G}\,\, gKg\up$.\el


\bp Observe that $K$ is the kernel of $W$.  By Clifford's theorem, $G$ permutes quasi\hcs of $V|_A$  transitively, and the kernel of $gW$ is  $gK_1g\up $ for $g\in G$. So the result follows from Lemma \ref{aa8}.
\enp
 
\bl{bc8} Let $G\subset S_n=Sym(\Om)$ be  a transitive subgroup on a set $\Om$ and let A be an abelian normal subgroup of G. Let $\Om=\Om_1\cup \cdots\cup  \Om_l$,
where $\Om_i $ is the sum of all  $A$-orbits   with the same kernel
  $K_i$ for $i=1\ld l$. Then A  has a  fixed point free element  \ii $A\neq (K_1\cup \cdots     \cup K_l)=\cup_{g\in G}\,\, gK_1g\up$. 
\el

\bp It is well known that a transitive abelian subgroup $X$ of $S_k$ has order $k$ and
every $1\neq x\in X$ is  fixed point free. Therefore,   $a\in A$ fixes a point on $\Om_i$ \ii $a\in K_i$. This implies the lemma. \enp

In contrast to Lemma \ref{bb8} $K_1 $ is not necessarily of index $p$ in $A$. 
\med


Let $A$ be an elementary  abelian $ p$-group and let $H\subset {\rm Aut}\,(A)$ be a  subgroup. In many situations below it is convenient to view $A$ as an $\FF_pH$-module 
and   $H$ as a subgroup of  the general linear group $GL(A)$. This allows one to use linear group terminology to express some properties of the action of $H$ on $A$ in a  more friendly fashion. On this way we simply write $H\subset GL(A)$. Note that subgroups of $A$ are interpreted as subspaces of the vector space in question, and $H$-invariant subgroups as $H$-submodules of $A$.

In the \f lemma we identify an elementary abelian $p$-group $A$ with the additive group of a vector space over $\FF_p$.

\begin{lemma}\label{gg2}  Let $r$ be a prime, $G=A\rtimes H$, a semidirect product of an elementary abelian r-group A and a  group $H$ such that $C_H(A)=1$.  Let F be an \acf of characteristic $\ell$ with $(\ell,r)=1$, let  $\lam:A\ra F$ be a non-trivial \rep of A and  $K={\rm ker}\, \lam$. Let $H_\lam:=C_H(\lam):=\{h\in H: \lam(hah\up)=\lam(a)$ for all $a\in A\}$,  and let $m=|H:H_\lam|$. Let $\mu$ be a one-dimensional \rep of $AH_\lam$ defined by $\mu(H_\lam)=1$ and $\mu|_A=\lam$. 
 Then the induced \rep $\mu^G$ is \ir (of dimension $m$) and the \f are equivalent:

$(1)$   $\mu^G$ is unisingular;

$(2)$ A is the union of the conjugates of K.

$(3)$  Let $K_0=\cap_{h\in H} hKh\up$ and $n=|H:N_H(K)|$. Then  $G/K_0$ is isomorphic to a transitive subgroup of the symmetric group $S_{r n}={\rm Sym}(\Delta)$ such that every A-orbit on $\Delta$ 
is of size r and every $a\in A$ fixes a point on $\Delta$. 
\el

\bp Observe that  $C_H(A)=1$ means that $H$ is isomorphic to a subgroup of ${\rm Aut} A$, so we can identify $H$ with a subgroup of $GL_d(\FF_p)$ and $A$ with $\FF_r^d$,
where $d$ is the rank of $A$. As $H\cap A=1$, we have $A\cap H_\lam=1$ and $[A,H_\lam]\subseteq K$. Therefore,  $AH_\lam/K\cong C_r \times H_\lam$, and hence $\mu$ is well defined. 
As  $\mu$ is one-dimensional, $\mu^G$ is monomial. Let $V$ be the module afforded $\mu^G$; to simplify notation we write $gv $ in place $\mu^G(g)v$ for $g\in G, v\in V$. So $V=V_1\oplus \cdots \oplus V_m $ is a direct sum of one-dimensional subspaces permuted by $H$. Let $0\neq v\in V$ be such that $xv=\mu(x)v$ for every $x\in AH_\lam$. Fix some representatives $h_1=1, h_2\ld h_m$ of the cosets $AH/AH_\lam$.
 Set $v_i=h_iv$; then $B=\{v_1\ld v_m\}$ is a basis of $V$ and we can assume that $v_i\in V_i$. Let  $\Omega$ be the set of lines $V_1=Fv_1,\ldots,  V_m=Fv_m$. 
Note that $AH_\lam$ is the stabilizer of $V_1$ in $G$ (as $v=v_1$ and $\mu(H_\lam)=1$). By \cite[Theorem 45.5]{CR},  $\mu^G$ is \ir if  \reps $\mu^g:A\ra F$, $a\ra gag\up, a\in A$ are not equivalent to $\mu$ for every $g\notin H_\lam$.  If $g\notin N_{H}(K)$ then this is the case as the kernel of $\mu^g$
is distinct from $K$. If $g\in N_{H}(K)$ then both $\mu$ and $\mu^g$
 are trivial on $K$, and hence can be viewed as \reps of a cyclic group $A/K$. Then 
 $g\in H_\lam$ \ii $\lam^g=\lam$, and hence $gxg\up=x$ for $x\in A/K$. So
 $N_{H}(K)/H_\lam$ acts faithfully on $A/K$, and we conclude that $\mu$ and $\mu^g$ are distinct, hence non-equivalent representations.

Let $V=W_1\oplus \cdots \oplus W_l$, where $W_1\ld W_l$ are quasi-homogeneous components of $A$ on $V$. Then $G$ transitively permutes them. Let $K_i$ be the kernel of $W_i$ for $i=1\ld l$. Let $a\in A$. Then $a$ is either trivial on $W_i$
or fixed point free. Therefore, $V=V^a\oplus V'$, where $V'$ is the sum of 
$W_i$'s with fixed point free action of $a$, and $V^a$ is the sum of  
$W_i$'s with trivial action of $a$. Therefore, $C_G(a)$ stabilizes both $V^a,V'$. 
\med

$(1)\implies (2)$ As $\mu^G$ is unisingular then so is $V|_A$. So the implication follows from Lemma \ref{bb8}.  
 
\med
$(2)\implies (1)$. Observe that elements $\mu^G(h)$ for $h\in H$ are permutational matrices and hence have \ei 1. By Lemma \ref{bb8}, all elements of  $\mu^G(A)$ have \ei 1.
So we have to deal with elements  $g=ah$, where $1\neq a\in A$ and $h\in H$. Moreover, we can assume  that $g$ is not $r'$-element as otherwise $g$ is conjugate to an element of $H$. We first prove a few auxiliary facts. 

\med
(i) Let $a\in A$  and set $V^a=C_V(a)$. Then $B^a:=B\cap V^a$ is a basis of $V^a$.

\med
Indeed, the matrix of $a$ in basis $B$ is diagonal, so $a(\sum f_ib_i)=\sum f_ib_i$ ($0\neq f_i\in F$) implies $ab_i=b_i$ for every $i$. Whence the claim.

\med
(ii) Let $g\in G. $ If $gv_i=ev_i$ for some $e\in F$ and  $0\neq v_i\in W_i$ for some  $i\in\{1\ld m\}$ then   $e^r=1$.

\med
Indeed, $gv_i=ev_i$ implies $h_i\up gh_i  v_1=h_i\up g v_i=h_i\up ev_i= ev_1$, whence $h_i\up gh_i\in AH_\lam$ by the above. Let $h_i\up gh_i=xy$ with $x\in A, y\in H_\lam$. As $yv_1=v_1$, we have $h_i\up gh_iv_1=xyv_1=xv_1=ev_1$, whence $e^r=1$ as $x^r=1$.

\med
 (iii) Let $g\in G$ be an $r$-element. If $gV_i=V_i$ for some $i\in\{1\ld m\}$ then
 $g\in A$.

\med
 Indeed, $gV_i=V_i$ implies $gh_iV_1=h_iV_1$ and $h_i\up gh_iV_1=V_1$. This implies $h_i\up gh_i\in AH_\lam $. As $A$ is normal in $H$ and  $AH_\lam /A$ is an $r' $-group, we have $g\in A$.
 
\med
 (iv) Let $g\in G$ be an $r$-element.  Then either $\lan g\ran\cap A=1$
 and all $g$-orbits on $\Omega$ are of size $|g|$, or $1\neq g^{d/r}=a\in A$ and
 all $g$-orbits on $\Omega$ are of size $|g|/r$.

 Indeed, let $O $ be a $g$-orbit on $\Omega$ and $\lan g_0\ran$ the stabilizer of some $V_i$ in $\lan g\ran$. By (iii), $g_0\in A $. If $g_0=1$ then $|O|=|g|$, otherwise $|O|=|g|/r$. The condition $g_0=1$ does not depend on the choice of $O$, whence the claim.

\med
Finally we show that $(2)\implies (1)$. By the above, we can assume that $g\notin A\cup H$.

 By (i),  $B^a:=B\cap V^a$ is a basis of $V^a$ (if $a=1$ then $V^a=V$.) Let $\Omega^a=\{V_i:V_i\subset V^a \}$.  As $gV^a=V^a$, $g$ permutes the elements of $\Omega^a$. Let $O\subseteq \Omega^a$. By (iv), $|O|=kt$ and hence $g^{kt}V_i=V_i$ for $V_i\in O$. Let $0\neq v\in V_i$. Then $g^{kt}v=g_1^{kt}v=a^l v=v$. Therefore $g$ fixes the non-zero vector $v+gv+g^2v+\cdots+ g^{kt-1}v$, as required.

\med
 $(2)\iff (3)$. We  first show that the group $G/K_0$ is isomorphic to a subgroup of $ S_r\wr S_n$  (the latter group is viewed as a subgroup of  $S_{rn}$) such that every $A$-orbit is of size $r$.  Let $X=A\cdot N_H(K)=N_G(K)$ and $Y=X/(K\cdot H_\lam)$. Then $Y$  is isomorphic to a subgroup of $S_{r}$. Indeed, 
if $x\in N_H(K)$ then $x$ acts on $A/K$. If this action is trivial then on $A/K$
then $\lam^x=\lam$ (as $\lam^x(a)=\lam(xax\up)$ for $a\in A$), and hence 
$a\in H_\lam$. So $X/Y=(A/K)\cdot (N_H(K)/H_\lam)$ is isomorphic to  a subgroup
$C_r\cdot L$ of $S_r$ with $L\subset \Out(C_r)$, $L\cong N_H(K)/H_\lam$. 
Then $G$ is isomorphic to a subgroup $R$, say, of $S_r\wr S_n$, where $n=G/X$, and every $A$-orbit is of size $r$. The group $R$ is transitive (this is well known and easy to show).

Let $T=T_1\times \cdots\times T_n$, where $T_1\cong T_2\cong \cdots \cong T_n\cong C_r $,  be an abelian subgroup of $ S_r\wr S_n$ of rank $n$. Let $\Delta=\Delta_{rn}$ be the set on which $S_{rn}$ naturally acts and let $\Delta_1\ld \Delta_n$ be subsets, each of size $r$, such that $T_i$ acts trivially on $T_j$ for every $j\neq i$ and $i=1\ld n$. Then the action of $G$ on $T$ by conjugation permutes $T_1\ld T_n$ transitively.
Let $\nu:A\ra T_1$ be a surjective \ho and $K$  the kernel of $\nu$.  Then the mapping  $\nu^g:A\ra gT_1g\up$ $(g\in G)$ defined by $ a\ra gT_1g\up$ for $a\in A$ is surjective. Therefore, $a\in A$ fixes a point on $\Delta_i$ \ii $a\in {\rm ker}\,\, \nu_i$ for $1\leq i\leq n$, equivalently, $a\in gKg\up$.
So $a$ fixes a point  on $\Delta$ \ii $a\in \cup_{g\in G}\,\, gKg\up$. \enp

Observe that $K_0$ is the kernel of the \rep $\mu^G$ in Lemma \ref{gg2}. One can view $\lam$ as an element of the group $A^*:={\rm Hom}(A,F^\times)$, which is isomorphic to $A$, and the action of $H$ on $A^*$ defines on it a structure of $\FF_pH$-module
dual to  $A$, when $A$ is viewed as an $\FF_pH$-module. The elements $\lam^h: h\in H$
form  an $H$-orbit on $A^*$, and the representations $\mu^G$ are parameterized by 
the $H$-orbits on $A^*$. The condition $K_0=1$ is equivalent to saying that the set
$\{\lam^h: h\in H\}$ spans $A^*$ (over $\FF_p$).

\med
 Remarks. 
(1) The character $\lam$ of $A$ can be viewed as an element of the dual 
$\FF_pH$-module $A^*$, and $\dim \mu^G$ equals the size of the orbit $H\lam$
on $A^*$. The \rep $ \mu^G$ is faithful \ii the orbit contains a basis of $A^*$
as a vector space over $\FF_r$. However, the condition that $A=\cup_{h\in H}hK $
cannot be expressed in terms of the orbit $H\lam$. 

(2) Let $G$ be as in Lemma \ref{gg2}. Suppose that $G$ is a transitive subgroup of 
$S_{r^mn}=\Sym(\Delta)$ such that every orbit of $A$ is of size $r^m$. Let $K$ be the kernel of some $A$-orbit. If every $a\in A$ fixes a point on $\Delta$ then $A=\cup _{h\in H}hKh\up$. Equivalently, viewing $A$ as a vector space $V$ over $\FF_r$ and $K$ as a subspace $W$ of codimension $m$, we have $V=HW$.

 \bp[Proof of Lemma {\rm \ref{eq9}}] Let $G=G_{r,p}$. In view of Lemma \ref{na3}, it suffices to prove the lemma for $F$ algebraically closed. Let $\phi:G\ra GL_n(F)$ be a faithful \irr of $G$. Let $H=G/A$ and let  $V$ be an $\FF_r H$-module arising from viewing $A$ as a vector space over $\FF_r$ and the action of $H$ on $V$ obtained from the conjugation action of $G$ on $A$. 

$(1)\ra (2)$  This is a special case of Lemma \ref{gg2}.

$(2)\ra (1)$  Let $\Om_1\ld \Om_k$ be the $A$-orbits on $\Om$.  As $\Om$ is $p$-restricted, we have $|\Om_1|=\cdots =|\Om_k|=p$ and $kp=|\Om|$. Let $K$ be a kernel of $\Om_1$. Then $|A/K|=r$. Therefore, there is a faithfil \rep  $A/K\ra F^\times$. We extends this to a \rep $\lam:A\ra F^\times$ with $\lam(K)=1$.
As $|H|=p$, we have $C_H(\lam)=1$, so $H_\lam=1$ and  $\mu=\lam$ in notation of Lemma \ref{gg2}. By Lemma \ref{gg2}, $\lam^G$ is unisingular.  \enp 

\bl{dc5} Let $W,V'\subset V$ be subspaces and $G\subset GL(V)$ a subgroup. Suppose that $GV'=V'$. If $GW=V$ then $G(W\cap V')=V'$. \el

\bp Let $v\in V'$. By assumption, $v=gw$ for some $g\in G$, $w\in W$. Then $w=g\up v\in W\cap V'$, so $v\in G(W\cap V')$, as required.\enp

 \bl{ny6} Let $n=p^b$, $b\geq 1$, where $p>2$ is a prime, and let $G \subset GL_n (r)=GL(V)$ be a cyclic p-group with $p|(r-1)$. 
Then $V\neq GW$ for every proper subspace $W$ of V.\el 

\bp It suffices to prove the lemma in the case where  $\dim W=n-1$. 
By Lemma  \ref{dc5},  we can assume that $G$ is \irt (Indeed, let $V'$ be a $G$-stable subspace of $V$. If $V'\cap W\neq V'$ then we can apply Lemma  \ref{dc5} to the pair $V', W\cap V'$. If $V'\subset W$ then we apply  Lemma  \ref{dc5} to the pair $V/V', W/ V'$.)

Observe that $G$ is conjugate in $GL_n(r)$ to a subgroup of the group of mononial matrices and $n=p^{b-c}$, where $p^c$ is maximal $p$-power dividing $r-1$. 
As $G=\lan g\ran$ is irreducible, $z:=g^{p^{b-c}}$ is a scalar matrix in $GL_n (r)$, and hence $zW=W $ for every subspace $W$ of $V$.  Set $A=V^+$ and $M=\lan z,A\ran$.

Let $H=AG$ be  a semidirect product. 
Then $M$ is a normal subgroup of $H$ of index $n$.

Suppose that contrary, that $GW=V$.  Let $\lam:V^+\ra \CC$ be a \rep of $V^+$ with kernel $W^+$. 
Then the induced \rep $\lam^G$ is of degree $|G|=r^b$, and $\lam^G=(\lam^M)^G$. By Lemma \ref{gg2}, $\lam^G$ is unisingular. Then the rank of $A$ equals $n=p^{b-c}$. Observe that \ir constituents  $(\lam^M)|_A$
have the same kernel. Let $\mu_1\ld \mu_l$ be   \ir constituents  of $(\lam^G)|_A$ whose kernels are pairwise distinct, and $l$ is maximal with this property. 
Then $l\leq |G/M|=n$ and   $\mu_1+\cdots +  \mu_l$ is a faithful unisingular \rep of $A$ of degree $l$. The rank  of $A$ equals $n$; so $n\leq l-2$ by  Lemma \ref{aa8}(2). Then $l\leq n\leq l-2$, a contradiction.\enp

Note that the assumption $r>2$ in Lemma \ref{ny6} cannot be dropped at least for $n=2$.
Indeed, if $G\subset GL_2(3)$ is the quaternion group or the cyclic group of order 8 then $G$ is transitive on the lines of $V=\FF_3^2$, so $GW=V$ for every proper subspace $W\neq 0$ of $V$. 

\med
We mention here the following fact: 

\bl{kf1} Let $G\subset \FF^\times_{q^2} $, q odd, be a  subgroup of order $q+1$, and let $W\subset V$ be a subspace of dimension $1$. Then $G\FF^\times_q$ is a subgroup of index $2$ in $\FF^\times_{q^2} $. \el

\bp Note that $G\FF^\times_q$ is a subgroup and $|G\cap \FF^\times_{q}| =2$. Hence
the order of $G  \FF^\times_{q}$ is $(q^2-1)/2$.  \enp

Suppose that $q>p$ is a   $p$-power, so $\FF_q$ is contained in an $\FF_p$-subspace $W$ of $V=\FF_{q^2}$ of codimension 1. Is it always true that $GW=V$? 
 Lemma \ref{ob1} hints   that this is true if $q=81$, where $|G/(G\cap \FF^\times_{81})|=41$ is a prime.   

\bl{tv2} Let $G$ be the  derived subgroup of $H=O^-_{2n}(q)$, q odd, and let $V=\FF_q^{2n}$ be the natural module for H. Then the $G$-orbits on V are $H$-orbits. Consequently, $GW=V$ whenever W is a subspace of dimension at least $n+2$ or a non-degenerate subspace of dimension at least $3$.
\el

\bp It is well known that $|H:G|=4$. Let $0\neq v\in V$ and $(v,v)=a$. By Witt's theorem, for every fixed $a\in \FF_q$, the $Hv=\{0\neq x\in V: (x,x)=a\}$.

Suppose that   $a\neq 0$. Then $gv=v$ implies $gv^\perp=v^\perp$,
and $\dim v^\perp$ is a non-degenerate space of dimension $2n-1$.
\itf that $C_H(v)\cong H_1=O_{2n-1}(q)$. Observe that $G$ consists of elements $h\in H$ whose spinor norm equals 1. If $g\in C_H(v)\cong H_1$ then the spinor norms of $g$ as an element of $H$ and $H_1$ coincide. So $G_1=C_G(v)$ is isomorphic to the subgroup of $H_1$ formed by elements of spinor norm 1 in $H_1$. So $|H|/|H_1|=4|G|/4|G_1|=|G|/|G_1|$, as required. 

Let $a=0$, that is, $v$ is totally isotropic. Then $v^\perp=\lan v \ran+V_1$, where
$V_1$ is a non-degenerate subspace of $V$, and $\dim V_1=2n-2$. Let $H_1=C_H(v)$. 
Then $H_1$ contains a subgroup $H_2$ isomorphic to $O^-_{2n-2}(q)$, and hence 
$H_2$ contains elements of arbitrary spinor norm. So for $h\in H$ there exists  $x\in H_1$  such  the spinor norm of $x\up$ equals the spinor norm of $h$. Then the spinor norm of $hx$ equals 1, so $hx\in G$ and $hv=hxv$. 

Let $W$  be a subspace of $V$. If $W$ is non-degenerate and $\dim W\geq 3$ then $W$ contains a non-zero vector of spinor norm $a$ for arbitrary $a\in F$. If $W$ is not non-degenerate, it contains a non-zero vector of norm 0 (that is, singular), and if $\dim W\geq n+2$ then $W$ contains a non-degenerate subspace of dimension 2, which in turn contains a non-zero vector of every non-zero spinor norm.\enp

\bl{nd3} Let $G\cong PSL_2(9)\cong \mathcal{A}_6$ be an \ir subgroup of $GL_4(3)=GL(V)$. Then G preserves a non-degenerate quadratic form  of Witt index $1$ on V, and if  $W\subset V$  is a non-degenerate subspace of dimension $3$ then $V=\cup_{g\in G}(gW)$.\el

\bp The   follows from Lemma \ref{tv2} as
$A_6$ coincides with the derived subgroup of  $ O^-_4(3)$, see \cite[p. 4]{Atlas}. \enp


\section{Unisingular representations of simple groups $L_2(q)$} The \f lemma sorts our the case with  $q$ is even:

\bl{sl2}   Let $G=L_2(q)$  with q even, and let $\phi $ be a non-trivial  \irr of G over a field F of characteristic $2$. Then  $\phi$ is not  unisingular. If F is algebraically closed then $\dim\phi$ is a $2$-power.  \el

\bp See {\rm \cite[Lemma 3.7]{CZ}}. The statement on dimensions  easily follows from Steinberg tensor product theorem \cite[Theorem 41]{St}, which implies that $\phi$ is a tensor product of \ir \reps of $G$, each of dimension 2. \enp 

We  assume until the end of this section that $q$ is odd.

Let $H=SL_2(q)$ and $G=L_2(q)$. It is known that every \ir   Brauer character of $H$ is liftable, that is, obtained from an ordinary character of $H$ by reduction of its values modulo 2 (see for instance \cite[Lemma 4.1]{DZ2}). 
If $\phi$ is an absolutely \irr of $H$ of degree $d$ then $\phi(H)$ is conjugate to $Sp_d(2)$  \ii the Brauer character  values of $\phi$ are integers, see  Lemma  \ref{tt1}.

\bl{20y} Let 
 $G=L_2(q)$, $q=r^a$, $q>3$ odd, r a prime. Let $\phi$ be a
non-trivial $2$-modular absolutely \ir \reps of $G$ and $d=\dim \phi$.

$(1)$  
$d\in\{(q-1)/2, q-1,q+1\}$. In addition, for each such $d$ there exists a $2$-modular absolutely \ir \reps of dimension $d$  unless  $q-1$   is a $2$-power and $d=q+1$.

$(2)$ $\phi$ is unisingular \ii  one of the \f holds:

$\,\,\,\,\,(i)$ $d=q+1;$ 

$\,\,\,\,\,(ii)$ $d=q-1$ and $q>r$; 

$\,\,\,\,\,(iii)$  $d=(q-1)/2$, $q>r^2$ and  $4|(q+1)$.
\el

\bp For (1) see  \cite[p. 31]{HM1}. (2)  Let $g\in G$ with $|g|$ odd. 
Then $|g|$ divides either $(q-1)/2$ or $(q+1)/2$ or $|g|=r$.  
Note that $\phi$ lifts to a \rep $\phi_1$, say, of $SL_2(q)$ over the complex numbers (see \cite[Lemma 4.1]{DZ2} or elsewhere). As $|g|$ is odd, there is $g_1\in SL_2(q)$ such that $|g_1|=|g|$ and $g$ is the image of $g_1$ under the mapping $SL_2(q)\ra SL_2(q)/Z(SL_2(q))$. In addition,  $\phi(g)$ is unisingular \ii so is $\phi_1(g_1)$.

(i) follows from \cite[Lemma 4.2]{DZ2}. Indeed, 1 is an \ei of $\phi(g)$ if $|g|\neq r$
by \cite[Lemma 4.2(4)]{DZ2}, and for $g$ of order $r$ one can use the character 
of $\phi_1$.

(ii), (iii).  Let $d\leq q-1$. 
If $|g|$ divides $q-1$ then 1 is an \ei of $\phi(g)$  by \cite[Lemma 4.2]{DZ2}.

 Suppose first that $g^{q+1}=1$. Using the character table of $SL_2(q)$ one concludes that
 $\phi_1(g_1)$ has \ei 1 if $d=q-1$ or $d=(q-1)/2>|g|$. The latter holds for every odd order element $g$ with $g^{q+1}=1$ \ii $4|(q+1)$.


 Suppose that $|g|=r$. 
If $q$ is a prime and $d\in\{q-1,(q-1)/2\}$  then   $\phi_1(g_1)$ is  fixed point free \cite{Z87}. Suppose that $q$ is not a prime.   Then $\phi_1(g_1)$ is  fixed point free for some  element $g_1$ of order $r$ \ii $d=(q-1)/2$ and $q=r^2$ \cite{Z87}, see also \cite[Proposition 1.2]{DZ8}.

We conclude that $\phi$ is unisingular \ii $d=q+1$, 
or $d=q-1$, $q>r$ and $4|(q+1)$, or $d=(q-1)/2$, $q>r^2$ and $4|(q+1)$.\enp

Remark. The condition $4|(q+1)$ in item (iii) of Lemma \ref{20y} is equivalent to saying that   $4|(r+1)$ and  $a$ is odd. Indeed, if $a$ even then  
  $4|(q-1)$ so $(q+1,4)=2$.  If $a=2b+1$ is odd then $r^a+1=(r+1)m$ with $m$ odd.  

\bl{22a} Let  $G=L_2(q)$, 
where $q>3$ is odd.

$(1)$ $G$ is isomorphic to a unisingular absolutely \ir  subgroup of  $Sp_{q+1}(2)$ \ii $3|(q-1);$

$(2)$ $L_2(q)$ is isomorphic to a unisingular absolutely \ir   subgroup of  $Sp_{q-1}(2)$ \ii q is not prime and $3|(q+1)$.


$(3)$ $G=PGL_2(q)$ is isomorphic to a unisingular absolutely \ir   subgroup of  $Sp_{q-1}(2)$ \ii  q is not a prime and $4|(q+1)$.
\el

\bp 
Let $\phi$ be an absolutely \irr of $G$ of degree $d>1$ and $\be$  the Brauer character of $\phi$. By Lemma \ref{tt1},  $\phi(G)\subset Sp_d(2)$ (in particular, $d$ is even) \ii  (*) $\be(g)\in\RR$   and $\be(g)\pmod 2\in\FF_2$ for every $g\in G$ of odd order. As $\phi$ lifts to characteristic 0 (as mentioned in the proof of Lemma \ref{20y}), we can assume that $\be$ is an ordinary character of $SL_2(q)$.

Suppose first that $G\cong L_2(q)$. We examine cases $(i), (ii), (iii)$ of Lemma \ref{20y} to verify the condition (*). 
 In fact, $(iii)$ is ruled out as the condition $4|(q+1)$ in $(iii)$ implies $d=(q-1)/2$ to be odd. As $\phi$ is liftable, we can use the character table of $G$. Below $1\neq g\in G$ and $|g|$ is odd. 

In case $(i)$  $d=q+1$. If $|g|$ divides $q+1$ then $\be(g)=0$, and $\be(g)\in \ZZ$ if $|g|=r$.
Let $|g|$ divide $q-1$. Let $m$ be the greatest odd divisor of $q-1$. Then $g\in T$ for a cyclic  subgroup $T\subset G$ of order $m$. Then $\be(g)=\zeta(g)+\zeta(g\up)$, where $\zeta$ is a non-trivial one-dimensional character of $T$. So $\be(g)\in \RR$.
If $|g|=m$ then $\be(g)\pmod 2\in \FF_2$ \ii $\zeta^3=1$.
(Indeed, let $\xi=\zeta(g)\pmod 2$. Then $\xi+\xi\up\in \FF_2$ implies $\xi+\xi\up=1$. So $\xi$ satisfies the equation $x^2+x+1=0$  (in $\overline{\FF}_2$), and hence  $\xi^3=1$ and $\zeta^3=1$.)

As $\zeta\neq 1_T$, this requires $3|(q-1)$. Conversely, if $3|(q-1)$ then there exist $\zeta\neq 1_T=\zeta^3$ and $\be$ such that $\be(g)=\zeta(g)+\zeta(g\up)$ for all $g\in T$, and then $\be(g)\in \ZZ$ for all $g\in T$.

If $|g|$ divides $q$ then $\be(g)=1$ then $\be(g)=1$ by the character table of $SL_2(q)$. 

 In case $(ii)$ $d=q-1$ and  $q$ is not a prime. 
Then $\be(g)\in \ZZ$ if $|g|$ divides $q(q-1)$. Suppose that $|g|$ divide $q+1$. Let $m$ be the greatest odd divisor of $q+1$. Then $G$ has a cyclic subgroup $T$ of order $m$. Then $\be(t)=-\zeta(t)-\zeta(t\up)$ for $1\neq  t\in T$, where $\zeta$ is a non-trivial one-dimensional character of $T$. As above we conclude that
 $\be(t) \pmod 2\in \FF_2$ implies $3|(q+1)$. Conversely, if $3|(q+1)$ then there is $\zeta\neq 1_T=\zeta^3$,
and $\be(t)\in \ZZ$ for $t\in T$. 

 Let $G=PGL_2(q)$ and set $G'=L_2(q)$. Let $\phi:G\ra Sp_{q-1}(2)$ be an absolutely \irr of $G$. Since $|G:G'|=2$, $\phi$ is unisingular \ii so is $\phi|_{G'}$ (Lemma \ref{tf1}).
  As every \irr of $G'$ of degree $d>(q-1)/2$ extends to that of $G$, by the above we are left to consider the case  where $\phi|_{G'}$ is reducible. By Clifford's theorem and Lemma \ref{20y}(1), this implies $d=q-1$. 
As $|g|$ is odd, $g\in G'$. Let $\be$ be the Brauer character of $\phi$, and 
$\be_1,\be_2$   the Brauer characters of the \ir constituents of $\phi|{G'}$. By  the Brauer character table of $G$, it suffices 
are to inspect $g\in G'$ of odd order  coprime to $q(q-1)$, and
then we have $\be_i(g)=-1$ for $i=1,2$ for such $g$.  
  \itf  $1$ is an \ei of $\phi_i(g)$ unless $|g|=(q+1)/2$.
Therefore, $\phi$ is unisingular \ii $(q+1)/2$ is not odd, that is, $4|(q+1)$.
\enp

Remark. Lemma \ref{22a} is not quite useful for showing  that $Sp_{q-1}(2)$ has
an absolutely \ir unisingular subgroup,  as already  $AGL_1(q)$, $q$ is not a prime, is isomorphic to  a unisingular $2$-modular \irr of  $Sp_{q-1}$ (Theorem \ref{af1}).  So this lemma  only serves  for proving that some $Sp_{2m}(2)$ has no \ir unisingular subgroup.
 
\bl{h47}
Let H be a group with normal subgroup $G\cong SL_2(q)$, $q$ odd, 
and let $\phi:H\ra Sp_{2n}(2)$ be a unisingular absolutely \ir \rept Let $d$ be the common dimension of the \ir constituents of $\phi|_G$ viewed as a \rep in $Sp_{2n}(\overline{\FF}_{2})$. 

$(1)$ If q is a prime then $d=q+1;$

$(2)$  If q is not a prime and $d=( q-1)/2$ then $d$ is odd.\el

\bp By Clifford's theorem, all \ir constituents of $\phi|_G$ are of the same degree $d$. 

(1) Suppose the contrary. Then  $d\in\{(q-1)/2, q-1\}$ by Lemma \ref{20y}.  By \cite{Z87}, 
1 is not an \ei of any element $g\in G$ of order $q$ in any \irr of $G$ over $\overline{\FF}_{2}$. Therefore, $\phi(g)$ does not have \ei  1 contrary to the assumption.  

$(2)$ 
If $d/2$ is even then $(q+1)/2$ is odd. In the proof of Lemma \ref{20y} we observe that
1 is not an \ei of any element $g\in G$ of order $(q+1)/2$ in any \irr of $G$ over $\overline{\FF}_{2}$ of degree $d$. Therefore, $\phi(g)$ does not have \ei  1, a contradiction.\enp

\section{Observations on representations of some groups of Lie type} 

In this section we assume readers to be familiar with general \rep theory of simple algebraic groups and finite groups of Lie type. The main reference is \cite{St}.

\med
The \f fact is well known, it was exploited in  \cite{Z90}, \cite{GT} and elsewhere:

\bl{u7u} Let $G$ be a simple agebraic group and $\phi$ a linear \rep of G. If  $\phi$ has weight zero then $\phi$ is unisingular.  \el

\med
\subsection{Some simple irreducible linear groups}

\bl{cc5} Let  G be a finite quasisimple group of Lie type in defining characteristic $2$ and let $\phi:G\ra GL_n(\overline{\FF}_2)$ be an \ir \rept Suppose that  $n\in\{2,3,9,18,27,54,81,162\}$. 

$(1)$ Suppose that G is  tensor-decomposable. Then $(n,G)$ is as in Table $2$. 

$(2)$ Suppose that G is  tensor-indecomposable.  Then 
$(n,G)$ is as in Table $3$,
where q is a $2$-power.\el

\bp If $G$ is a group of Lie type in defining characteristic $2$ then the result follows by inspection in \cite{Lu} and general facts of \rep theory of groups of Lie type. For instance, if $G\cong SL_2(q)$ then $\dim\phi$ is well known to be a 2-power (Lemma \ref{sl2}).
So $n=2$ in this case. \enp
 

\begin{table}[ht]\label{tab2}

Table 2: Tensor-decomposable simple \ir subgroups of $GL_n(\overline{\FF}_2)$,\\ $n|162$

\bigskip

\begin{center}
\noindent
\begin{tabular}{|c|c|}
\hline
9&$L_3(q), q>2$, $PSU_3(q), q>2$\\
\hline
27&$L_3(q), q>4$, $PSU_3(q), q>4$ \\
\hline
81&$L_3(q), q>8$, $PSU_3(q), q>8$, $L_9(q), q>2, PSU_9(q), q>2$\\
\hline
\end{tabular}
\end{center}
\end{table}



\begin{table}[ht]\label{tab3}

\noindent
Table 3: Tensor-indecomposable simple \ir subgroups of $GL_n(\overline{\FF}_2),$\\  $n|162$
\begin{center}
\bigskip
\noindent
{\small
\begin{tabular}{|c|c|}
\hline
$n$&$G$ \\
\hline
2&$L_2(q)$, $q>2$  \\
\hline
3&$L_3(q),~ 3\not|q-1,SU_3(q),~ 3\not|(q+1), q>2$\\
\hline
6&$L_6(q),3\not|(q-1),SU_6(q),3\not|(q+1),Sp_6(q), q>2,L_4(q),SU_4(q),G_2(q)$\\
\hline
9&$L_9(q),3\not|(q-1),SU_9(q),3\not|(q+1)$  \\
\hline
18&$L_{18}(q),3\not|(q-1),SU_{18}(q),3\not|(q+1), Sp_{18}(q), \Omega^\pm_{18}(q)$~ 
\\
\hline
 27&$L_{27}(q),3\not|(q-1),SU_{27}(q),3\not|(q+1),$\\ &$~E_6(q),3\not|(q-1),~{}^2E_6(q),3\not|(q+1)$   \\
\hline
54&$L_{54}(q),3\not|(q-1),SU_{54}(q),3\not|(q+1), Sp_{54}(q), \Omega^\pm_{54}(q)$ 
\\
\hline
81 &$L_{81}(q),3\not|(q-1),SU_{81}(q),3\not|(q+1)$  \\
\hline
162&$  L_{162}(q),3\not|(q-1),SU_{162}(q),3\not|(q+1), Sp_{162}(q), \Omega^\pm_{162}(q)$
\\
\hline
\end{tabular}}
\end{center}
\end{table}

\def\ag{algebraic group }

\bl{cb5}  Let $\GG$ be a simple \ag in   characteristic p and $\phi$ an \irr  of $\GG$
with highest weight $\om$. Let $G\subset \GG$ be a finite group.
Let $\phi_i$ be an \irr of $\GG$ with \hw $p^i\om$,  $\psi_i=\phi_i|_G$ for an integer $i\geq 0$, and let $\be_i$ be the Brauer character of $\psi_i$. Then $\be_i=\gamma(\be_0)$ for some Galois \au of $\QQ(\zeta)/\QQ$ for some root $\zeta$ of unity.\el

\bp It is well known that every $p'$-element  $g\in\GG$ lies in a maximal torus $T$, say, and all maximal tori of $\GG$ are conjugate. The weights of $\phi_i$ are $p^i\mu$ when $\mu $ runs over the weights of $\phi$. The weights are homomorphisms $T\ra F^\times$, so $\mu(g)$ is a $|g|$-root of unity for $g\in T$. In addition, $(p^i\mu)(g)=\mu(g)^{p^i}$. Therefore, $\be_0(g)=\sum_j a_j \xi^{p^j}$ for some primitive $|g|$-root of unity $\xi$, and $\be_i(g)=\sum_j a_j \xi^{p^{i+j}}$,  where  
 $a_j$'s are non-negative integers.  
As $(|g|,p)=1$, the mapping   $\xi\ra \xi^{p^i}$ yields a Galois \au $\gamma, $ say, of $\QQ(\xi)/\QQ$, so  $ \xi^{p^i}=\gamma(\xi)$ and $\be_i(g)=\gamma (\be_0(g))$.\enp

\bl{e619} Let $G\cong E_6(q)$ or ${}^2E_6(q)$, q even,  and let $\phi$ be a \rep of G whose all composition factors are  of degree $27$. Then $\phi$ is not unisingular. \el

\bp Suppose first that $\phi$ is \irt It is well known that $\phi$ extends to a 
 \rep $\Phi$ of the simple algebraic  group $\mathbf{G}$ of type $E_6$. The \hw of $\Phi$ is well known to be $2^i\om_1$ or $2^i\om_6$ for some integer $i\geq 0$.  
  One easily observes that  it suffices to examine   $\Phi$ of \hw $\om_1$. Observe that $\Phi$ is faithful. If $Z(G)\neq 1$ then the statement is trivial, otherwise $(3,q-1)=1$ if
$G\cong E_6(q)$  and $(3,q+1)=1$ if
$G\cong {}^2E_6(q)$. 

  The group $ E_6(\overline{\FF}_2)$ contains a subgroup isomorphic to $X:=SL_6(\overline{\FF}_2)$. Let $\lam_1\ld \lam_5$ are the fundamental weights of $X$.  By \cite[Table 8.7, p.108]{LS},   the composition factors of $V_{\om_1}|_X$ are  $V_{\lam_1}$ (with \mult 2) and $V_{\lam_4}$. Note that $E_6(q)$, respectively, ${}^2E_6(q)$ contains a subgroup $X(q)$ such that $X(q)\cong  SL_6(q) $, respectively, $X(q)\cong SU_6(q)$. (As $Z(G)= 1$, we have $Z(X(q))=1$ too.)  Therefore, it suffices to show that there is an element $g\in X(q)$ acting  fixed point freely on $V_{\lam_1}$ and  $V_{\lam_4}$. Since $V_{\lam_4}$
 is dual to  $V_{\lam_2}$, we may deal with $V_{\lam_2}$ in place of $V_{\lam_4}$.
 Let $g\in  X(q)$ be an element of order $d=(q^6-1)/(q-1)$ or $(q^6-1)/(q+1)$ if $G\cong SL_6(q)$ or $G\cong  SU_6(q)$, respectively. 

The Jordan form of $g\in X$ at $V_{\lam_1}$  is $\diag(\zeta,\zeta^q,\zeta^{q^2},\zeta^{q^3},\zeta^{q^4},\zeta^{q^5})$ if $X(q)\cong SL_6(q)$, and $\diag(\zeta,\zeta^{-q},\zeta^{q^2},\zeta^{-q^3},\zeta^{q^4},\zeta^{-q^5})$ if $X(q)\cong SU_6(q)$, where $\zeta\in F$ is a primitive $d$-root of unity. So $g$ is  fixed point free on $V_{\lam_1}$. As $V_{\lam_2}$ is the exterior square of $V_{\lam_1}$, the \eis of $g$ on $V_{\lam_2}$ are in the set $\{\zeta^{\pm q^i\pm q^j},~i,j\in\{0,1\ld 5, ~i<j\}$. Note that  $\zeta^{\pm q^i\pm q^j}=1 $ \ii $\pm q^i\pm q^j\equiv 0\pmod{d}$, equivalently, $\pm 1\pm q^{j-i}\equiv 0\pmod{d}$.  As $d>q^6-q^5>q^5+1\geq q^{j-i}\pm 1$, this implies  $q^{j-i}\pm 1=0$, which is false. \enp

\bl{c44} Let $G\subset GL_n(\overline{\FF}_2)=GL(V)$ be a group listed in Table $3$  and $\psi,\psi'$ two \ir \reps of G of degree n. 

$(1)$ $\psi'$ is a Galois conjugate of either $\psi$ or the dual of $\psi$. Consequently, $\psi(g)$ has \ei $1$ \ii $\psi'(g)$ has \ei $1$.

 $(2)$  G contains fixed point free elements.\el

\bp 
(1) The result follows from \cite{Lu} and the theory of automorphisms of finite simple groups of Lie type. The additional statement follows from Lemma \ref{ek1}.

(2) For \ir \reps of classical groups in Table 3 with $n$ equal the dimension of their natural
 \rep the claim is well known. For $E_6(q)$ and ${}^2E_6(q)$ see Lemma \ref{e619}.
Let $G=SL_4(q)$ and $SU_4(q)$, $q$ even, and $n=6$. Note that  $H=SL_4(2)\subseteq SL_4(q)$.  By general theory, $\psi(H)$ is \irt Then elements of order 7 in   $SL_4(2)$ and   do not have \ei 1 in an \irr of   $H$ of dimesion 6, see \cite{JLPW}. Let $G\SU_4(q)$. Then
an \irr of degree 6 identify $PSU_4(q)$ with the orthogonal group $\Om_6^-(q)$. 
Element of order $q^3+1$ are well known to be \ir on the natural module of the latter group,
hence do not have \ei 1. \enp

\subsection{The \mult of $1_G$ in the 2-modular reduction of the Weil \reps of $SU_n(q)$, $q$ even} The main purpose of this section is to prove Lemma \ref{tp2} used in the proof 
of  Lemmas \ref{x54} and \ref{na1}. On the way we obtain a more general result (Proposition \ref{pu1}) which we hope to use in future. 

Let $H=U_n(q)$, $n$ odd. We define Weil \rep $\Phi$ of $H$ following \cite{Ge}.

\bl{ww1}  Let $H=U_n(q), G=SU_n(q)$, $n>1$ odd, q even, and let $Z=\lan z\ran=Z(H)$. Let $\Phi$ be the Weil \rep of H and M the $\CC H$-module afforded by $\Phi$. Define $M_s=\{x\in M: zx=\zeta^sx\}$ for $s=0,1\ld q$, where $\zeta$ is a primitive $(q+1)$-root of unity. Let $d=(n, q+1)$  and $D_s(1_G)$ be the \mult of $1_G$ in $M_s\pmod 2$. Then

 $$D_s(1_G)=\begin{cases}
 d-1&if\,\,s=0\,\,;\\
 d&if\,\,s>0 \,\, {\rm and} \,\,d|s;\\
 0&if\,\,s>0\,\, {\rm and} \,\,d\not|s.\\
 \end{cases}$$
\el

\bp We use \cite{z91}, where the computation of the decomposition numbers of $D_s(\phi)$ for every \ir 2-modular \rep  $\phi$ is reduced to
computing the number of solution of a  certain equation in the ring $\ZZ/(q+1)\ZZ$. For $\phi=1_G$
 the equation takes the form
\begin{equation}\label{eq1}n(1+\sum _{a=1 }^lx_a2^{a-1})\equiv s \pmod{q+1},\end{equation}
and the main theorem of \cite{z91} asserts that $D_s(1_G)$ equals the number of solution of (\ref{eq1}) in $x_1\ld x_l\in\{0,1\}$.

Let $e=(q+1)/d$.
Suppose first that $s\not\equiv 0\pmod d$, in particular, $d>1$. Then the congruence has no solution. Indeed, the lhs lies in $d\ZZ$ and $s\notin d\ZZ$. As $d\ZZ\subset (q+1)\ZZ$, we have $s \pmod{q+1}\notin d\ZZ$. (This can be also observed straightforwardly as follows. Let $z_1=z^e$; then  $z_1\cdot \Id\in G$ and $z_1^s\neq 1$. Then $z_1^s\pmod 2\neq 1$. As $z_1$ acts scalarly on $M_s$, it follows that $z_1$ acts scalarly on every composition factor of $M_s\pmod p$. This yields a contradiction.)

Suppose that  $s\equiv 0\pmod d$. Then we have 
\begin{equation}\label{eq2}\frac{n}{d}(1+\sum _{a=1}^lx_a2^{a-1})\equiv \frac{s}{d} \pmod{e}
\end{equation} 

Let $m=\sum _{a=1}^lx_a2^{a-1}$; this can be viewed as the 2-adic expansion of $m$. \itf
the set $\{\sum _{a=1}^lx_a2^{a-1}:x_a=0,1\}$ is in bijection with $\{0\ld q-1\}$, and hence the set $\{1+\sum _{a=1}^lx_a2^{a-1}:x_a=0,1\}$ is in bijection with $\{1\ld q\}$.
These can be viewed as  $\ZZ_{q+1}\setminus 0$, where $\ZZ_{q+1}=\ZZ/(q+1)\ZZ$ is the residue ring modulo $q+1$. As $(\frac{n}{d},q+1)=1$, we have $\frac{n}{d}\{1\ld q\}=\{1\ld q\}$.

\itf the number of solutions of (\ref{eq2}) equals the number of integers $k\in \{1\ld q\}$ congruent to $s/d$ modulo $e$. Let $s'=(s/d)\pmod e$.
 If $s'=0$ then these are $e,2e\ld (d-1)e$. If $s'\neq 0$ then these are $s',s'+e\ld s'+(d-1)e $ as $s'+(d-1)e=q+1-(e-s')\leq q$. So we obtain $d-1$ solutions in the former case and $d$ solutions in the latter case.\enp

\begin{propo}\label{pu1} Let $G=SU_n(q)$, $q>2$ even, n odd, and let $\Phi$ be the Weil \rep of G. Let $\phi\neq 1_G$ be an \ir constituent of $\Phi\pmod 2$.
Let $g\in G$ with $|g|=(q^n+1)/(q+1)$. Then $1$ is not an \ei of $\phi(g)$. In particular, $\phi$ is not unisingular.\end{propo}

\bp Let $V$ be the natural $\FF_{q^2}G$-module. Then $g$ is \ir on $V$. Let $1\neq h\in \lan g\ran$. As $gV^h=V^h$, it follows that  $V^h=0$.  Let $\chi$ be the character  of $\Phi$. Then  $\chi(1)=q^n$ and $\chi(h)=-1$ \cite{Ge}. Therefore, $(\chi|_{\lan g \ran}, 1_{\lan g \ran})=q$, so the \mult of \ei $1$ of $\Phi(g)$ equals $q$.
By Lemma \ref{ww1}, the \mult of $1_G$ in $\Phi\pmod 2$ is $d-1+xd$, where $x$ is the number of integers $s\in\{1\ld q\}$ such that $d|s$. So $x=e-1$, whence $d-1+xd=q$. Therefore, the \mult of \ei 1 of $\Phi\pmod 2$ equals the number of trivial composition factors, and hence 1 is not \ei of  $\phi(g)$.\enp

\bl{rr3} {\rm  \cite[Proposition 1.10 and  Lemma 3.12]{z91}} Let $q=2^l$, let $F$ be an \acf of characteristic $2$ and $\mathbf{G}=SL_n(F)$.  Let $\Phi$ be a Weil \rep of $G=SU_n(q)$ and $\phi$ is an \ir F-\rep of $\mathbf{G}$. Then $\phi$ is a composition factor of $\Phi\pmod 2$ \ii $\phi=\lam|_G$, where $\lam$ is an \irr of  $\mathbf{G}$ with highest weight $\lam_1+2\lam_2+\cdots +2^{l-1}\lam_l$, where $\lam_i\in\{0,\om_1\ld \om_{n-1}\}$ for $i=1\ld l$ and $\om_1\ld \om_{n-1}$ are the fundamental weights of $SL_n(F)$.\el
 
\bl{tp2}  For q even let $G\cong SL_{3}(q)$ or $SU_{3}(q), q>2$. Let
$\phi$ be an \ir $2$-modular \rep of G of degree  
$  3^m$, $m>0$. Let $g\in G$ be an \ir element of  order $q^2+q+1$ or $q^2-q+1$, respectively. Then $1$ is not an \ei of  $\phi(g)$.\el

\bp Let $\mathbf{G}=SL_{3}(F)$ and let $\tau$ be an \irr of  $\mathbf{G}$ such that $\phi=\tau|_G$. Then $\tau=\tau_1\otimes \cdots\otimes \tau_m$,
where $\tau_1 \ld \tau_m$ are non-trivial  \ir \reps of $\mathbf{G} $, each is  of   3-power dimension.  Then each $\tau_i$, $i=1\ld m$, is a Frobenius twist of a 2-restricted \irr  $\mathbf{G} $.  By \cite{Lu}, $\mathbf{G}$ has no 2-restricted \irr  of degree $3^k$ for $k>1$, so each degree equals 3, and those of degree 3 are of  \hw  $\om_1$ or $\om_2$.     \itf the \hw of  $\tau$ is as in Lemma \ref{rr3}. 

If $G=SU_{3}(q)$ then the result follows from Proposition \ref{pu1}.

Let $G=SL_{3}(q)$. Then $\tau$ satisfies the assumption of \cite[Lemma  4.14]{z16}. Observe that  $\tau$ does not have weight 0 (this can be shown straightforwardly and also with use of  Proposition \ref{pu1} and Lemma \ref{rr3}.) By \cite[Lemma  4.14(1)]{z16}, 1  is not an \ei of  $\phi(g)$, unless $\lam$, the \hw of $\tau$, is of shape $(q-1)\om_i$ for $i=1,2$. In the latter case the result follows from \cite[Corollary  4.15(1)]{z16}.\enp

Similarly, we have:

\bl{tp9}  For q even let $G\cong SL_{9}(q)$ or $SU_{9}(q), q>2$. Let
$\phi$ be an \ir $2$-modular \rep of G of degree  
$  3^m$, $0<m<10$. Let $g\in G$ be an \ir element of  order $(q^9-1)/(q-1)$ or $(q^9+1)/(q+1)$, respectively. Then $1$ is not an \ei of  $\phi(g)$.\el

\bp Let $\mathbf{G}=SL_{9}(F)$ and let $\tau$ be an \irr of  $\mathbf{G}$ such that $\phi=\tau|_G$. Then $\tau=\tau_1\otimes \cdots\otimes \tau_m$,
where $\tau_1 \ld \tau_m$ are  \ir \reps of $\mathbf{G} $, each is  of   3-power dimension.  Then each $\tau_i$, $i=1\ld m$, is a Frobenius twist of a 2-restricted \irr  $\mathbf{G} $. 
Let $\om_i$, $i=1\ld 8$, be the fundamental weights of $\mathbf{G}$, and let $\om=\sum a_i\om_i$ be the \hw of a 2-restricted \irr  $\mathbf{G} $, where $a_i\in\{0,1\}$ as  $q$ is a 2-power.  By \cite{Lu}, $\mathbf{G}$ has no 2-restricted \irr  of degree $3^k$ for $2<k<10$, so each degree equals 9, and those of degree 9 are of  \hw  $\om_1$ or $\om_8$.     \itf the \hw of  $\tau$ is as in Lemma \ref{rr3}. 

If $G=SU_{9}(q)$ then the result follows from Proposition \ref{pu1}.

Let $G=SL_{9}(q)$. Then $\Phi$ satisfies the assumption of \cite[Lemma  4.14]{z16}. Observe that  $\tau$ does not have weight 0 (this can be shown straightforwardly and also with use of  Proposition \ref{pu1} and Lemma \ref{rr3}.) By \cite[Lemma  4.14(1)]{z16}, 1  is not an \ei of  $\phi(g)$, unless $\om$, the \hw of $\tau$, is of shape $(q-1)\om_i$ for $i=1,8$. In the latter case the result follows from \cite[Corollary  4.15(1)]{z16}.\enp

\bl{12a} Let $F$ is \acf of characteristic $r>0$,   $\mathbf{G}=SL_2(F)$,  and let $\rho$ be the \irr of $\mathbf{G}$ with \hw $a\om$, where  $0<a<q=r^b$ and $\om$ is the fundamental weight of $\mathbf{G}$.  Let $G=SL_2(q)\subset \mathbf{G}$ and $\phi=\rho|_G$. 

$(1)$ $\phi$ is unisingular \ii $\rho$ has weight $0$.

$(2)$ Let $g\in G$ be   of order $q+1$. Then $\rho(g)$ has \ei  $1$ \ii
$\rho$ has weight $0$.
 
$(3)$ Let $a=a_0+a_1r+a_2r^2+\cdots+ a_{b-1}r^{b-1}$ $(0\leq  a_0\ld a_{r-1}<p)$ be r-adic expansion of a. Then $\rho$ has weight $0$ \ii all numbers $a_0\ld a_{b-1}$ are even. In particular,  if q is even then  $\rho$ has weight $0$ \ii $a=0$.\el

\bp If $\rho$ has weight 0 then $\rho$ is unisingular (Lemma \ref{u7u}). 

(1), (2) It is well known that every semisimple element of $G$ is contained in a maximal torus of $G$. There are two conjugacy classes of maximal tori of $G$, say, $T_1$, $T_2$, each of them is a cyclic group. 
Let $\be$ be the Brauer character of $\rho$. Then $\rho$
is  is unisingular \ii 
$(\be|_{T_i},1_{T_i})>0$, where 
$1_{T_i}$ is trivial character of $T_i$ and $(.,.)$ is the inner product of 
characters. By \cite[Theorem 1.3]{z16}, if $a<q-1$ then
 $(\be|_{T_i},1_{T_i})>0$ \ii $\rho$ has weight 0, whence the result in this case. 
Let $a=q-1$. Then    $(\be|_{T_i},1_{T_i})>0$ if either $\rho$ has weight 0 or  $|T_i|=q-1$, again by \cite[Theorem 1.3]{z16}.  This implies the claims.

(3) Let  $\rho_i$ be the \irr of $\mathbf{G} $ with \hw $a_i\om$.  By  Steinberg's tensor product theorem $\rho=\otimes_{i=0}^{r-1} Fr^i(\rho_i)$, where $Fr$ is the standard Frobenius endomorphism $\mathbf{G} \ra \mathbf{G} $ \cite[Theorem 41]{St}. Then the weights of $\rho$
are $\mu_0+r\mu_1+\cdots +r^{b-1}\mu_{r-1}$, where $\mu_0\ld \mu_{r-1}$  are weights of $\rho_0\ld \rho_{r-1}$, respectively. The weights $\mu_i$ are of the form $ c_i\om$, where $-(r-1)\leq c_i\leq r-1$. \itf $\sum _{i=0}^{r-1} c_ir^i=0$ \ii $c_0=\cdots =c_{r-1}=0$, whence the result.\enp

\begin{corol}\label{1k4} Let $G=SL_2(q)$ and let $\phi:G\ra GL_4(F)$ be an \irr of $G$,
where F is an \acf of characteristic $r|q$.
Suppose that  $\dim\phi=4$. Then $\phi$ is not unisingular.  In addition, if $g\in G$ is   of order $q+1$ then $\phi(g)$ does not have \ei $1.$
\end{corol}

\bp Let $q=r^b$. If $\phi$ is tensor indecomposable then $p>3$ and  $\phi$ is faithful, so $\phi(G)$ contains $-\Id$, whence the result in this case. Suppose that $\phi$ is tensor-decomposable. Then $\phi(Z(G))=1$ and $q\geq p^2$, so $\phi(G)$ is a simple group.  
By Steinberg's theorem \cite[Theorem 43]{St},   $\phi$ extends to a \rep $\rho$ of $\mathbf{G}=SL_2(F)$ with \hw $(r^i+r^j)\om$, where $0<i<j<b$.   So the result follows from Lemma \ref{12a}.\enp

\subsection{The Steinberg characters of unitary groups $SU_{n}(q), n$ odd, $q$ even, and $E_7(q)$, $q$ even}
 
Let $G$ be a simple group of Lie type in defining characteristic $r$. Then $G$ has a unique \irr over the complex numbers of degree $|G|_r$, called the Steinberg \rep of $G$. We denote it by $\si(G)$. It is well known that $\si$ remains \ir under reduction modulo $r$; we denote this by $\si_r$. One easily observes that $\si$ is unisingular \ii
so is $\si_r$. (Indeed, this is trivial for $r'$-elements (=semisimple elements) of $G$.
If $x\in G$ is $r$-singular element then $\chi(x)=0$, where $\chi$ is the character
of $\si$. Let $X=\lan x\ran$ and let $Y$ be the subgroup formed by $r'$-elements of $X$. Then $|X|(\si|_X,1_X)=|X|\cdot \sum_{y\in Y} \si(y)=(|X|/|Y|)\cdot (\si|_Y,1_Y)$.) 

It is shown in \cite[Theorem 3]{Z90} that $\si$ and $\si_r$ are unisingular whenever $r> 2$. This is not   true for   $G=SL_2(q)$ with $q$ even, see \cite[Remark 1]{Z90}. In addition, $\si_2$ is unisingular if $\si_2$, extended to a \rep of a suitable algebraic group $\mathbf{G}$,  has weight 0 (Lemma \ref{u7u}). This holds for $r=2$ and the groups  of type $E_8,E_6,F_4,G_2$, for  $A_n$ for $n$ even, for $C_n$ with $4|n(n+1)$, for $D_n$ with $4|n(n-1)$,
  see the proof of \cite[Theorem 3]{Z90}. 

In a recent paper  \cite{CZ} it was proved that $\si_2$ is unisingular for classical groups of type $ C_n(q),D_n(q),{}^2D_n(q)$ and $A_n(q)$. For $G=E_7(q)$ and $G={}^2A_n(q)$ with $n$ odd the question   remained open.   

  These cases are settled below; the argument is 
based on a recent work by Malle and Robinson \cite{MR}.   

\begin{theo}\label{st5} Let G be a finite simple group of Lie type in defining characteristic $r$, $\si$ the Steinberg \rep of G over the complex numbers and $\si_r$ the reduction of $\si$ modulo $r$. Then $\si$ and $ \si_r$ are unisingular.  
\end{theo}

\bp It is well known that $\si$ is absolutely \irt  One has to show that  $\si(g)$  has \ei 1 for every semisimple element $g\in G$ (as mentioned above, this is equivalent to saying that $\si_r(g)$ has \ei 1.) 
Let $\chi$ be the   character of $\si$.   Recall that $Z(G)$ is in the kernel of $\si$.

\med
Let $X=\lan g \ran$.  We show that $(\chi|_X,1_X)>0$. We have

$$(\chi|_X,1_X)=\frac{1}{|X|}(\chi(1)+\sum _{1\neq x\in X}\chi(x)), $$
which is non-zero \ii $\chi(1)+\sum _{1\neq x\in X}\chi(x))>0$. We have 

$$\chi(1)+\sum _{1\neq x\in X}\chi(x)\geq \chi(1)-\sum _{1\neq x\in X}|\chi(x)|\geq  |G|_r-(|X|-1)\cdot\max |\chi(x)|,$$
where $|\chi(x)|$ means the absolute value of a complex number $\chi(x)$. By general theory, $|\chi(x)|=|C_G(x)|_r$, in our case $r=2$. 

In \cite[Proof of Proposition 6.1]{MR}, the authors show that 
$|C_G(x)|_r\leq \si(1)/q^m$, where $m$ is given in \cite[Table 4]{MR}. Specifically, $m=2n-2$ for ${}^2A_n(q)$ with $n>1$ and $m=32$ for $G$ of type $E_7(q)$.
So for $G={}^2A_n(q)$ we have  $|\chi(x)|\leq |G|_r/q^{2n-2}$ for $1\neq x\in \lan g \ran $ and
 $$|G|_r-(|X|-1)\cdot\max |\chi(x)|\geq |G|_r-(|X|-1)\cdot |G|_r/q^{2n-2}=|G|_r(1-\frac{|X|-1}{q^{2n-2}}>0$$

\noindent
if $|X|-1< q^{2n-2}$. One easily observes that, for $n$ even, the maximum order of $g$ does not exceed $q^n+1$,  
so $|X|< q^{2n-2}$ if  $q^{n}+1<q^{2n-2}$, equivalently, $q^n(q^n-2)>1$. This is  true for $n>2$.

\medskip
Let $G=E_7(q)$. Then $|\chi(x)|\leq |G|_r/q^{32}=q^{31}$ for $1\neq x\in \lan g \ran $, see \cite[Proof of Proposition 6.1]{MR}. In addition, $|x|\leq (q^5 -1)(q^2 +q+1)<q^8$, see for instance \cite[p. 898]{DF}, so $|\chi(x)|\cdot |X|< q^{39}<|G|_r$, and  the result follows. \enp

Observe that the above argument does not work for $r$-modular \ir \reps $\phi$ of $G$ other than $\si_r$, whereas the approach used in \cite{CZ} provides a sufficient condition for
for $\phi$ to be unisingular.

\section{Examples of irreducible unisingular subgroups}

In this section we provide examples of absolutely \ir unisingular subgroup of $Sp_{2n}(2)$
 for $n\leq 125$.

\subsection{Sources of examples}

Corollary \ref{in5} yields examples for every $n=4k$, totally 31 examples for $n<125$, and additionally 
9 examples for $n=6k$ and 11 examples for $n=7k$ due to Lemma \ref{c3a4} for $n=6$ and Lemma \ref{22a} with $q=13$ for $n=7$. In fact,  Lemma \ref{22a} gives another large group of examples for $n=(p+1)/2$ with $p>2$ a prime, including $n\in\{7,10,19,22,31,34,37,61,79,93,106,115\}$
and their multiples. Lemma \ref{22a} and Theorem \ref{af1} lead to a few other example $n=(q\pm 1)/2$ with $q$ is not  a prime.

Further examples arise from an \ir 2-modular \rep of simple groups $G$ over $\FF_2$. There are two cases to be differed: $G$ is isomorphic to a simple group of Lie type in defining characteristic 2 and remaining simple groups.

Suppose first that $G$ is a simple group not isomorphic to a group of Lie type in defining characteristic 2. The degrees $d$ of \ir \reps  of these groups are listed for $d<250$ in \cite[Table 2]{HM}, with omission of $L_2(q)$ and the alternating groups $A_{d+1}$ and $A_{d+2}$; the latter cases
can be ignored for our purpose due to Lemmas \ref{22a} and \ref{aa3}. Moreover, \cite[Table 2]{HM} indicates the minimal field of realization of a \rep in question, so we omit the groups  that are not realized over $\FF_2$. Next we have to verify whether the remained \reps  are unisingular. Usually this can be checked with the Brauer character tables in \cite{JLPW}. In particular, 
the entries with $2n\in\{ 78,174,202,218\}$ are worked out on this way, see  Table 4. 

Suppose   that $G$ is a simple group of Lie type in defining characteristic 2. To obtain \ir \reps over with $\FF_2$
we restrict ourselves with Chevalley groups  over $\FF_2$    such as $L_{m}(2)$, 
$Sp_{2m}(2)$, $\Omega^+_{2m}(2)$, $E_6(2)$, 
$E_7(2)$, $E_8(2)$, $F_4(2)$,   $G_2(2)$. The \ir \reps of these groups can be obtained from \ir \reps of the corresponding algebraic
group $\mathbf{G}$. More precisely, if $\mathbf{G}$ is of rank $r$ then
the \ir \reps of $\mathbf{G}$ are parameterized by integral vectors $\om=(a_1\ld a_r)$, with non-negative $a_1\ld a_r$, and those
of $G$ are parameterized by such vectors under condition $ a_1\ld a_r\in\{0,1\}$. In fact, under this condition  an \irr
 $\mathbf{G}$ remains \ir on  $G$, and all \ir \reps of $G$ are obtained in this way.
 If $\phi$ is an \irr of $\mathbf{G}$ corresponding to a string $(a_1\ld a_r)$ with entries $a_1\ld a_r\in\{0,1\}$ is of degree $n$, say, then $\phi(G)$ is equivalent to a \rep into $GL_n(2)$. If $\phi$ is self-dual then  $n$ is even and $\phi(G)$ is equivalent to a \rep into $Sp_n(2)$.

 The \ir \reps of $\mathbf{G}$ of relatively small degrees  are listed in \cite{Lu} (in term of the strings
$(a_1\ld a_r)$).  Our strategy is to extract from the list in \cite{Lu} the groups $\mathbf{G}$ in defining characteristic 2, and exclude those for which $\phi|_G$ are not unisingular. 

This method is used for $2n\in\{34,118,132,142,188,194,230,246\}$. (Note that we cannot ignore
other groups of Lie type in defining characteristic 2 for proving that $Sp_{2n}(2)$ has no unisingular \ir subgroup for certain $n$.) The Lemma \ref{ad6} illustrates the use of the method for special values of $n$.

\bl{ad6}  
Let $\mathbf{G}$ be a simple simply connected  algebraic group of rank $>1$ in   characteristic $2$.
Let $\phi$ be the non-trivial \ir constituent of the adjoint \rep of $\mathbf{G}$ and $n=\dim \phi$.  

$(1)$  $\phi(\mathbf{G})\subset Sp_{n}(F)$ is unisingular. In addition, $\phi(G)$ is unisingular for every finite subgroup G of $\mathbf{G}$.

$(2)$ If $G$ is a finite simple group of Lie type over $\FF_2$ of the same type as $\mathbf{G}$ then $\phi(G)$ is \ir and $\phi(G)\subset Sp_n(2)$. 

$(3)$ If $\mathbf{G}=SL_d(\overline{\FF}_2)$ then $\dim\phi=d^2-2$ for d  even and $d^2-1$ for d odd.

$(4)$ If $\mathbf{G}=SO_{2d}(\overline{\FF}_2)$ then $\dim\phi=2d^2-d-2 $ for d  even and $2d^2-d-1$ for d odd. 
\el

\bp It is well known that the weight 0 \mult  in the adjoint \rep of $\mathbf{G}$ is at least 1,
except for $\mathbf{G}\cong SL_2(\overline{\FF}_2)$. (Moreover, the \mult in question is listed in \cite{TeZ}  for instance). Whence the result.\enp

\begin{corol}\label{22y} $(1)$ The groups $SL_d(2)$ with $3\leq d\leq 15$ has  a unisingular  absolutely \ir \rep in $Sp_n(2)$ for $n$ in Table $4.$  

$(2)$ The groups $E_7(2)$   has  a unisingular  absolutely \ir
\rep in $Sp_{132}(2)$.

$(3)$ The simple group   $\Om^+_{10}(2)$    has a unisingular   absolutely \ir
\reps in $Sp_{118}(2)$. 

$(4)$ The simple group  $\Om ^+_{20}(2)$ has a unisingular   absolutely \ir
\reps in $Sp_{188}(2)$ of degree $188$.
\end{corol}

\bp For (2) see for instance \cite{Lu}. Recall that simple group $\Om^+_{2n}(2)$, $n\geq 4$,   is a  subgroup of index 2 in $O^+_{2n}(2)$.
\enp

\bigskip
\begin{center}

Table 4: Degrees $n$ of some unisingular \reps of $SL_d(2)$ for $d\leq 15$

\bigskip
\begin{tabular}{|c|c|c|c|c|c|c|c|c|c|c|c|c|c|}
\hline
$d$&$3$&4&5&6&7&8&9&10&11&12 &13&14&15\\
\hline
$n$&$8$&14&24&34&48&62&80&98 &120&142&168&194&224  \\
\hline
\end{tabular}
\end{center}
\bigskip

\subsection{The cases $2n=12,30$}

\bl{c3a4} Let $G=C_3^3{\rtimes}{\mathcal A}_4$ with faithful  conjugation action of the alternating group ${\mathcal A}_4$ on $C_3^3$. Then G is isomorphic to a unisingular absolutely \ir subgroup of  $Sp_{12}(2)$.\el

\bp   In Table 
6 is given the character table of $G$, which also shows that $G$ has 11 conjugacy classes of odd order elements. Therefore, $G$ has 11 \ir Brauer characters. The character $\rho_4$ has $C_3^3$ in its kernel, so this can be viewed as a character of ${\mathcal A}_4$; this is reducible modulo 2 as ${\mathcal A}_4$ has a non-trivial normal 2-subgroup.

In addition, the characters  $\rho_{11},\rho_{12}$ coincide modulo 2. \itf that all
characters except $\rho_4$ are \ir modulo 2.

Next we observe that $\rho_{13} $ is the character of a unisingular \rept For this we show that $(\rho_{13}|_X,1_X )>0$ 
for every cyclic group $X$. If $|X|=3$ then $|X|(\rho_{13}|_X,1_X )=3(12+\rho_{13}(x)+\rho_{13}(x^2))>0$
as   $-3\leq \rho_{13}(x)\leq 3$ for every $x\in G$ of order 3. If $|X|=\lan x\ran=9$ then $|X|(\rho_{13}|_X,1_X )=9(12+\rho_{13}(x^3)+\rho_{13}(x^6))$  as $\rho_{13}(x^i)=0$ for $i\in\{1,2,4,5,7,8\}$ and we conclude as above. \itf that $\rho_{13}\pmod 2$ is unisingular. (Similarly, $(\rho_{13}|_X,1_X )>0$ for   $X$ of order 6 and 2.) 

Finally, by Lemma \ref{tt1}, $\rho_{13}(G) \pmod 2$ is contained in a group isomorphic to $Sp_{12}(2)$.\enp
 
Now we consider  the case where $2n=30$.

\bl{31a} Let $G$ be a semidirect product of $V^+=\FF_4^3$ and the alternating group $H\cong \mathcal{A}_6$, such that $V^+\neq Z(G)$. Then there exists an absolutely \irr $G\ra Sp_{30}(2)$. \el

\bp In notation of Lemma \ref{nd3} let $A=V^+=\FF_3^4$ and let $\lam$ be a non-trivial one-dimensional \rep of $V^+$ with kernel  $K=W^+$. By Lemma \ref{nd3}, $HW=V$. Then  $H_\lam$, the stabilizer  of $\lam$ in $H={\mathcal A}_6$ is isomorphic to ${\mathcal A}_4$ (as $Y$, the stabilizer  of $W$ in $H$ is isomorphic to $S_4$ and the latter group acts non-trivially on $V/W$). Therefore, $|H:H_\lam|=|{\mathcal A}_6:{\mathcal A}_4|=30$. So the \rep $\lam^G$ constructed in Lemma \ref{gg2} is of dimension 30. By Lemma \ref{gg2}, this is \ir and unisingular. 

Furthermore, let $N=N_H(K)$. Then $\lam^G=(\lam ^K)^G$ and $\lam^K|_A$
is the sum of two non-trivial \reps of $A$ with kernel $K$, so the Brauer character of  $\lam^K$ is integrally valued. \itf the Brauer character of $\lam^G$ is integrally valued.
 Then the group $\lam^G(G)\subset GL_{30}(F)$ is conjugate to a subgroup of $Sp_{30}(2)$ by Lemma \ref{tt1}.  \enp

Remark. One can deduce Lemma \ref{31a} from the character table of the group $C_3^4\rtimes \mathcal{A}_6$. The characters of degree $30$ are listed in Table 7.

\subsection{The cases with $2n\in\{22,46, 82,146\}$}

\bl{ob1} Let $G\subset GL_n(3)=GL(V)$ be an \ir cyclic subgroup. Suppose that 
$(|G|,n)\in \{(11,5), (23,11), (41,8),(73,12)\}$. Then $GW=V$ for some subspace $W$ of V of codimension $1$.\el

\bp This result is due to Eamonn O'Brian (University of Auckland, New Zeeland) obtained using his program on the computer program package Magma. 
(Note that $(|G|,n)\in \{(11,5), (23,11)\}$ have been also settled by A. Hulpke.) \enp

\bl{46v} $Sp_{2n}(2)$ contains absolutely \ir unisingular subgroup for $n\in\{11,23$, $41,73\}$.
\el

\bp Let $G=A\rtimes  H$ be a non-trivial semidirect product, where $A$ is an elementary abelian $3$-group of order $3^m$, $m=5,11,8,12$ respectively, for $n=11,23,41,73$, and  $H\cong C_{2n}$ is a cyclic group of order $2n$. Note that the involution $t$, say, of $H$ acts on $A$ as an inversion, that is, $tat\up=a\up$ for $a\in A$. By Lemmas \ref{ob1} and \ref{gg2}, for a non-trivial \rep $\lam:A\ra \FF_4$ the induced \rep $\lam^G$ is a unisingular \irr of degree $2n$. We first show that  $\lam^G$ can be realized over $\FF_2$. For this observe that the kernel of $\lam$ is invariant under $t$,
and hence for $G_1=\lan A,t\ran$ the 2-dimensional \rep  $\lam ^{G_1}$ of $G_1$ is realized over $\FF_2$. Then, of course, $\lam^G=(\lam ^{G_1})^G$ is realized over $\FF_2$. This also implies that the values of the Brauer characters of $\lam ^{G_1}$ and $\lam ^{G}$  are integers. By Lemma \ref{tt1}, $\lam ^{G}$ is conjugate in $GL_{2n}(2)$ to a subgroup of $Sp_{2n}(2)$, as required. 
By Lemma \ref{ct1}, all \ir constitients of $\lam^G(G_1)$ are non-equivalent, and so are 
the \ir constitients of $\lam^G(A)$. Note that $|g|\leq 19$ for $g\in G$. This implies $\lam^G$ to be absolutely \irt  \enp

\subsection{Cases with $2n\in\{174,198,202,218, 246 \}$}


\bl{aa9} Let $G=\mathcal{A}_{9}$. Then G has a unisingular absolutely \irr into $Sp_{78}(2)$.\el

\bp By \cite{JLPW}, $G$ (and also $S_{9}$) has a self-dual \ir 2-modular \rep $\phi$ of degree 78. The Brauer character values of $\phi$ are integers. One easily checks that 1 is an \ei of every element of odd order in $G$. So $\phi$ is unisingular. As $\phi$ self-dual, this is symplectic. By Lemma \ref{tt1} or \cite{HM},   $\phi(G)$ is equivalent to a \rep $G\ra   Sp_{78}(2)$.  \enp


\bp By \cite{JLPW}, $G$ (and also $S_{9}$) has a self-dual \ir 2-modular \rep $\phi$ of degree 78. The Brauer character values of $\phi$ are integers. One easily checks that 1 is an \ei of every element of odd order in $G$. So $\phi$ is unisingular. As $\phi$ self-dual, this is symplectic. By Lemma \ref{tt1} or \cite{HM},   $\phi(G)$ is equivalent to a \rep $G\ra   Sp_{78}(2)$.  \enp

\bl{s47} Let $G=PSp_{4}(7)$. Then G has a unisingular absolutely \irr into $Sp_{174}(2)$.\el

\bp By \cite{HM}, there is an absolutely \irr $\phi:G\ra Sp_{174}(2)$. Note that
$G$ has a cyclic Sylow 5-subgroup of order $25$. By \cite[Proposition 6.2]{TZ23}, the elements of order 25 of $G$ have \ei 1 in every \irr   of degree $n>24$.  

The odd prime divisors of $G$ are $3,5,7$, and the exponent of Sylow $p$-subgroups are $3,25,7$. So it suffices to inspect the case with  $|g|=3,7$.
If $|g|=3$ then $g$ is contained in a parabolic subgroup of $G$; then the result follows from \cite{DZ1}. The case with $|g|=7$ is ruled out  by \cite[Theorem 1.1]{DZ8}, saying that if $\lam$ be an \irr of $G$ such that 1 is not an \ei of $\lam(g)$ then $\tau$ is a Weil \rep of degree $24$.  \enp

\bl{he2} Let $G=He\cdot 2$, where $He$ denotes the Held sporadic simple group. 
Then G has a unisingular absolutely \irr into $Sp_{202}(2)$.\el

\bp By \cite{Lu2}, $G$ has an absolutely \irr $\si:H\ra GL_{202}(F)$, and the Brauer character $\be$ of $\si$ equals $\chi_{4^+}-\chi_{2^+}-2\cdot 1_G$, where $\chi_{2^+}:=\chi_2+\chi_3$, $\chi_{4^+}:=\chi_4+\chi_5$ in notation of \cite{Atlas}. One easily deduces that $\si$ is unisingular. Moreover, the values of $\be$ are integers. Therefore, $\si$ is self-dual, and hence $\si(G)\subset Sp_{202}(F)$. By Lemma \ref{tt1}, $\si(G)\subset Sp_{202}(2)$. \enp

\bl{3d43} Let $G={}^3D_4(3)$. Then G has a unisingular absolutely \irr into $Sp_{218}(2)$.\el

\bp The existence of the \rep in question follows from \cite[Table 2]{HM}.
By \cite{TZ22}, the minimum \po degree of every element $g\in G$ equals $|g|$, 
in particular, 1 is an \ei of $g$ in this representation.\enp

\bl{d83} Let $G=\Omega^+_8(2)$. Then G has a unisingular absolutely \irr into $Sp_{246}(2)$.\el

\bp The existence of the \rep in question follows from \cite{Lu}, where it is  observed that the algebraic group $\mathbf{G}=Spin_8(F)$ has an \ir 2-modular 2-restricted \rep $\phi$ of degree 246, and the \mult of weight 0 equals 6. By Lemma \ref{u7u}, $\phi$ is unisingular. As $\phi$ is 2-restricted, this remains \ir under restriction to $Spin_8^+(2)$ by Steinberg's theorem. It is also well known that $\phi(G)\subset GL_{246}(2)$, in fact,  $\phi(G)\subset Sp_{246}(2)$ as it is self-dual (Lemma \ref{tt1}).  \enp
 
\med
\section{Absence of irreducible unisingular subgroup}

\subsection{Cases $2n=\{10,58,86,106,178, 226\}$}

Let $H=Sp_{2n}(2)$.

\bl{1012} $H$ 
 has no unisingular  \ir subgroup for
$2n=10,58,86,106,178,\\ 226$.\el

\bp Suppose the contrary,  let $G$ be a unisingular \ir subgroup of $H$ and $N$ a minimal normal subgroup of $G$. Note that   in the above list  $n$ are primes.
By Lemma \ref{2po}, for $n$ a prime either $N$ is  simple and  has an \irr of degree $d\in \{n,2n\}$ over  $\overline{\FF}_2$ or
$N$ is an elementary abelian 3-group. 

 By Lemma \ref{sl2}, \cite{Lu} and Lemma \ref{c44}, no simple group of Lie type in defining characteristic 2 has \ir unisingular \rep of degree $n$ and $2n$ for $n$ listed above. So $N$ is not such a group. Suppose that $N$ is not isomorphic to a group of Lie type in characteristic 2. We can ignore the alternating groups $\mathcal{A}_{2n+1}$ and $\mathcal{A}_{2n+2}$ due Lemma \ref{aa3}. Suppose that $N\cong L_2(q)$ for $q>3$ odd. If $d=2n$ then $2n=q+1$ and $3|(q-1)$ or $2n=q-1$ and $q$ is not a prime. One checks that none of these conditions holds for $2n$ in the above list. If $d=n$ then, by   Lemma \ref{h47},   $ 2n=q-1$ and $q$ is not a prime. This is obviously false.

Let $N$ be a simple group isomorphic to neither alternating group nor $L_2(q)$.  The  degrees $d$ of  absolutely  \ir \reps of such simple groups $N$ is tabulated in \cite{HM} up to  $d<250$, together with the minimal fields of realization.  

Let $n=5$. Then $N\in\{ M_{12},M_{22}\}$ \cite{HM}. 
None of these groups is unisingular in $GL_{10}(\overline{\FF}_2)$ as  elements of order 11 do not have \ei 1.  In addition,  none of these groups has an \irr of degree 5 over $ \overline{\FF}_2$. 

Let $2n\in \{58,86,106,178,226\}$. Then $N$ has no absolutely \irr of degree  $n$ and $2n$ over $ \overline{\FF}_2$. 

Suppose that $N$ is  an elementary abelian  3-subgroup. By Lemma \ref{2po}, $N$ contains an \ir subgroup isomorphic to $G_{3,p}$.  Let $m={\rm ord}_3(  n)$.   By Table 1, 
 $m=n-1$ for these $n$. Then, by Lemma \ref{mi2}, $G$ is not unisingular. \enp

\subsection{Cases with $2n\in\{18, 54,162\}$}

\def\glt{ group of Lie type }

\def\airr{ absolutely irreducible representation }

Suppose the contrary, and let $G\subset Sp_{2n }(2)$. We first show that $G$
has no non-trivial abelian normal subgroup: 

\bl{as3} Let $G\subset Sp_{2n}(2)$ with $n|81$ be an \ir unisingular subgroup. Then G has no non-trivial abelian normal subgroup.\el

\bp Suppose the contrary, and let $A$ be a  minimal non-trivial  normal subgroup of $G$. Then $A$ is an elementary abelian $p$-group for some odd prime $p$.  By Lemma \ref{zz1}, $A$ is not cyclic, and $p>3$ by  Lemma \ref{it4}. Let $V$ be the underlying space for $Sp_{2n}(2)$ and   $W_1\ld W_l$ the  quasi-homogeneous components of $V|_A$.  Let $d=\dim W_1$. By Lemma \ref{ct2}, $W_1\ld W_l$ are non-degenerate, so $d$ is even, and    transitively permuted by $G$, hence $ld|162 $ and $l$ is a 3-power. As $p>3$, we have $l< n$ and $l>3$   by Lemma \ref{aa8}.   Therefore, $l\in\{9,27,81\}$. 

 Let $e$ be the  common dimension of all \ir constituents of $W_1$. Then  $p|(2^e-1)$ and   $e|d$. By Lemmas \ref{ny6} and \ref{gg2}, $(3,p-1)=1$,  so $p\neq 7$, hence  $\dim W_1>6$ and $l\neq n/3$, in particular, $l\neq 27$.  
 
So   $l=9$. Then $d=18$ and  $e=9$ or 18  (as $e\neq 2,3$). 
If $e=9$ then $p$ divides $2^9-1=7\cdot 73$, whence $p=73$. 
If $e=18$ then $p|2^e+1$ and $2^9+1=9\cdot (2^6-2^3+1)=27.19$, hence $p=19$.
Both the cases are ruled out by Lemma \ref{ny6} as $3|(p-1)$.\enp

Next we consider the case where $G$ has no non-trivial abelian normal subgroup.
Then, by Lemma \ref{c21}, $G$ has a normal subgroup $N$, say, that is a direct product of non-abelian simple groups.  

  \bl{x54}  Let $G\subset GL_d(\overline{\FF}_2)$, $d|162$, be a unisingular \ir subgroup. Suppose that G is a direct product of simple groups. 
Then   $G= L_2(53) $, $d=54$. In addition, the group $G= L_2(53)$ is not a subgroup of $Sp_{54}(2)$.\el

\bp Suppose the contrary. 
For $d=2,3$ this is known \cite{gl3}.   So we are to inspect the cases with  $d\in\{6,9,18,27,54,81,162\}$.

Suppose first that $G$ is simple.  

 By  Lemmas \ref{sl2} and \ref{20y},  $G$ is not isomorphic to $ L_2(q)$, unless   $q=53,d=54$.  In addition, $G$ is not isomorphic to $\mathcal{A}_{d+1}$ or $\mathcal{A}_{d+2}$ by Lemma \ref{a22}. 

 If $G$ is neither isomorphic to  a group of Lie type in defining characteristic 2 nor to  
$\mathcal{A}_{d+1}$ or $\mathcal{A}_{d+2}$,
then, by \cite{HM1}, we conclude that $d=6$ and  $G\in\{  U_3(3), J_2\}$. By \cite{JLPW}, 1 is not an \ei of $\phi(g)$ for $g\in G$  of order 7, so these groups are ruled out.

 Suppose that $G$ is of Lie type in characteristic 2. 
Then $G\subset \mathbf{G}\subset GL_d(\overline{\FF}_2)$, where $\mathbf{G} $ is a simple algebraic group   so $\phi$ extends to a \rep $\Phi$, say, of $\mathbf{G} $. If $\Phi$ is tensor-indecomposable then, by Lemmas \ref{cc5} and \ref{c44},  either $d\in\{6,27\}$ or
  $G$  is  isomorphic to a finite classical group and $\phi$ is a Frobenius twist of the natural \rep of $G$ or the dual of it;  these \reps are not unisingular by Lemma \ref{c44}).

 Let $d=6$. Then $G\in \{ L_4(q)\cong \Omega^+_6(q), PSU_4(q)\cong \Omega^-_6(q), G_2(q)\}$.  Elements of order 7 of  $G=G_2(2)\subset G_2(2)$ 
are   fixed point free in an \irr  of $G$ of degree 6. The other groups are ruled out by Lemma \ref{c44}. 

Let $d=27$. Then,  by \cite{Lu}, 
 $\mathbf{G}\cong E_6(\overline{\FF}_2)$ and $G\in\{E_6(q),{}^2E_6(q)\}$. By Lemma \ref{e619}, these  examples do  not yield  unisingular representations.  

So $\Phi$, and hence $\phi$,  is tensor-decomposable.  Then $\Phi=\otimes_{j=1}^m \Phi_j$, where $\Phi_j $, $j=1\ld m$, are tensor-indecomposable \ir \reps of $\mathbf{G}$, that are \ir on the restriction to $G$. Let   $d_j=\dim \Phi_j$; we can assume that $d_1\leq d_2<\cdots \leq d_m$. Note that   $d_1\neq 2$ as otherwise 
$G\cong  L_2(q)$ for $q$ even. This has been ruled out above. 

If $d_1=3$ then  $\mathbf{G}\cong SL_3(\overline{\FF}_2)$, and hence $G\cong SL_3(q)$ or $PSU_3(q)$ for $q$ even. This is ruled out by Lemma \ref{tp2}. So $d_1\geq 6$.

 If $d_1=6$ then  $\mathbf{G}\in \{SL_4(\overline{\FF}_2), Sp_6(\overline{\FF}_2), SL_6(\overline{\FF}_2),G_2(\overline{\FF}_2)\}$  and $d_2|27$.
 By \cite{Lu}, these groups have no tensor-indecomposable \irr of degree $9,27$.

 So $d_1\geq 9$ and then $d_1=9$ and $\mathbf{G}\cong SL_9(\overline{\FF}_2)$. Then $d_2\in\{9,18\}$, and in fact $d_2=9$
 as  $\mathbf{G}$ has no \irr of degree 18 \cite{Lu}. This implies $\Phi_i$, $i\in\{1,2\}$ to be a Frobenius twist of an \irr with highest weight $\om_1$ or $\om_8$. By Lemma \ref{tp9}, $\Phi|_G$ is an \ir constituent of the Weil \rep of $G$, which is
     not unisingular by Proposition \ref{pp1}.

 Now suppose that $G$ is not simple, and let $G=S_1\times S_2$, where $S_1\neq 1$ is a simple subgroup of $G$. Then $V=V_1\otimes V_2$, a tensor product of \ir $G$-modules $V_1,V_2$
 where $V_1$ is trivial on $V_2$ and $S_2$ is trivial on $V_1$. By the above, $\dim V_1$ is a multiple of $ 54$, and hence  $\dim V_2=3$. This yields a contradiction.

 The additional claim of the lemma follows from Corollary \ref{22a}.\enp

\bl{no7}  Let $G\subset GL(V)\cong GL_n(\overline{\FF}_2)$ be a unisingular \ir subgroup, where $n|162$, $n>1$, and let  N be a minimal normal \sg of G.  Then $V|_N$
is not homogeneous, unless $N=L_2(53)\subseteq G\subseteq PGL_2(53)$ and $n=54$.
 \el
 
\bp  Suppose the contrary. Then $N$ is isomorphic to a  unisingular \ir subgroup of $GL_d(F)$ with $d|n$. By Lemma \ref{x54},  
   we conclude that $N=L_2(53)$ and $d=54$. Set $M=N\cdot C_G(N)$. Then $G/M$ is isomorphic to a subgroup of ${\rm Out}\,(N)$. As $N$ is simple, $M=N\times C_G(N)$. If $M$ is \ir on $V$ then, by Schur's lemma,  $C_G(N)$ consists of scalar matrices, hence 
$N$ is \ir as claimed.    

Suppose that $M$ is reducible; then $n=162$ and $V|_M$ is a sum of three \ir modules of dimension 54. Therefore, $C_G(N)$ is abelian,
 and  non-scalar (otherwise $M=N$, as a unisingular subgroup contains no non-identity scalar matrix, and then $G\subseteq {\rm Out}\,(L_2(53))=PGL_2(53)$, 
but this group has no \irr of degree 162). Then $G$ is imprimitive, and hence permutes three
\ir $FM$-submodules transitively. Therefore, there is $g\in G$  
which permutes these three submodules. As $|{\rm Aut}\,(L_2(53))|=2$, there exists $h\in N$ such that $gxg\up=hxh\up$
for every $x\in N$, and hence $h\up g\in C_G(N)$. Then $g\in M$, a contradiction. \enp

\bl{332} Let F be an \acf of characteristic $2$ and let $G\subset GL_3(F)$ be a non-abelian simple subgroup. Then $G\in\{SL_3(q), SU_3(q)\}$ with q even and $(3,q-1)=1$ in the former case and $q>2,(3,q+1)=1$ in the latter case.
\el

\bp One easily observes that $G$ is neither an alternating group nor $L_2(q)$ with $q$ odd,
except $q=7$ where $L_2(7)\cong SL_3(2)$. By \cite{HM}, $G$ is a group of Lie type in 
 defining characteristic $2$. Then the lemma follows from Lemma \ref{cc5}.\enp

\bl{na1}  Let $G\subset Sp_{2n}(\overline{\FF}_2)$ be a  unisingular \ir subgroup, where $n|81$, $n>1$, and let  N be a minimal normal \sg of G. Then  either $N\cong L_2(53)$ and $n=27$, or $n=81$ and $N$ is the direct product of three copies of $L_2(53)$.\el

\bp Observe that the exceptions are genuine. If $n=27$ then this follows from Lemma \ref{20y}; if $n=81$ then one can use Lemma \ref{in5a}. 

Let $V$ be the underlying space of $Sp_{2n}(F)$. By Lemma \ref{as3},  $N$ is non-abelian.
Then $N= S_1\times\cdots\times S_k$, where $S_1\ld S_k$ are non-abelian simple groups isomorphic to each other. By Clifford's theorem, 
 $V=V_1\oplus \cdots\oplus V_l$, where  $V_1\ld  V_l$ are the homogeneous components of  $V|_N$, transitively permuted by $G$. We can assume that $S_1$ is non-trivial on $V_1$.
Clearly, $V_1$ is a homogeneous $FS_1$-module. Let $d$ be the dimension of 
an \ir $FN$-submodule of $V_1|_N$ and $d_1$  the common dimension of  \ir  $FS_1$-submodules of $V_1|_{S_1}$. Then $d_1|d$, $d|\dim V_1$ and $l\cdot \dim V_1=2n$.  

\med
(i) The lemma is true if    $ S_1\cong L_2(53)$. 

\med
By Lemma \ref{20y}(1), $d_1\in\{26,52,54\}$,  and $d_1|2n$ implies   $d_1=54$
and $2n\in\{54,162\}$. The former case  is stated in the lemma conclusion, so let $2n=162$. Then $l=3$ by dimension reason,  and   $V_1,V_2,V_3$ are non-isomorphic \ir $FN$-modules by Lemma \ref{no7}. (In fact, by dimension reason, $S_i$ for $i>1$ acts trvially on $V_1$, so $k\leq 3$.)  
Therefore, there is $g\in G$   permuting $V_1,V_2,V_3$. If $k=1$ then $|{\rm Out}\,(N)|=2$, if $k=2$ then $|{\rm Out}\,(N)|=8$. So $g\in N\cdot C_G(N)$ (as $G/(C_G(N)N)\subset {\rm Out}\,(N)$). \itf the $FN$-modules $V_1,gV_1,g^2V_1$ are isomorphic, a contradiction. Then   $k=3$ as stated.
  
So we assume in what follows that $S_1\neq L_2(53)$. By Lemma \ref{x54}, $N$
is not unisingular on every $V_i$, $i=1\ld l$.  

\med
(ii) The lemma is true if $S_1\cong L_2(q)$ with $q$ even or if $d_1=2$.

\med 
As \ir \reps of $S_1$ are of 2-power degrees (Lemma \ref{sl2}) and $d_1|162$, we have $d_1=2$. Moreover, each $S_i$ with $i>1$ is in the kernel of $V_1$ by the same reason. Therefore, for every $V_i$ ($i\in\{1\ld l\}$) there is exacly one
$S_j$ acting on it non-trivially. By reordering them, we can assume that $i=j$ and hence $k=l$. Let $s_i\in S_i$ be of odd order and $s=s_1\cdots s_k$. Then $s_i$ acts  fixed point freely on $V_i$ and trivially on $V_j$ for $j\neq i$. Therefore, $s$ acts  fixed point freely on $V$, a contradiction.  This additionally implies $d_1>2$ due to Lemma \ref{cc5}.

\med
 Let $W_j$, $j=1\ld t$, be the quasi-homogeneous components of $V|_N$ and let $K_j$ be the kernel of $W_j$. We can assume that $V_1\subseteq W_1$. 
 Then $gK_1g\up$ is the kernel of $gW_1$ for $g\in G$.  
As $G$ permutes $W_1\ld W_t$ transitively, we observe
  that $K_1\ld K_t$ are conjugate in $G$, so  $N/K_j\cong N/K_1$ for $j=2\ld t$.

By Lemma \ref{ct2}, $\dim W_1$ is even, so $t$ is odd. 
  So $\dim W_1\in\{6,18,54,162\}$.  

\med
(iii) The lemma is true if $N_1=N/K_1$ is simple.

\med
 In this case  $N/K_j$ is simple for $j=1\ld t$. By reordering $S_1\ld S_k$ we can assume that   $N/K_j\cong S_j$ so $t=k$ here (as $\cap K_j=1$). So $S_j$ acts non-trivially on $W_j$ and trivially on $W_{j'}$ for $j'\in \{1\ld t\}$,  $j'\neq j$. 

Observe that $N_1$  is  unisingular on $W_1$. Indeed, if not, then there are elements $g_j\in S_j$ acting  fixed point freely on $W_j$ and trivially on $W_{j'}$ for $j'\neq j$.
Then $g=g_1\cdots g_m$ acts  fixed point freely on $V$, a contradiction. 

 \smallskip
 
Thus,   $S_1\cong N_1=N/K_1$ is not unisingular on $V_1$ and is unisingular on $W_1$. We show that this leads to a contradiction. In fact, we shall show that some element $x\in S_1$
does not have \ei 1 in every \irr of $S_1$ of degree $d_1$, and hence on $W_1$. 

Suppose that $S_1$ is a group of Lie type in defining characteristic 2. Then  $(S_1,d_1)$ occurs in Tables 3,4. By Lemma \ref{c44}, each of those groups has an element  $g$, say, which does not have \ei 1 on an \irr of $G$ of degree $d_1$,  in particular, $g$ does not have \ei 1 on $W_1$. 

\def\st{Suppose that }

Suppose that $S_1\cong L_2(q)$, $q$ odd. Then $(q,d_1)\in\{(13,6),\,(19,18),\,(37,18 ),\, (53,54)$, $(109,54),\, (163,81), \, (163,162)\}$ by Lemma \ref{20y}. In these cases $q$ is a prime and 
$d_1<q$, except for the pair $(53,54)$ included in the statement. In the remaining case
an element $g\in G$ of order $q$ does not have \ei 1 in every \irr of degree $d_1$.
So  $g$ does not have \ei 1 on $W_1$. 

 \st $S_1$ is not isomorphic to $L_2(q)$. If $S_1\cong \mathcal{A}_m$ then 
$(m,d_1)\in\{(7,6), \, (8,6)$, $(19,18),\,(20,18),\,(55,54 ),\, (56,54),\, (163,162),\, (164, 162)\}$, and elements of order $d_1+1$ do not have \ei 1 on an \irr of degree $d_1$ by Lemma \ref{aa3}. So we conclude as above.

By \cite{HM1}, the only other simple group $S_1$ with an \irr of degree $d_1\in\{6,9,18,27,54,81,162\}$ are $SU_3(3)$ and $J_2$ for $d_1=6$. In these groups
elements of order 7  do not have \ei 1 on an \irr of degree 6 \cite{JLPW}.

\med
(iv) This reasoning in fact shows that if $S\neq L_2(53)$ is a simple group which has an \irr of degree $d_1|162$ then $S$ has an element $x$ acting fixed point freely on every \ir $FS$-module of dimension $d_1$.

\med
 Suppose that $N_1/K_1$ is not simple.

Then $k>1$. 
By reordering of $S_1\ld S_k$ we can assume that   $S_1\ld S_{s}$ acts non-trivially on $V_1$, and $S_{s+1}\ld S_k$ acts trivially (if $k>s$). Then $S_1\times \cdots\times S_s$ acts faithfully on every \ir constituent of $V_1$. Note that $s>1$
and hence $d\geq 9$, which implies     $\dim W_1\geq  18$, $t\leq 9$. Moreover,
it follows from (iv) that $t>1$. Indeed, if $t=1$ then $N$, and hence $S_1$ acts faithfully
on every $V_1\ld V_l$, hence $S_1$ is not unisingular on $V_1$ by (iv). This also implies $s<k$. 

 Let $Y$ be an \ir submodule of $V_1|_N$, so $d=\dim Y\leq 54$. 
 Let $Y=Y_1\otimes \cdots\otimes Y_{s}$, where $Y_i$ is a non-trivial \ir  $FS_i$-module for $i=1\ld s$. 
 Set $d_i=\dim Y_i$, so $d=d_1\cdots d_{s}$. By the above, $d_i>2$ for every $i$,
so $s\leq 3$ (otherwise  $d\geq 81$, whence $\dim W_1=162$ and then $t=1$). We can assume $d_1\leq \cdots \leq d_{s}$.

Thus we assume from now on that $1<t\leq 9$ and $1<s<k$.

\med
Suppose that $d_1=9$. One easily observes that $\mathcal{A}_n$ is not isomorphic to 
an \ir subgroup of $GL_9(\overline{\FF}_2)$. A similar conclusion holds for groups  $L_2(q)$ with $q$ odd by Lemma \ref{20y} and for simple groups by \cite{HM} that are  not isomorphic to groups of Lie type in characteristic 2.

Therefore, $S_1$ is isomorphic to a simple group of Lie type in defining characteristic 2.
Then  either $S_1\cong PSL^\ep_3(q)$ and 
the \ir constituents of $Y|_{S_1}$ are of dimension 9 or $S_1\cong PSL^\ep_9(q)$, see Lemma \ref{cc5}.   Then, by dimension reason, $s=2$ and $\dim Y=81$.
This implies $Y=V_1$, and $V=W_1=V_1+V_2$, that is $t=1$.  Then $V_1|_{S_1}$ is homogeneous and $V_2|_{S_1}$ is the dual of $V_1|_{S_1}$. So we have a contradiction by Lemma \ref{tp2} in the former case and straightforwardly in the latter case. 
 
\med
  Suppose that $d_1=3$. 

\med

By Lemma \ref{332}, $S_1\in \{L_3(q)$, $q$ even, $PSU_3(q)$, $q\geq 2$ even$\}$, and 
  $(3,q- 1)=1$, $(3,q+ 1)=1$, respectively. For uniformity we define $SL^\ep_3(q)$ with $\ep\in\{+,-\}$ to be $SL_3(q)$ if $\ep=+$ and $SU_3(q)$ otherwise.

Groups $SL_3^\ep(q))$, $q$ even, have no 2-modular   \irr of degree 6 or 18 \cite{Lu}, so $d_i\neq 6,18$ for $i=1\ld s$.  
\itf    $d\in\{9,27\}$ is odd.  Therefore $\dim V_1$ is  odd. Indeed, if $\dim V_1$ is even then $V_1$ is non-degenerate by Lemma \ref{ct2}.  As $d=\dim Y$ is odd, $Y$ is totally isotropic, and then $V_1/(Y^\perp\cap V_1)$ is dual to $Y$. As $V_1$ is homogeneous, $Y$ is self-dual, so the Brauer character of $Y$ is real-valued, which contradicts \cite[Ch. IV, Corollary 11.2]{Fe}.
 
\med
(a) Suppose that  $t=3$. Then $d\in\{9, 27\}$ and $V=W_1\oplus W_2\oplus W_3$.
 
By reordering of $S_1\ld S_k$ we can assume that  $N/K_1$
is isomorphic either to  $ S_1\times S_2\times S_3$, or to $ S_1\times S_2$, so  $S_1$ and $S_2$ are non-trivial on $W_1$.  We differ two cases: $(a_1)$ $S_1$ is trivial on $W_2,W_3$ and $(a_2)$ $S_1$ is non-trivial on $W_2$ and trivial on $W_3$ (up to reordering $W_2,W_3$).

\med
$(a_1)$ Since $gS_1g\up$ for $g\in G$ acts trivially on $gW_2$ and $gW_3$, each $S_i$, $i=1\ld k$, acts non-trivially on exactly one of $W_1,W_2,W_3$, in particular, $S_2$ acts trivially on
$W_2,W_3$. One deduces that there are $S_i\neq S_j$ $(i,j\in\{3\ld k\})$ such that $S_i$ is non-trivial on $W_2$ and trivial on $W_1,W_3$ and $S_j$ is non-trivial
on $W_3$ and trivial on $W_1,W_2$. Then, by Lemma \ref{tp2}, we can take $g_1\in S_1$, $g_i\in S_i$, $g_j\in S_j$ such that $g_1g_ig_j$ is  fixed point free on $V$.

\med

  $(a_2)$  We show that $s=2$ and $g_1g_2$ does not have \ei 1 on $V$ for some $g_i\in S_i$, $i=1,2.$

\med
As above, for every $i\in\{1\ld k\}$ there is exactly one $j\in\{1,2,3\}$ such that $S_i$ is trivial on $W_j$.  If $s=3<k$ then one easily checks that this is impossible. So $s=2$ and we assume that $N/K_1\cong S_1\times S_2$. Then $k=3$, and $S_3$ acts non-trivially on $W_2,W_3,$ 
so we can assume that  $S_2$ is trivial on $W_2$ and  non-trivial on $W_3$, $S_1$ is trivial on $W_3$ and  non-trivial on $W_2$. Then $g_1g_2$ with $g_1\in S_1$, $g_2\in S_2$
acts on $W_2$ as $g_1$ does, and on $W_3$ acts as $g_2$ does. 

  Observe that $Z(N)=1$ implies $3\not |(q-\ep 1)$. 

 Suppose first that $S_1\neq SL_3(2)$. Then we choose $g_1,g_2$ so that
 $|g_1|=q-\ep$ and $|g_2|=q^2+\ep q+1$. Then  $(|g_1|,|g_2|)=1$ (as $(3,q-\ep)=1$).
By Lemma \ref{tp2}, $g_2$ does not have \ei 1 on every \irr of $S_2$ degree 3 or 9. We can choose $g_1$ so that $g_1$ does not have \ei 1 on every \irr of $S_1$ degree 3.
By the above, $g_1g_2v=g_1v$ if $v\in W_2$ and $g_1g_2v=g_2v$ if $v\in W_3$. So $g$ daes not have \ei 1   on $W_2$ and $W_3$. As $Y=Y_1\otimes  Y_2$ and $(|g_1|,|g_2|)=1$, we conclude that none of \eis of $g$ on   $W_1$ equals 1. So this case is ruled out.

Suppose that $S_1\cong SL_3(2)$.  Then $d=9$ and $\dim Y_1=\dim Y_2=3$. There are $g_i\in S_i$ of order 7, $i=1,2$, such that the \eis of $g_i$ on $Y_i$ are $\{\zeta,\zeta^2,\zeta^4\}$, where $\zeta^7=1\neq \zeta$. Then $g_1g_2$ has \ei 1 neither on $W_2$ nor on $W_3$. In addition,
 1 is not an \ei of $g_1g_2$ on $Y$, and hence on $V_1$. If $\dim W_1=2\dim  V_1$ then $W_1=V_1+V_2$, where $V_2$ is the dual of $V_1$. 

Suppose that $\dim W_1>2\dim V_1$.  Then $\dim W_1\geq 6\dim V_1\geq 6\dim Y\geq 54$. As $t>1$, we have  so $\dim W_1=54$, $\dim V_1=9$, so $Y=V_1$, $r=6$, that is, $W_1=V_1+\cdots +V_6$. As $S_3$ is trivial on $W_1$ and $V_1\ld V_6$ are non-equivalent $FN$-modules, these are non-equivalent $F(S_1S_2)$-modules. Note that $S_1$ (as well as $S_2$) has only two non-equivalent \ir \reps of dimension 3; so $S_1S_2$ has exactly 4 non-equivalent \ir \reps of dimension 9. This is a contradiction.




This completes our analysis of the case with $t=3$.

\med
(b) Let $t>3$. Recall that $t$ is odd. Then $t\geq 9$. As $d\geq 9$, we have   $\dim W_1\geq 18$, whence $t=9$,  $\dim V=162$, whence $\dim W_1= 18$, $\dim V_1=9$, so $V_1=Y$. Therefore, every $V_1\ld V_l$ is \ir and tensor-decomposable. This implies 
the non-trivial composition factors of $V_{S_i}$ are of dimension 3 for $i=1\ld k$.
In addition, $W_1=V_1+V_2$, where $V_2$ is the dual of $V_1$. Then $N/K_1=S_1S_2$, and $S_i\cong PSL_3(q)$, $(3,q-1)=1$ or $PSU_3(q)$, $3\not|(q+1)$ for every $i$.
Therefore, by Lemma \ref{tp2},  $S_i$ is either trivial or not unisingular on $W_j$ for all $i,j$.

 
Observe that

\med
$(*)$ $S_1$ is non-trivial on some  $W_j$ with $j\neq 1$. 

\med
\noindent  Indeed, otherwise (by reordering $S_1\ld S_k$) we can assume that $S_j$ is non-trivial on $W_j$ and trivial on $W_m$ for $m\neq j$. Then we can take   $h_j\in S_j$ such that $h_1\cdots h_t$ is  fixed point free on $V,$ a contradiction.

 Set $G_1=\{g\in G: gS_ig\up =S_i$ for every $i=1\ld k\}$. Then $xW_j=W_j$ for every $x\in G_1$ and $j=1\ld t$ as $xS_1S_2x\up=S_1S_2$ acts faithfully on $xW_1$. If $h:G\ra Sym(W_1\ld W_t)$ then $G_1\subseteq $ker$(h)$. We show that ker$(h)=G_1$.
Indeed otherwise there is $g\in G$ such that $gW_j=W_j$ for $j=1\ld t$ and $gS_ig\up\neq S_i$ for some  $i\in\{1\ld k\}$. We can assume $i=1$ here. Then $gW_1=W_1$ implies $gS_1g\up= S_2$. This in turn implies $S_1$ to act trivially on $W_j$ for $j>1$. (Indeed,  $gW_j=W_j$ implies $S_2=gS_1g\up$ to acts non-trivially on $W_j$, which violates the definition of $W_j$.) This contradicts $(*)$.



$(b_1)$ We show that $k\leq t$. Set $X=G/G_1$. Then $X$ acts faithfully and transitively on $S_1\ld S_k$ and on $W_1\ld W_t$.
By the definition of $W_1\ld W_t$, for every $j\in \{1\ld t\}$ there is a unique pair $S_a,S_b$
$(1\leq a\neq b\leq k)$ such that $S_aS_b$ acts faithfully on $W_j$. If $gW_1=W_j$
then $S_aS_b=g(S_1S_2)g\up$. By the above, this yields a faithful transitive \rep
$X\ra Sym(W_1\ld W_t)$, which is permutationally isomorphic to the action of $X$
on the orbit of $X$ on the unordered pairs $S_aS_b$ containing $S_1S_2$. By Lemma \ref{a33},
either $k=2t$ or $k\leq t$. In the former case all pairs $S_aS_b$ in the orbit are disjoint,
so $S_1$ acts trivially on every $W_j$ with $j\neq 1$. This violates (*).

\med
$(b_2)$ Next we show that $k=t$. Suppose  that $k<t=9$. We have seen in the proof of Lemma \ref{a33} that $9=k\cdot |X_1:X_{1,2}|/|X_{(1,2)}:X_{1,2}|$, where $X_{1}$ is the stabilizer of $S_1$ in $X$, $X_{(1,2)}$ is the stabilizer of  $S_1S_2$ and $X_{1,2}$ is the stabilizer of both $S_1,S_2$. Note that $k>3$ as otherwise $t\leq 3.$
Therefore, $k=6$, $|X_{(1,2)}:X_{1,2}|=2$ and $|X_1:X_{1,2}|=3$. Then after reordering 
$W_1\ld W_9$ and $S_1\ld S_6$ we can assume that the $j$-th pair in the list  $S_1S_2$,  $S_1S_4$,  $S_1S_6$,  $S_3S_2$,  $S_3S_4$,  $S_3S_6$,  $S_5S_2$,  $S_5S_4$,  $S_5S_6$ acts faithfully on  
$W_j$ for $j=1\ld 9$. Let $g_i\in S_i$ for $i=1,3,5$ be an element of order $q^2+\ep q+1$,
and $g=g_1g_3g_5\in S_1S_3S_5\subset N$. Then $g$ acts  fixed point freely on every $W_j$, $1\leq j\leq 9$, and hence on $V$. This is a contradiction.   

\med

 So we have  $k=t=9$.  Then $G/C_G(N)$ has a transitive subgroup $R$  of order 9 (Lemma \ref{pp1}). This naturally acts on the pairs $S_aS_b$ ($a,b\in 1\ld k, a\neq b$); specifically  $x\in R$ sends   $S_aS_b$ to $ S_{x(a)}S_{x(b)}$. As $R$ permutes $S_1\ld S_k$ transitively, $\{x(a):x\in R\}=\{x(b):x\in R\}=\{S_1\ld S_k\}$. Therefore, $x(b)=a$ for some $x\in R$ with $x\neq 1$, so the above list contains  a pair $S_{x(a)}S_{a}$ for some $x\in R$ with $x(a)\neq b$. 

Observe that $R$ is transitive on $W_1\ld W_9$. If not, then $R$ has an orbit of size 3, hence some some $1\neq x\in R$ fixes the $W_i$'s in this orbit. Then $x(W_i)=W_i$ implies
$x(S_aS_b)=S_aS_b$, where $S_a\neq S_b$ act non-trivially on $W_i$. As $|x|$ is a 3-power, we have $x(S_a)=S_a$, a contradiction.  \itf no element  $1\neq x\in R$ fixes $W_i$. 

Let $1\neq x\in R$. Then we can assume that $x$ permutes $W_1,W_2,W_3$,
$W_4,W_5,W_6$ and $W_7,W_8,W_9$. By reordering $S_1\ld S_9$ we can assume that $S_1S_a,S_2S_b,S_3S_c$ acts faithfully on $W_1$, $W_2$, $W_3$, respectively, and that $x(S_1)=S_2, $ $x(S_2)=S_3$, $x(S_3)=S_1 $. Then $x(S_a)=S_b$, 
$x(S_b)=S_c, $ $x(S_c)=S_a $. 

\med
$(b_3)$ We show that  $\{a,b,c\}\cap \{1,2,3\}$ is not empty. Suppose the contrary. Then we can assume that $a=2$. As $x(a)=b,x(b)=c$, we have $\{a,b,c\}= \{1,2,3\}$. Then   $S_1S_2$ is faithful on $W_1$, $S_{2}S_{3}$ is  faithful on $W_2$, $S_{3}S_{1}$ is  faithful  on $W_3$. Set $g=g_2g_3$, where $g_2\in S_2$,  $g_3\in S_3$ with $|g_2|=| g_3|=q^2-\ep q+1$.  

Recall that $W_1=V_1+V_2$, where $V_1,V_2$ are dual.
 Therefore, $W_1|_{S_1} $ is a sum of \ir \reps of dimension 3. Therefore, either $W_j|_{S_i}$ is a sum of \ir \reps of dimension 3, or $S_i$ acts on   $W_j$ trivially. Let $Y_j$ be an \ir constituent of $W_j|_N$. 
 Then $g$ acts on $Y_1$ as $\Id\otimes g_2$, as  $g_2\otimes g_3$  on $Y_2$ and as $g_1\otimes \Id$ on $Y_3$.  \itf $g$ acts fixed point freely on $W_1+W_2+W_3$ \ii $g$ acts  fixed point freely on $Y_1+Y_2+Y_3$. By Lemma \ref{c44}, $g$ acts 
fixed point freely on $Y_1$ and $Y_3$. Let $\lam\otimes \mu$ be the \rep of $S_2\times S_3$ on $Y_2$, where $\lam , \mu$ are some   3-dimensional \ir \rep of $S_2,S_3$, respectively.
Then it is clear that we can choose  $g_2,g_3$ such that $\lam(g_2)\otimes \mu(g_3)$ would not have \ei 1. Then $g$ is  fixed point free on $W_1+W_2+W_3$.

Similarly, we can assume that $S_4,S_5,S_6$ are non-trivial on $W_4+W_5+W_6$ and 
trivial $W_j$ with $j\neq 4,5,6$, and $S_7,S_8,S_9$ are non-trivial on $W_7+W_8+W_9$ and trivial $W_j$ with $j<7$. The above argument shows that there are elements $g'\in S_8S_9$ acting fixed point freely on $W_4+W_5+W_6$, and $g''\in$

Let $ y\in R$, $y\neq 1,x,x^2$. One easily observes that there are exactly three groups $y(S_1),y(S_2),y(S_3)$ that act non-trivially on $y(W_1)+y(W_2)+y(W_3)$, and similar argument work, leading to the conclusion that  some element of $y(S_1)y(S_2)y(S_3)$ acts  fixed point free on $y(W_1)+y(W_2)+y(W_3)$. In turn, this implies that $N$ is not unisingular on $V=W_1+\cdots +W_9$.  

\med
$(b_4)$  Suppose first that  
$\{a,b,c\}\cap \{1,2,3\}$ is empty.

\med
 As $x$ permutes $S_a,S_b,S_c$, it follows that 
$S_a,S_b,S_c$ acts non-trivially on exactly one space $W_4+W_5+W_6$ or $W_7+W_8+W_9$; we can assume that $S_a,S_b,S_c$ acts non-trivially on $W_4+W_5+W_6$. Similarly, 
$S_1,S_2,S_3$ acts non-trivially on exactly one of these spaces, which is  $W_7+W_8+W_9$ as otherwise we arrive at the case (a).
Then we take  $g_i\in S_i$ with $i=1,2,3$ to be of order  
$q^2+\ep q+1$, and also $g_a\in S_a$, $g_b\in S_b$, $g_c\in S_c$ of order $  q-\ep 1$. If $S_1$ is not isomorphic to $ SL_3(2)$ then we can choose $g_a,g_b,g_c$
to be  fixed point free at the natural $SL_3^\ep(q)$-module, and then these are fixed point free  on every 3-dimensional \irr of  $SL_3^\ep(q)$. Set $g =g_1g_2g_3g_ag_bg_c\in N$. 
 Then $g$ obviously acts  fixed point freely on each space $W_4,W_5,W_6,W_7,W_8,W_9$. 
The action of $g$ on $W_1$ is realized via $\lam(g_1)\otimes \mu(g_a)$ for some 3-dimensional \reps $\lam,\mu$ of $S_1$, $S_a$, respectively,   and $g$ acts on $W_2,W_3$ similarly. As $(3,q-\ep 1)=1)$,   it follows that $(q^2+\ep q+1, q-\ep 1)=1$.  Therefore, $(|g_1|,|g_a|)=1$, and hence $\lam(g_1)\otimes \mu(g_a)$ does not have \ei 1.  We conclude that  $S_1S_a$ is not unisingular on $W_1$.
A similar conclusion follows for $W_2$ and $W_3$. Therefore,   $g$ acts  fixed point freely on $W_1+W_2+W_3$, and also on $V$ by the above. 

Suppose that $S_1\cong SL_3(2)$. This group has exactly 2 \ir \reps of degree 3, the natural one and its dual.  Then the \eis of $\lam_1(g_1)$ are $\nu,\nu^2,\nu^4$
for some primitive 7-root of unity $\nu$. If $|g_a|=7$ then the \eis of $\mu_1(g_a)$ or $\mu_1(g_a\up)$ are $\nu,\nu^2,\nu^4$.  So  $\lam_1(g_1)\otimes \mu_1(g_a)$ does not have \ei 1 for some $g_a\in S_a$. Similarly, there are $g_2\in S_2,g_b\in S_b$ and $g_3\in S_3,g_c\in S_c$, each of order 7,  such that of $\lam_2(g_2)\otimes \mu_2(g_b)$ and $\lam_1(g_3)\otimes \mu_1(g_c)$ do not have \ei 1. As above, set $g=g_1g_2g_3g_ag_bg_c\in N$.
Then $g$     does not have \ei 1 on $W_1+W_2+W_3$ by the choice of 
$g_1,g_2,g_3,g_a,g_b,g_c$. In addition, as $g\in S_1S_2S_3S_aS_bS_c$, the action of $g$ on each $W_i$ with $i>3$ reduces to the action of exactly one of 
$g_1,g_2,g_3,g_a,g_b,g_c$, hence this is  fixed point free. This completes our consideration of the case with $t=9$, and hence $d_1=3$. \enp

\bl{a16}  Group $Sp_{2n}(2)$, $n|81$, has no unisingular  \ir subgroup.\el

\bp Suppose the contrary, let $G$ be such a subgroup. By Lemma \ref{c21},
a minimal normal subgroup $N$ of $G$ is either elementary abelian or a direct product of non-abelian simple groups. The former is not the case by Lemma \ref{as3}. 
In the latter case let $S$ be a minimal subnormal subgroup of $G$. Suppose first that
$G$ is absolutely \irt Then, by Lemma \ref{na1}, we have only to inspect the case with 
$S\cong L_2(53) $ and $2n\in\{54,162\}$. By Lemma \ref{x54}, $2n\neq 54$. We show that $2n\neq 162$.
 By Lemma \ref{na1}, $N$ is a direct product of three copies of $S$. Let $V$ be the underlying space for $Sp_{162}(2)$. Obviously, $N$ is reducible, and $V=V_1\oplus V_2\oplus V_3$, where $V_1,V_2,V_3$ are non-degenerate \ir $FN$-modules. As $S$
is not a subgroup of $Sp_{54}(2)$ (Lemma \ref{na1}), it follows that $S$ is irreducible,
which contrdicts  Schur's lemma as $N\neq S$.

Suppose that $G$ is not absolutely \irt Let $m$ be the composition length of $G$ viewed as a
subgroup of $GL_{2n}(\overline{\FF}_2)$. Then $G$ is isomorphic to a unisingular subgroup of $GL_{2n/m}(2^m)$, see Lemma \ref{na3}. As $m>1$, we have $2n/m\geq 54$ by Lemma \ref{20y}, whence $m=3$. However, by \cite[Table 2]{HM1}, we observe that the minimum realization field  of  $L_2(53)$ is 
 the index 2 subfield of  $\FF_2(\zeta)$, where $\zeta$ is a primitive 13-root of unity (as  $|q-1|_{2'}=13$ for $q=53$). One easily checks that $|\FF_2(\zeta):\FF_{2}|=12$, so $m=6$, which gives a contradiction.    \enp

\def\hmm{homogeneous }
 

  \vskip2cm
\section{Unisingular irreducible linear groups in odd characteristic}

\bp[Proof of Theorem {\rm \ref{fo8}}] Suppose that $n>1$ is odd. Let $D$ be the group of all diagonal matrices $x$ of determinant $1$ such that $x^2=1$, and let $C$ be the cyclic group of order $n$ consisting of permutational matrices. Then $G:=DC$ is absolutely \ir (which is well known) and unisingular.  Indeed, let $g\in G$ and $|g|=m$. If $m=2$ or  odd then the result is obvious. Otherwise, $h:=g^{m/2}$ is an involution, and hence lies in $D$ and in $  C_G(C)$. \itf that $h\in Z(G)$. As $G$ is absolutely irreducible, $Z(G)$ consists of scalar matrices, so $h=-\Id$, which is false as $\det h=1$. 

Suppose that $n$ is not a 2-power. Then $n=2^lm$, where $m>1$ is odd.  Then $GL_n(F)$ has a unisingular absolutely \ir  subgroup by Lemma \ref{in5a}.  

Suppose that $n$ is a 2-power.  If $n=8$ and $r\neq 3$ then the group $AGL_1(9)\cong C_3^2\cdot C_8\subset GL_8(F)$ is unisingular. This follows  
by inspection of the character table of $G:=C_3^2\cdot C_8$
available in \cite{s31}, and can be also deduced from Theorem \ref{af1}. Indeed, let $\nu$ be the character of $G$ of degree 8. 
As $C_8$ acts transitively on the set of non-trivial elements of $A=C_3^2$, it follows that every non-trivial character of $A$ occurs in $\nu|_A$. Therefore, 
$\nu|_A=\rho_A^{reg}-1_A$, and hence $\nu|_A$ is unisingular and integer-valued. As $C_8$ permutes the one-dimensional constituents of $\nu|_A$,
it follows that the same true for $\nu|_{C_8}$ and consequently for $\nu$ itself.  
Furthermore, the reduction of $\nu$ modulo every prime $r\neq 3$ is irreducible,
and takes values in $\FF_r$. Therefore, $\nu(G)$ can be realized over $\FF_r$,
that is, $\nu(G)$ is conjugate to $GL_8(r)$.  
 
Consider the cases $n=2,4$ and $n=8,r=3$. Let $V$ be the underlying space for $GL_n(F)$.  
Suppose the contrary, and let $G\subset GL_8(F)$ be a unisingular \ir subgroup. By  Lemma \ref{c21}, either $G$ has a non-trivial abelian normal subgroup or a simple subnormal subgroup.

(i) Suppose first that  $G$ has a non-trivial abelian normal subgroup $A$, say. 
We can choose for $A$ an elementary abelian $p$-group for some prime $p$.
Observe that $p>2$. Indeed, suppose the contrary. Then $V=\sum_{i=1}^e V_i$, where 
$V_i$'s are homogeneous components of $V|_A$.  By Clifford's theorem, they are of the same dimension $d$, say, and transitively permuted by $G$. Let $N=\{g\in G:gV_i=V_i\}$ for every $i$. Then $G/N$ is isomorphic to a transitive subgroup of $S_{e}$, so $e$ is a 2-power.  By Lemma \ref{pp1}, $G/N$ contains a transitive 2-subgroup $B$, say, and hence $G$ contains a 2-group $B_1$ which transitively permutes $V_1\ld V_e$. We can assume that $B_1$ is a Sylow 2-subgroup of $G$, and hence $A\cap Z(B_1)\neq 1. $ This contradicts Lemma \ref{am1}.


So $|A|$ is odd, and $p\geq 5$ as $p\neq r=3$.  We choose for $A$ a minimal non-trivial $G$-invariant abelian subgroup. Let $W_1\ld   W_l$ be the quasi-homogeneous components
of $V|_A$. As $H$  is irreducible, we have $l|8$ by Clifford's theorem.
be the decomposition of $V|_A$ as in 
By Lemma \ref{aa8}, $l>p\geq 5$, whence $l=8$ and $p\in\{5,7\}$. Therefore,
$\dim W_1=\cdots =\dim W_8=1$ so $W_1\ld W_8$ are non-isomorphic $FA$-modules.

Let $S_2$ be a Sylow 2-subgroup of $G$. By Lemma \ref{pp1}, $S_2$ acts transitively on $W_1\ld W_8$; as these are non-isomorphic $FA$-modules, the group  $G_1:=AS_2$
is \irt Let $A_1$ be a minimal non-trivial $G$-invariant abelian subgroup of $G_1$. 
As above, by Lemma  \ref{aa8},  the number of quasi-homogeneous components
of $V|_{A_1}$ equals  $8$, so thise are $W_1\ld W_8$.  

We show that $G_1$ is not unisingular. 
For this, let $1\neq z\in Z(S_2)$ be of order 2. Then $z_1$, the projection of $z$ in 
$GL(A_1)$ is diagonalizable and $z_1^2=1$.  So either $zaz\up =a\up$ for every 
$a\in A_1$ or $C_{A_1}(z)\neq 1$. In the latter case  $C_{A_1}(z)\neq 1$ is $S_2$-invariant, which contradicts the minimality of $A_1$.  In the former case we have $zW_i\neq W_i$ for every $i$. (Indeed,  $zW_i= W_i$ for some $i$ implies  $zgW_i= gW_i$ for every $g\in S_2$, and moreover the $F\lan z\ran$-modules $gW_i$ are isomorphic as $zg=gz$. As $\dim W_i=1$, we have $z=-\Id$, so $z$ is  fixed point free on $V$, a contradiction.)  
So $zW_i\neq W_i$ for every $i$. Let $a\in A_1$ and $w\in W_i$. If $aw=w$ then $azw=z\cdot z\up azw=za\up w = zw$, so $W_i$ and $zW_i$ has the same kernel,
which contradicts the definition of the decomposition $V=W'_1\oplus \cdots \oplus W'_{l'}$.  
 
\med
(ii) Suppose first that $G$ has a simple subnormal subgroup $S$, and let $N$ be  the minimal normal subgroup of $G$ containing $S$.  Let $d$ be the dimension of an 
\ir \ccc of $S$ on $V$. Then $d\neq 2$ as every quasi-simple subgroup of $GL_2(F)$ has non-trivial center (for $r>2$).  In addition, by Clifford's theorem,  $d|n$, so $d$ is a 2-power. 
\med

Suppose $n=4$.  Then $d=n$ and $S\in\{ PSL_2(r^a)$, $a>1$, $\mathcal{A}_5$ for $r\neq 5$,  $\mathcal{A}_6$ for  $r=3\}$, see \cite{HM}. (Each group $SL_4(r^a),SU(r^a),Sp_4(r^a)\subset GL_4(F)$  has non-trivial center, so is not simple. This also rules out the $r$-restricted \ir \rep of $SL_2(r^a)$ in dimension 4, which exists for $r>3$.) Groups $\mathcal{A}_5,p\neq 5$ and $\mathcal{A}_6$ are rules out as element of order 5 is not unisingular.  Let $S\cong PSL_2(r^a)$, $a>1$.
Then $S$ arises as the tensor product of two non-equivalent \ir \reps of $SL_2(q)$
of degree 2.    As $S$ is unisingular, we have a contradiction by  Corollary \ref{1k4}. 
 

So $n>4$ and hence   $n=8,r=3$. Let $V$ be the underlying space of $GL_8(F)$, where $F$ is a field of characteristic 3. Observe first that $H$ is absolutely \irt Indeed, otherwise 
$H$ is isomorphic to an absolutely \ir unisingular subgroup of $GL_k(\overline{F})$,
where $k|8$ and $\overline{F}$ is an \acf containing $F$. By the above, this is false. So we can assume that $F$ is algebraically closed. By Clifford's theorem, $V|_S$ is a completely reducible $FS$-module. Let $d$ be the dimension of a non-trivial  \ir constituent of $S$ on $V$. Then $d|8$. In addition $d\ne2$ as $GL_2(F)$ has no \ir  simple subgroup   in characteristic $r>2$. 

One of the \f holds: (i) $S$ is a group of Lie type in defining characteristic 3; (ii)
$S\cong L_2(q)$ for some $q$ with $(3,q)=1$; (iii) $\mathcal{A}_9$ for $r\neq 3$, $\mathcal{A}_{10}$ for $r= 5$, (iv) $S$ is another simple group. Case (iii) does not appear as $r=3$.
 Inspection of \cite{HM} rules out case (iv) as there are no simple \ir subgroups of $GL_d(F)$ with $d=4,8$ that are not of Lie type, $\mathcal{A}_m$  or $L_2(q)$.
If $S\cong L_2(q)$ and $(3,q)=1$ then $q\leq 17$ as $SL_2(q) $ has no non-trivial \irr of degree less than $(q-1)/2$. So $q\in \{5,7,8,9,11,13,16,17\}$.  By inspection in \cite{JLPW}, none of groups $L_2(q)$ for  these values of $q$ has an \ir 3-modular \rep of a 2-power degree (but $SL_2(q)$ may have). Therefore, $S$ is a group of Lie type in defining characteristic 3. 

Suppose first that $d=8$. Then $S$ is \ir and hence either $S$ is one of the classical groups of degree 8, or ${}^3D_4(3^a)$, $a$ odd, 
or $S\cong SL_2(3^a)$. It follows from the general theory of \reps of groups of Lie type that
the groups $SL_8(q)$, $SU_8(q)$ and  $Sp_8(q)$ can be realized as subgroups of $GL_8(F)$
only via their natural \rep (the other \reps differ from the natural one by a Frobenius (or Galois) twist). As $q$ is odd, the center of each of these groups is non-trivial, and not in the kernel of any their \reps of degree 8. Let $S$ be of type $D_4(q)$, ${}^2D_4(q)$, or ${}^3D_4(q)$, where $q$ is a 3-power. The simple group $D_4(3)$ has a subgroup $X\cong \Omega_8^+(2)$, see \cite[p.144]{Atlas}. Let $\phi$ be an \irr of $S$ of degree 8.
Then $\phi|_X$ has an \irr over $F$ of degree 8. This is, however, false by \cite[p. 233]{JLPW}. (Note that $S$ and $X$ have projective \ir \reps of degree 8.) This rules out the case with $S$ of type $D_4(q)$. 

Let $S$ be of type  ${}^2D_4(q). $ Then $S$ contains a cyclic subgroup of order  $q^4+1=3^{4a}+1$. By Zsigmondy's theorem \cite[Theorem 5.2.14]{KL}, there is a prime $t$, say, dividing $q^8-1$ and coprime to $3^i-1$ for every $i<3^{4a}-1$. Let $g\in S$ be of order $t$. Then $g$ is irreducible on a vector space of dimension 8 over $\FF_q$, hence 1 is not an \ei of $g$.
In addition, every \irr of degree 8 can be realized over $\FF_q$. 

  Let $S$ be of type  ${}^3D_4(q). $ Then $S$ is a subgroup of $O_8^+(q^3)$. In addition, every \irr of of $G$ of degree 8 can be realized over $\FF_{q^3}$. It is well known that $S$
contains a subgroup of order $q^4-q^2+1=(q^{6}+1)/(q^2+1)=(q^{12}-1)/(q^6-1)(q^2+1)$. As above, there is a prime $t$ dividing $q^{12}-1$ and coprime to $3^i-1$ for every $i<3^{12a}-1$. Let $g\in S$ be of order $t$. Then $g$ is irreducible on a vector space of dimension 8 over $\FF_{q^3}$, hence 1 is not an \ei of $g$.
 
Note that $3$-restricted \ir \reps of $\mathbf{S}=SL_2(F)$ are of degree $1,2,3$. By the Steinberg tensor product theorem, $S$ is a tensor product of three \ir \reps of degree 2, but this violates the fact that $Z(S)=1$. (Alternatively, one can use   Lemma \ref{12a}.) 



Next assume that $d<8$. Then $N$ is reducible. Indeed, otherwise $N\neq S$, and hence $V$ is a reducible  homogenious $FS$-module. Then $V$ has a unisingular \ir  $FS$-submodule, contrary to the above.   

   Suppose that $d=4$. Then $V|_N=V_1+V_2$, where 
$V_1,V_2$ are \ir $FN$-modules. If $N\neq S$ then $N=S_1\times S_2$
with $S_1\cong S_2\cong S$, and $S_i$ is trivial on $V_i$ for $i=1,2$ up to reordering of $S_1,S_2$. This contradicts the above. So $N=S$. \itf $V_1,V_2$
are non-isomorphic \ir $FS$-modules of dimension 4. Then $S$ is either a classical group of dimension 4 (such as $PSL_4(3^a)$, $PSU_4(3^a)$, $PSp_4(3^a)$)
or $SL_2(3^a)$ in some \irr of dimension 4. If $S$ is classical then $S$ has no 
\irr of dimension 4. (For instance, every \irr of $SL_4(3^a)$ of dimension 4 is faithful on  $Z(SL_4(3^a))$.)  

Let $S=L_2(3^a)$. By Corollary \ref{1k4}, the elements of order
$(3^a+1)/2$ do not have  \ei in every \irr of degree 4, in particular, on $V_1$ and on $V_2$, and hence on $V$.\enp

 \med

\section{The tables}


Table 5 
lists the values of  $n$ such that $Sp_{2n}(2)$ contains
an absolute \ir unisingular subgroup. The table is organized as follow: the first column gives
some values of $2n$ with $1\leq n< 125$. The second column either exposes a group $G$ that is isomorphic to an  absolutely \ir unisingular subgroup of $Sp_{2n}(2)$ or states 
"none" if such group does not exist, or states "open" if the question on the existence of a group in question remains open.  In view  of Corollary \ref{in5}, if a subgroup $G\subset Sp_{2k}(2)$ with the property required exists, then it exists in  $Sp_{2n}(2)$ for $n$ a multiple of $k$. So the third column prints the values of $2n$ that are multiples of $2k$,
with omitting those already appeared in any $2l$-row for $l<k$. 
The forth column refers to the result justifying $G$ in the second  column to indeed exist.

Table 6  borrowed from \cite{s31}.

 
\vskip1cm
{\large Acknowledgement}. A part of this work has been performed during author's visit to the Isaacs Newton Institute for Mathematical Sciences, Cambridge, for his participation in the programme "Groups, representations and applications: new perspectives" (July 2022). I would like to thank the Isaac Newton Institute  for the hospitality during my visit supported by EPSRC grant no EP/R014604/1. A fruitful discussion of Problem  11 of
 Section 1   during the programme with Eamonn O'Brian  leaded to a proof of Lemma \ref{ob1}. I am undebted to Alexander Hulpke who proved a special case of Lemma \ref{ob1}. 

I express my thanks to John Cullinan who introduced me to the circle of problems
related with the unisingularity problem for subgroups of $Sp_{2n}(2)$.

\newpage


\begin{center}  ${}~~~~~~~\,\,\,\,\,\,\,$ Table 5. The values of  $n$ such that $Sp_{2n}(2)$ contains\newline
an absolute \ir unisingular subgroup 
{\small
\begin{tabular}{|c|c|l|l|}
\hline
$2n$&$G$&\,\,\,\,\,\,\,\,\,\,\,\,\,\,\,\,\,\,\,\,\,\,\,\,\,multiples of $2n$ &reference\\
\hline
$2,4,6$& none& &\cite{cu10}
\\ \hline
$8$& $AGL_1(9)$&16,24,32,40,48,56,64,72,80,88,96,104, &Lemma \ref{in5}\\
&&112,120,128,136,144,152,160,168,176,&
\\
&&184,192,200,208,216,224,232,240,248&\\
\hline
$10$ & none& &Lemma \ref{1012}
\\ \hline
$12$&  $C_3^3{\rtimes}A_4$&  36, 60, 84, 108, 132, 156, 180, 204, 228 &Lemmas \ref{c3a4} and \ref{in5}\cr
 \hline
$14$& $L_2(13)$& 28,42,70,98,126,140,154,182,196,210,238&Lemmaa \ref{22a} and \ref{in5}
\\ \hline
$18$& none &$-$&Lemma \ref{a16}
\\ \hline
$20$& $L_2(19)$&100, 220&Lemmas \ref{22a} and \ref{in5}\\
 \hline
$22$& $3^5\rtimes C_{11} $& 44, 66, 110, 198, 242& Lemmas \ref{46v} and \ref{in5} \\
\hline
$26$& $L_2(25) $&  52, 78, 130, 234&  Lemmas \ref{22a} and \ref{in5}
\\ \hline
$30$& $C_3^4\wr S_5$&90, 150 &Lemmas \ref{31a} and  \ref{in5}
\\ \hline
$34$& $L_6(2)$ & 68, 102, 170, 204 &Lemma \ref{ad6} and \ref{in5}
\\ \hline
$38$& $L_2(37)$& 76, 114, 190, 228&Lemmas \ref{22a} and \ref{in5}
\\ \hline
$46$& $C_3^{11}{\rtimes}C_{23}$& $92, 138, 184,230$&Lemmas \ref{46v} and \ref{in5}
\\ \hline
$50$& $L_2(49)$& &Lemma \ref{22a}
\\ \hline
$54 $& none& &Lemma \ref{a16}
\\ \hline
$58$& none& &Lemma \ref{1012} 
\\ \hline
$62$& $L_2(61)$& 124, 186&Lemmas \ref{22a} and \ref{in5}
\\ \hline
$74$& $L_2(73)$&148, 222&Lemmas \ref{22a} and \ref{in5}
\\ \hline
$78$& $A_9$&156& Lemma \ref{aa9}
\\ \hline
$82$& $ C_3^8\rtimes C_{41}$& 164&  Lemma \ref{46v}
\\ \hline
$86$&  none& &Lemma \ref{1012}
\\ \hline
$94$& $open$& & 
\\ \hline
$106$&  none& & Lemma \ref{1012} 
\\ \hline
$116$& $open$&  & 
\\ \hline
$118$& $ SO^+_{16}(2)$& 236&  Corollary \ref{22y} 
\\ \hline
$122$& $ L_2(11^2)$&244& Lemmas \ref{22a} and \ref{in5}
\\ \hline

$132$& $ E_7(2)$& &  Corollary \ref{22y}  
\\ \hline
$134$& $open$& & 
\\ \hline
$142$& $ L_{12}(2)$& & Corollary \ref{22y}\\%
 \hline
$146$& $G= C_3^{12}.73$ & &Lemma  \ref{46v}
\\ \hline
$ 158$&  $L_2(157)$& & Lemma \ref{22a}
\\ \hline
$162$& none& &Lemma \ref{a16}
\\ \hline
$166$& $open $& &
\\ \hline
$172$& $open$  & & 
\\ \hline
$174$& $G=PSp_4(7) $  & & Lemma \ref{1012} 
\\ \hline
$178$&  none&  & Lemma \ref{mi2}
\\ \hline
$ 188$&$SO^+_{20}(2) $  & & Corollary \ref{22y} 
\\ \hline
$194$& $L_{14}(2) $&& Corollary \ref{22y} 
\\ \hline
$202$& $He:2$&&Lemma \ref{he2}
\\ \hline
$206$ &$open$ &  & 
\\ \hline
$212$& $L_2(211)$ &&  Lemma  \ref{22a} 
\\ \hline
$214 $& $open$ & &
\\ \hline
$218 $ &${}^3D_4(3 )$&&Lemma \ref{3d43} 
\\ \hline
$226$&  none& & Lemma \ref{1012} 
\\ \hline
$230$&$SO^+_{11}(2)$&&  Lemma \ref{22y} 
\\ \hline
$246$&   $SO^+_8(2)$&& Lemma \ref{d83} 
\\ \hline
\end{tabular}}
\end{center}

\newpage
 

\centerline{Table 6. Character table of $C_3^3{\rtimes}A_4$}
$$
\begin{array}{c|rrrrrrrrrrrrr}
  \rm class&\rm1&\rm2&\rm3A&\rm3B&\rm3C&\rm3D&\rm3E&\rm3F&\rm6&\rm9A&\rm9B&\rm9C&\rm9D\cr
  \rm size&1&27&4&4&6&12&36&36&54&36&36&36&36\cr
\hline
  \rho_{1}&1&1&1&1&1&1&1&1&1&1&1&1&1\cr
  \rho_{2}&1&1&1&1&1&1&\zeta_3&\zeta_3^2&1&\zeta_3^2&\zeta_3^2&\zeta_3&\zeta_3\cr
  \rho_{3}&1&1&1&1&1&1&\zeta_3^2&\zeta_3&1&\zeta_3&\zeta_3&\zeta_3^2&\zeta_3^2\cr
  \rho_{4}&3&-1&3&3&3&3&0&0&-1&0&0&0&0\cr
  \rho_{5}&4&0&\frac{-1-3\sqrt{-3}}{2}&\frac{-1+3\sqrt{-3}}{2}&-2&1&\zeta_3&\zeta_3^2&0&\zeta_3&1&\zeta_3^2&1\cr
  \rho_{6}&4&0&\frac{-1+3\sqrt{-3}}{2}&\frac{-1-3\sqrt{-3}}{2}&-2&1&\zeta_3^2&\zeta_3&0&\zeta_3^2&1&\zeta_3&1\cr
  \rho_{7}&4&0&\frac{-1+3\sqrt{-3}}{2}&\frac{-1-3\sqrt{-3}}{2}&-2&1&\zeta_3&\zeta_3^2&0&1&\zeta_3&1&\zeta_3^2\cr
  \rho_{8}&4&0&\frac{-1-3\sqrt{-3}}{2}&\frac{-1+3\sqrt{-3}}{2}&-2&1&1&1&0&\zeta_3^2&\zeta_3&\zeta_3&\zeta_3^2\cr
  \rho_{9}&4&0&\frac{-1+3\sqrt{-3}}{2}&\frac{-1-3\sqrt{-3}}{2}&-2&1&1&1&0&\zeta_3&\zeta_3^2&\zeta_3^2&\zeta_3\cr
  \rho_{10}&4&0&\frac{-1-3\sqrt{-3}}{2}&\frac{-1+3\sqrt{-3}}{2}&-2&1&\zeta_3^2&\zeta_3&0&1&\zeta_3^2&1&\zeta_3\cr
  \rho_{11}&6&2&-3&-3&3&0&0&0&-1&0&0&0&0\cr
  \rho_{12}&6&-2&-3&-3&3&0&0&0&1&0&0&0&0\cr
  \rho_{13}&12&0&3&3&0&-3&0&0&0&0&0&0&0\cr
\end{array}
$$

\vskip1cm

\centerline{Table 7. Irreducible characters of $C_3^4{\rtimes}A_6$ of degree 30}

\med
\tiny{
$$
\begin{array}{|c|ccrrrrrcrrrrrrccc|}\hline
  \rm class&\rm3A&\rm3B&\rm3C&\rm2A&\rm6A&\rm6B&\rm6C&\rm3D&\rm3E&
  \rm3F&\rm3G&\rm3H&\rm3I&\rm9A&\rm{9BCD}&\rm4A&\rm{5AB}\cr\hline

\hline
\rho_{10}&3&-6&3&2&2&-1&-1&0&6&0&-3&0&-3&0&0&0&0\cr
\rho_{11}&3&3&-6&2&-1&2&-1&6&0&-3&0&-3&0&0&0 &0&0\cr

\rho_{18}&3&3&-6&2&-1&2&-1&-3&0&-3&0&6&0&0&0 &0&0 \cr
\rho_{19}&3&3&-6&2&-1&2&-1&-3&0&6&0&-3&0&0&0&0&0 \cr

\rho_{20} &3&-6&3&2&2&-1&-1&0&-3&0&-3&0&6&0&0 &0&0 \cr
\rho_{21} &3&-6&3&2&2&-1&-1&0&-3&0&6&0&-3&0& 0&0&0 \cr

\hline
\end{array}
$$

\vskip1cm


}


\begin{thebibliography}{hhhh}
\bibitem{AGJ} J. Aaronson, C. Groengl,  T. Johnston, Cyclically covering subspaces in $\FF_2^n$, J. Combin. Theory, Ser. A 181 (2021), Paper No. 105436, 41 pp.

\bibitem {Asch} M. Aschbacher, {\it Finite group theory}, 2nd edition, Cambridge Univ. Press,  Cambridge, 2000.




\bibitem {Be} \'A. Bereczky, Fixed-point-free $p$-elements in transitive permutation groups, Bull. London Math. Soc. 27(1995), $447-452$.


\bibitem{Bo} N. Bourbaki, {\it Groupes et algebres de Lie}, ch. IV-VI, Masson, Paris, 1981.

\bibitem{Bo8} N. Bourbaki, {\it Groupes et algebres de Lie}, ch. VII-VIII, Springer, Berlin, 2006.


\bibitem{Lu2}   T. Breuer, Decomposition matrices, \\   http://ftp.math.rwth-aachen.de/homes/MOC/decomposition/tex/He/He.2mod2.pdf


\bibitem {Bu} R. Burkhardt, Die Zerlegungsmatrizen der Gruppen $PSL(2,p^f)$, {\it J. Algebra} 40(1976), $75-96$.

\bibitem{Ca} P. Cameron, D. Ellis and W. Raynaud, Smallest cyclically covering subspaces of $F^n_q$, and lower bounds in Isbell's conjecture, {\it Europ. J. Combinatorics} 81 (2019), 242–255. 

\bibitem{cfk} P. Cameron, P. Frankl, W. Kantor, Intersecting families of finite sets and fixed-point-free 2-elements, {\it European J. Combin.} 10 (1989), $ 149-160$.

\bibitem{C}  R. Carter, `{\it Finite Groups of Lie type: Conjugacy Classes and Complex Characters}', Wiley, Chichester, $1985$.

\bibitem{Atlas} J.H. Conway, R.T. Curtis, S.P. Norton,  R.A. Parker and R.A. Wilson, {\it An ATLAS of Finite Groups}, Clarendon Press, Oxford, $1985$.

\bibitem{cu9} J. Cullivan, A computational approach to the $2$-torsion of abelian threefolds, {\it Mathematics of Computations} 78(2009), $1825-1836$.

\bibitem {cu10} J. Cullivan, Points of small order on three-dimensional abelian varieties; with an appendix by Yuri Zarhin, {\it J. Algebra} 324 (2010), $565-577.$

\bibitem{cu12}   J. Cullinan, Symplectic stabilizers with applications to abelian varieties, {\it Intern.  J.  Number Theory}  8, No. 2 (2012), $321-334$.

\bibitem {gl3} J.~Cullinan, Fixed-point subgroups of $\mathrm{GL}_3(q)$, {\it J. Group Theory} \textbf{22}(2019), $893 - 914$.

\bibitem {CZ} J. Cullivan and A. Zalesski, Unisingular Representations in Arithmetic and Lie Theory,  {\it Europ. J. Math.}  7(2021), no.4, $1645 - 1667$.

\bibitem{CR} Ch. Curtis and I. Reiner, {\it Representation theory of finite groups and associative algebras}, Interscience, New York, 1962.

\bibitem{CR1} Ch. Curtis and I. Reiner, {\it Methods of representation theory with applictions to  finite groups and orders}, vol.I, , Interscience, New York, 1981.

\bibitem{DF}  D.I. Deziriotis and A.P. Fakiolas, Maximal tori in the finite Chevalley groups of types $E_6$, $E_7$ and $E_8$,  {\it Comm.  Algebra} 19(1991), $889 - 903$.


\bibitem{DZ1}  L. Di Martino and A. Zalesskii, Minimum polynomials and lower bounds for eigenvalue
multiplicities of prime-power order elements in representations of  quasi-simple groups,
{\it J. Algebra} {\bf 243} $(2001)$, $228 - 263$.
Corrigendum:   {\it J. Algebra} {\bf 296} $(2006)$, $249 - 252$.

\bibitem{DZ8} L. Di Martino and A. Zalesski, Eigenvalues of unipotent elements
in cross-characteristic representations of finite classical groups.
{\it J. Algebra}  319(2008),   $2668 - 2722$.

\bibitem{DZ2}  L. Di Martino and A. Zalesskii, Almost cyclic elements in Weil \reps of finite classcal groups, {\it Comm. Algebra} 46(2018), $2767 - 2810.$

\bibitem{Di} J. Dixon and B. Mortimer, {\it Permutation groups}, Springer-Verlag, Berlin, $1996$.

\bibitem{EZ} L. Emmett and A.E. Zalesski, On regular orbits of elements of classical groups in their permutation representations, {\it Comm. Algebra} 39(2011), 3356 -- 3409.

\bibitem{Fe}  W. Feit, {\it The Representation Theory of Finite Groups}, North-Holland, Amsterdam, $1982$.


\bibitem{Ge}   P.  Gerardin, Weil \reps associated to finite fields,
{\it J. Algebra} 46(1977), $54 - 101$.


\bibitem{s31} Group Names, an online source of properties of finite groups of order at most $500$,  \begin{verbatim} https://people.maths.bris.ac.uk/~matyd/GroupNames/321/C3%5E3sA4.html\end{verbatim}

\bibitem{rr1} Group Names, an online source of properties of finite groups of order at most $500$, Transitive groups of degree up to 31,
~\begin{verbatim} https://people.maths.bris.ac.uk/~matyd/GroupNames/T31.html
\end{verbatim}

\bibitem{GT} R. Guralnick and Pham Huu Tiep, Finite simple unisingular groups of Lie type, {\it J. Group Theory} 6(2003), $271 - 310$.


\bibitem{HeB} D.R. Heath-Brown, Artin's conjecture for primitive roots, {\it Quart. J. Math. Oxford} 37 (1986), $27-38$.


\bibitem{HM1}  G. Hiss and G. Malle,  Low-dimensional representations of quasi-simple groups, {\it LMS J. Comp. Math}. 4(2001), $22 -63.$

\bibitem{HM}  G. Hiss and G. Malle, Corrigenda: Low-dimensional \reps of quasi-simple groups, {\it LMS J. Comp. Math}. 5(2002), $95 -126.$






\bibitem{HB}  B. Huppert and N. Blackburn, ``{\it Finite Groups} III", Springer-Verlag, Berlin etc., 1982.

\bibitem{Is} J. Isbell, Homogeneous games, {\it Math. Stud.} 25(1957), $123-128$.

\bibitem{Is2} J. Isbell, Homogeneous games, II, {\it Proc. Amer. Math. Soc}. 11(1960), $159-161$.


\bibitem{JLPW}  C. Jansen, K. Lux, R.A. Parker, and R.A. Wilson, `{\it An ATLAS
 of Brauer Characters}', Oxford University Press, Oxford, $1995$.

\bibitem{katz} N.M.~Katz,  Galois properties of torsion points on abelian varieties, {\it Inventiones Math.} 62(1981), $481-502.$

\bibitem{KL}  P. B. Kleidman and M. W. Liebeck, `{\it The Subgroup Structure of the Finite Classical Groups}', London Math. Soc. Lecture Note Ser. no.
$129$, Cambridge University Press, $1990$.

\bibitem{zk} A.S. Kleshchev and  A.E. Zalesski, Minimal \pos of elements of order $p$
in $p$-modular projective \reps of alternating groups,   
{\it Proc. Amer.  Math. Soc. } 132(2003),   $1605-1612$.

\bibitem{LS} M. Liebeck and G. Seitz, Reductive subgroups of exceptional algebraic groups, {\it Memoirs  Amer. Math. Soc.} 121(1996), no.580, $1 -111$.


\bibitem{Lu}  F. L\"ubeck, Small degree \reps of finite Chevalley groups in defining characteristic,  {\it  LMS J. Comp. Math}. 4(2001), $135 - 169$.


\bibitem{Luh} J. Luh, On the representation of vector spaces as a finite union of subspaces, {\it Acta Math. Acad. Sci. Hungar}. 23 (1972), $341-342.$

\bibitem{MR} G. Malle and G. Robinson, Projective indecomposable permutation modules, {\it Vietnam J. Math.}   51(2023), $633 - 656$. 

 

\bibitem{Ne} B.H. Neumann, Groups covered by permutable subsets, {\it J. London Math. Soc.} 29 (1954), $236-248.$

\bibitem{OL1} The Online Encyclopedia of   integer sequences, https://oeis.org/A062117

\bibitem{OL2} The Online Encyclopedia of   integer sequences, 
https://oeis.org/A211241

\bibitem{sp} P. Spiga, $p$-{\it elements in permutation groups}, PhD thesis, Queen Mary College, University of London, 2004.

\bibitem{St} R. Steinberg, {\it Lectures on Chevalley groups},  Amer. Math. Soc. Univ. Lect. Series, vol. 66, Providence, Rhode Island, 2016.


\bibitem{TeZ} D. Testerman and A.E. Zalesski, Almost cyclic semisimple  regular elements in irreducible representations of simple algebraic groups, {\it Transformation Groups} 28(2023), $1299-1324$. 

\bibitem{TZ8}  Pham Huu Tiep and A.E. Zalesski, Hall-Higman type theorems for semisimple elements of finite classical groups. {\it Proc. London Math. Soc}. (3) 97(2008), 623 -- 668.

\bibitem{TZ22}   Pham Huu Tiep and A.E. Zalesski, Hall-Higman type theorems for exceptional groups of Lie type, I, {\it J. Algebra}  607(2022), Part A, $755-794$.

\bibitem{TZ23}   Pham Huu Tiep and A.E. Zalesski, Hall-Higman type theorems for simple groups with cyclic Sylow $p$-subgroups, (in preparation) 


\bibitem{W} A. Wagner, The faithful linear \reps of least degree of $S_n$ 
and $A_n$ over a field of characteristic $2$, {\it Math. Zeitschr.} 151(1976), $127-137$.  




\bibitem{Z71}  A.E. Zalesski, Subgroups of classical groups, {\it Siber. Math. J.} 12(1971), $90 - 94$.

\bibitem{Z86}A.E. Zalesski,  Spectra of elements of order $p$ in representations of  Chevalley      groups of characteristic $p$ (in Russian), {\it Vesti  AN  BSSR, ser. fiz.-mat. navuk}, 1986, no.6, 20 - 25.

\bibitem{Z87} A.E. Zalesski, Fixed points of elements of order $p$ in complex \reps of finite Chevalley groups in characteristic $p$ (in Russian), {\it Doklady Akad. Nauk, Belorussian SSR} 31(1987), $104 - 107$.

\bibitem{Z88}  A.E. Zalesski, Eigenvalues of matrices of complex representations of finite groups of Lie type. In: ``{\it Algebra, some current trends}, Lecture Notes in Math. (Springer-Verlag), vol. 1352, 1988'', 206 - 218.



\bibitem{Z90}  A.E. Zalesski, The eigenvalue $1$ of matrices of complex representations of finite Chevalley groups, {\it Proc. Steklov Inst. Math.} 1991, issue 4, $109 - 119.$

\bibitem{z91} A.E. Zalesski\u\i, A fragment of the decomposition matrix of the  special  unitary   group over a  finite field. English translation: {\it Math. USSR, Izvestija} 36(1991),      $23 -  39$. (Russian original: {\it Izvestia Akad. Nauk SSSR, ser. matem.}
      54(1990), $26-41$.)



\bibitem{z08}  A.E. Zalesski, Minimal polynomials of the elements of prime order in complex irreducible representations of quasi-simple groups. {\it J. Algebra} 320(2008), $2496 - 2525$.


 \bibitem{z16}  A.E. Zalesski, Invariants of maximal tori and unipotent constituents of some quasi-projective characters for finite classical groups, {\it J. Algebra} 500(2018), $517 - 541$.

\bibitem{z18}  A.E. Zalesski, Singer cycles in $2$-modular \reps of $GL_n(2)$,  {\it Archiv der Math.} 110(2018), $433 - 446$.


\bibitem{z21}  A.E. Zalesski, Unisingular representations of finite symplectic groups,  {\it Comm. Algebra} 50(2022), $1697 - 1719$.

\bibitem{z22}  A.E. Zalesski, Upper bounds for eigenvalue multiplicities of almost cyclic elements in \ir \reps of simple algebraic groups, {\it European J. Math.} 9(2023),
no.3,  doi: 10.1007/s40879-023-00680-7 

\bibitem{z23}  A.E. Zalesski, Rational elements in representations of simple
algebraic groups, I, {\it Vietnam J. Math.} 52(2024), https://doi.org/10.1007/s10013-023-00629-z


 \end{thebibliography}
\end{document}